\documentclass[a4paper,12pt]{book}

\usepackage{amsmath,amsfonts,amsthm,amssymb,amscd,enumerate,graphics,graphicx}
\usepackage[all]{xy}

\usepackage{amentashortcuts}
\usepackage{amentastyle}
\usepackage[draft=false, bookmarks, colorlinks=true, pdftitle={The Geometry of Orbifolds via Lie Groupoids}, pdfauthor={Alexander Amenta}]{hyperref}

\newcommand{\blanknonumber}{\newpage\thispagestyle{empty}}

\begin{document}
\setcounter{chapter}{0}
\frontmatter
\begin{titlepage}
\begin{center}

  \vspace*{\fill} \Huge
The Geometry of Orbifolds \\ via Lie Groupoids
\\
\vfill\vfill\Large
                          Alexander Amenta
\\
\vfill\vfill
                          May 2012
\\
\vfill\vfill \normalsize
A thesis submitted in partial fulfilment of the requirements of the degree of \\
Bachelor of Science (Honours) in Mathematics\\
of the Australian National University
\vfill
         \includegraphics{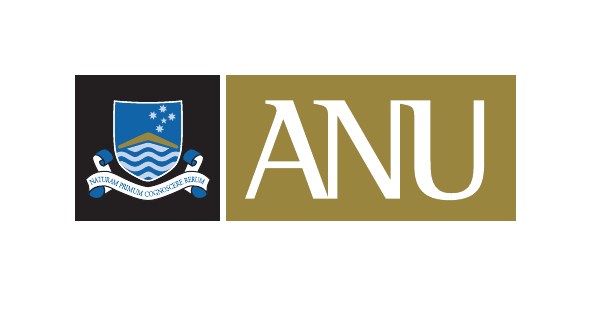}

\end{center}

\end{titlepage}

\blanknonumber\ \blanknonumber

\vspace*{\fill}

\begin{quote}
  Ithaca has given you the beautiful journey. \\
  Without her you would not have set out on the road. \\
  Nothing more does she have to give you. \\
\end{quote}

\begin{flushright}
  -Constantine P. Cavafy, \emph{Ithaca}
\end{flushright}

\vspace{5cm}

\begin{quote}
  ``NO TIME FOR THAT! WE'RE LOOKING AT POLYHEDRA!!!!''
\end{quote}

\begin{flushright}
  -POKEY THE PENGUIN, \\ \emph{ONE DAY AT THE MUSEUM OF GEOMETRY}
\end{flushright}

\vspace*{\fill}
\blanknonumber
\chapter*{Acknowledgements}\label{acknowledgements}
\addcontentsline{toc}{chapter}{Acknowledgements}

First of all, I would like to thank Dr Bryan Wang for two years of mathematical supervision, starting with a seminar course on Dirac operators and culminating in this thesis.
His guidance and enthusiasm have well and truly cemented my interest in mathematics as a whole, and I am constantly impressed by his tolerance of my increasing stubbornness.
Similarly I thank all the other mathematicians who have guided me along the way: there are too many to name (and indeed I have probably forgotten some), but they know who they are.

Secondly I thank my friends: again, too many to name (and again I have probably forgotten some), but likewise they know who they are.
In particular I thank Ben and Chian (for their almost painful enthusiasm), Chris (for his inhuman ability to live on rice and his enthusiasm for Twitter), Jon (for being Jon), Lachie (for being a top bloke and a below-average bassist), and Ping (for her constant support and for putting up with my rants about orbifolds).

I also thank my small proofreading team---all of whom have already been mentioned---for being my comrades in the unwinnable war against typos.

Finally I thank my father Adam-Joseph, for teaching me the value of fighting the system.

\blanknonumber
\footnotesize
\tableofcontents\blanknonumber
\normalsize
\mainmatter
\chapter*{Introduction}
\addcontentsline{toc}{chapter}{Introduction}

Orbifolds were first defined by Ichir\^{o} Satake in the mid-fifties under the name `$V$-manifolds' \cite{iS56,iS57}.\footnote{Conflicting accounts on the meaning of the letter `V' have been given. Some suggest it stands for `Virtual', virtual meaning `up to a finite cover'. Others suggest that the `V' is a typographic representation of a cone point.}
In his initial work, Satake proved $V$-manifold analogues of some classical theorems in differential geometry, including the de Rham theorem, Poincar\'e duality, and the Gauss-Bonnet theorem.
In doing this, Satake showed that $V$-manifolds are quite similar to smooth manifolds.
Over the next twenty years, $V$-manifolds were considered to be a straightforward generalisation of smooth manifolds \cite[p. ix]{ALR07}.
Although they were useful geometric tools (indeed Satake introduced $V$-manifolds in order to study the Siegel modular variety), they were not considered as a legitimately `new' concept in this time.

In the late seventies and early eighties, Tetsuro Kawasaki wrote a small series of articles presenting $V$-manifold versions of the foundational theorems of index theory: the Hirzebruch signature theorem \cite{tK78}, the Hirzebruch-Riemann-Roch theorem \cite{tK79}, and the Atiyah-Singer index theorem \cite{tK81}.
Although these results are generalised versions of well-known results concerning smooth manifolds, the proofs of these theorems emphasised features of the $V$-manifold theory which are not present for that of smooth manifolds.
In particular, Kawasaki introduced what is now called the \emph{inertia orbifold} associated to a given $V$-manifold.
When applied to a smooth manifold, this construction simply returns the original manifold; this construction is only significant when applied to nontrivial $V$-manifolds.
Furthermore, this construction highlights the need to consider \emph{ineffective} orbifolds, which were not discussed in Satake's work.
Thus it was realised that $V$-manifolds have a geometric theory which is not merely a generalisation of that of smooth manifolds.

Around the same time, William Thurston made use of orbifolds in his study of the geometry of three-manifolds.
Given that $V$-manifolds deserved to be considered as being distinct from smooth manifolds, the search for a new name was initiated.
Eventually the name `orbifold' was settled upon, and over time this name became standard.\footnote{The name `orbifold' is meant to evoke imagery of a manifold-like object which can be obtained as an orbit space.}

In the mid-eighties, orbifolds were being used by physicists in their study of conformal field theory \cite{DHVW85,DFMS87}.
This work was important not only for its physical significance, but also for the physical intuition it brought to the purely mathematical theory of orbifolds.
A more informed overview of this aspect of orbifold history can be found in \cite{ALR07}.

In this thesis we are concerned with the methods used to represent orbifolds as geometric objects.
Satake's work was done in terms of \emph{orbifold atlases}; this approach is geometrically intuitive, but does not readily tell us much about the global structure of orbifolds.
More recently orbifolds have been represented by \emph{orbifold groupoids}; this is the approach we will ultimately follow.
This approach can be traced back to Andr\'e Haefliger's work on foliations \cite{aH84}, and was more recently developed in the context of orbifolds by Ieke Moerdijk and Dorotea Pronk \cite{iMdP97}.
Unlike orbifold atlases, orbifold groupoids can be used to study the global structure of orbifolds.
However, they have no particular advantage with respect to local matters.
The purpose of this thesis is to use the language of orbifold groupoids to describe the geometry and topology of orbifolds, highlighting advantages and disadvantages of this language as they arise.


\chapter{Orbifolds: the classical viewpoint}\label{c1}

In this chapter, we shall introduce orbifolds as they were originally studied: as topological spaces enriched by `orbifold atlases'.
This is the most intuitive way of dealing with orbifolds, but it is not always the most convenient, and has some interesting flaws when used naively.
We shall introduce the concepts of smooth maps between and vector bundles over orbifolds, highlighting the aforementioned flaws in the process.
Finally, we will discuss ineffective orbifolds, further motivating the shift towards modern methods (both atlas-theoretic and otherwise).

The exposition in this chapter is essentially a combination of Satake's original exposition in \cite{iS56,iS57} and that of the recent text \cite{ALR07}.
We have adopted a new choice of notation which, despite its reliance on multiple typefaces, makes for simpler definitions and clearer exposition than these references.

\section{Effective orbifolds}\label{effectives}

Put briefly, an $n$-dimensional orbifold is a space which is locally given by the orbit space $G \sm \wtd{U}$ of the action of a finite group $G$ on a connected open subset $\wtd{U} \subset \RR^n$ of Euclidean $n$-space.
However, it is not enough to say that an orbifold is a topological space which is locally homeomorphic to such an orbit space, as we want to keep track of these `local actions'.
Thus we make the following definition.

\begin{dfn}\label{chart}
  Let $X$ be a topological space.
  An $n$-dimensional \emph{orbifold chart}\footnote{Orbifold charts were originally called \emph{local uniformising systems}, as in \cite{iS56} and \cite{iS57}. This terminology still appears in current literature.}  $\mathbf{U} = (\wtd{U},G_U,\gf_U)$ on $X$ is given by
  \begin{itemize}
  \item
    a connected open subset $\wtd{U} \subset \RR^n$ of Euclidean $n$-space,
  \item
    a finite group $G_U$ along with a smooth effective left action of $G_U$ on $\wtd{U}$, and
  \item
    a continuous $G_U$-map $\map{\gf_U}{\wtd{U}}{X}$ which induces a homeomorphism of the orbit space $G_U \sm \wtd{U}$ onto its image $U:= \gf_U(\wtd{U}) \subset X$.
  \end{itemize}
  We will often refer to an orbifold chart more concisely as a \emph{chart}.
\end{dfn}

Note that the action associated to an orbifold chart is assumed to be effective; that is, the only element of $G_U$ which acts trivially on $\wtd{U}$ is the identity element.

\begin{rmk}
  We have used different typefaces to refer to different aspects of an orbifold chart.
  A chart $\mb{U}$ is written in boldface, the corresponding open subset of $\RR^n$ is written $\wtd{U}$ (lightface with tilde), and the image $U := \gf_U(\wtd{U})$ of this space under the map $\gf_U$ is written in lightface with no tilde.
  The group $G_U$ and map $\gf_U$ associated to this chart are indexed by the `name' of the chart (written in lightface).
  So a chart $\mb{V}$ would be a triple $\mb{V} = (\wtd{V},G_V,\gf_V)$, and we would write $V = \gf_V(\wtd{V})$.
  We will occasionally break this convention (indeed, we do so in the following example), but it should be clear from context what we mean.
\end{rmk}

\begin{figure}[t]
  \centering
  \includegraphics[keepaspectratio=true,scale=1.8]{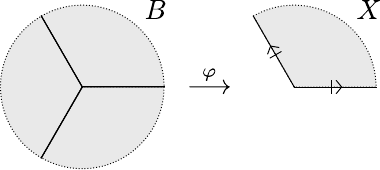}
  \caption{The orbifold chart of Example \ref{disc} (with $n=3$).}
  \label{cone1}
\end{figure}

\begin{example}\label{disc}
For a simple example of an orbifold chart, let $B$ denote the open unit ball in $\RR^2$, and let the cyclic group $\ZZ_n$ act on $B$ by rotation.
More explicitly, identify $\RR^2$ with the complex numbers $\CC$, identify $\ZZ_n$ with the group of $n^{\text{th}}$ roots of unity
\begin{equation*}
  \ZZ_n = \langle \gz_n \rangle < \CC, \quad \gz_n := e^{\frac{2\gp i}{n}},
\end{equation*}
and let $\ZZ_n$ act on the complex unit ball $B \subset \CC$ by multiplication.
Set $X$ to be the sector of $B$ with angle $2\gp/n$ with edges identified as in Figure \ref{cone1}, and let $\map{\gf}{B}{X}$ be the natural $n$-fold cover of $X$.
Since $\gf$ induces a homeomorphism from the orbit space $\ZZ_n \sm B$ to $X$, we find that $(B,\ZZ_n,\gf)$ is a $2$-dimensional orbifold chart on $X$.
\end{example}

Thanks to our newly-expanded vocabulary, we can be more specific: an orbifold is a topological space which is covered by orbifold charts.
In this way, the notion of an orbifold is a generalisation of the notion of a smooth manifold.
We immediately reach another problem: if we have two orbifold charts $\mathbf{U}, \mathbf{V}$ on a space $X$ which overlap, in the sense that $U \cap V$ is nonempty, how do we express `local compatibility' of the two charts?
Invoking the analogy of an atlas on a smooth manifold, if a point $x \in X$ lies in the intersection of $U$ and $V$, there should be a third chart $\mathbf{W}$ with $x \in W$ which somehow `sits inside' the other two charts.
This consideration leads us to the next definition.

\begin{dfn}
  Let $\mathbf{U}$ and $\mathbf{V}$ be two orbifold charts on a topological space $X$.
  An \emph{embedding}\footnote{Embeddings were originally called \emph{inclusions} in \cite{iS56} and \cite{iS57}.} $\map{\gl}{\mathbf{U}}{\mathbf{V}}$ is given by an embedding (in the usual sense) $\map{\gl}{\wtd{U}}{\wtd{V}}$ such that $\gf_V \circ \gl = \gf_U$.
  Note that embeddings can be composed: if $\map{\gl}{\mb{U}}{\mb{V}}$ and $\map{\gm}{\mb{V}}{\mb{W}}$ are embeddings, then the composition $\map{\gm \circ \gl}{\wtd{U}}{\wtd{W}}$ determines an embedding $\map{\gm \circ \gl}{\mb{U}}{\mb{W}}$.
\end{dfn}

We are now well on our way to pinpointing what an orbifold actually is.
Just as a smooth structure on a topological manifold is given by an equivalence class of atlases, an orbifold structure on a topological space is given by an equivalence class of `orbifold atlases', as follows.

\begin{dfn}
  \begin{enumerate}[(1)]
  \item
    An $n$-dimensional \emph{orbifold atlas} $\mc{U}$ on a topological space $X$ is a set of $n$-dimensional orbifold charts such that
    \begin{itemize}
    \item
      the charts in $\mc{U}$ cover $X$, in the sense that
      \begin{equation*}
        \bigcup_{\mathbf{U} \in \mc{U}} U = X, \quad \text{and}
      \end{equation*}
    \item
      if $\mathbf{U}$ and $\mathbf{V}$ are charts in $\mc{U}$, then for all $x \in U \cap V$ there exists a chart $\mathbf{W}$ in $\mc{U}$ with $x \in W \subset U \cap V$ such that there exist embeddings of $\mathbf{W}$ into both $\mathbf{U}$ and $\mathbf{V}$.
    \end{itemize}
  \item
    A \emph{refinement} $\mc{U} \to \mc{V}$ of orbifold atlases on $X$ is given by a function $\map{r}{\mc{U}}{\mc{V}}$ and, for each chart $\mathbf{U} \in \mc{U}$, an embedding $\map{\gl_{\mathbf{U}}}{\mathbf{U}}{r(\mathbf{U})}$.
    In this situation we say that $\mc{U}$ is a refinement of $\mc{V}$.
    We define a symmetric and reflexive relation $\sim$ on the set of orbifold atlases on $X$ by taking $\mc{U} \sim \mc{V}$ to mean that there exists a common refinement of $\mc{U}$ and $\mc{V}$.\footnote{Some authors refer to this relation as \emph{direct equivalence}.}
    We say that two orbifold atlases are \emph{equivalent} if they are related by the transitive closure of $\sim$.
  \item
    An \emph{orbifold structure} on a topological space $X$ is given by an equivalence class $[\mc{U}]$ of orbifold atlases on $X$.
  \end{enumerate}
\end{dfn}

\begin{rmk}
  Equivalence of orbifold atlases is our first technical hurdle.
  Unlike equivalence of manifold atlases, it is not at all obvious that the relation $\sim$ as defined above is transitive.
  Originally this issue was avoided by taking the transitive closure as we have \cite[p467]{iS57}.
  However, in \cite[Lemma 1.34]{mT10} it is shown that $\sim$ is transitive by using an equivalent definition of the relation.
  The proof is too much of a diversion to be included in this work, so we shall be content with taking the transitive closure.
\end{rmk}

\begin{dfn}\label{eff}
  An \emph{effective orbifold} $\mf{X} = (X,[\mc{U}])$ is given by a paracompact Hausdorff space $X$ along with an orbifold structure $[\mc{U}]$ on $X$.
  We call $X$ the \emph{underlying space} of $\mf{X}$, and denote it by $|\mf{X}|$.
  The effective orbifold $\mf{X}$ is said to be \emph{compact} (resp. \emph{connected}) if the space $X$ is compact (resp. connected).
\end{dfn}

\begin{rmk}\label{ineff}
  An \emph{ineffective orbifold} is one in which the groups $G_U$ in charts $\mb{U} \in \mc{U}$ are not assumed to act effectively on $\mc{U}$.
  Defining ineffective orbifolds via orbifold atlases causes (solvable) problems; we consider this in Section \ref{ineffective}.
  For now we shall use the terms `effective orbifold' and `orbifold' interchangeably.
\end{rmk}

\begin{figure}[t]
  \centering
  \includegraphics[keepaspectratio=true,scale=1.8]{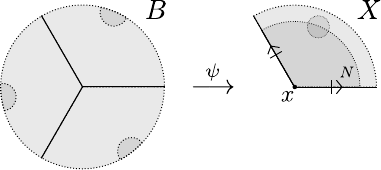}
  \caption{An orbifold chart $(B,\ZZ_3,\gy)$ on $X$, and embeddings $V_\ga \cap N \to B$. ($V_\ga$ is the small open ball shown intersecting $N$.)}
  \label{cone2}
\end{figure}

In Example \ref{disc} we saw an example of an orbifold chart $(B,\ZZ_n,\gf)$ on a paracompact Hausdorff space ${X \cong \ZZ_n \sm B}$.
This chart defines an orbifold atlas $\mc{U}$ on $X$ consisting only of the chart itself; hence $\mf{X} := (X,[\mc{U}])$ is an orbifold.
We shall construct an equivalent atlas for $\mf{X}$.
Let $x$ denote the point $\gf(0) \in X$.
The space $X \sm \{x\}$ is homeomorphic to a cylinder, and thus inherits a smooth manifold structure, so there exists a smooth manifold atlas for $X \sm \{x\}$.
Let $V = \{(V_\ga,f_\ga)\}_{\ga \in A}$ denote such an atlas (with coordinate maps $\map{f_\ga}{V_\ga}{\RR^2}$ and indexing set $A$).
Then letting $\mathbf{1}$ denote the trivial group, $(f_\ga(V_\ga),\mathbf{1},f_\ga^{-1})$ is an orbifold chart on $X$ for each $\ga$.
Let $N$ be a neighbourhood of $x$ in $X$ obtained by `scaling down' $X$, as in Figure \ref{cone2}, and let $\map{\gy}{B}{N}$ be the obvious modification of the map $\map{\gf}{B}{X}$.
Then $(B,\ZZ_n,\gy)$ is an orbifold chart on $X$ with $\gy(B) = N$.
For each $\ga \in A$ such that $V_\ga \cap N$ is nonempty, $(f_\ga(V_\ga \cap N),\mb{1},f_\ga^{-1})$ is an orbifold chart on $X$.
In order to construct an embedding from the chart $(f_\ga(V_\ga \cap N),\mb{1},f_\ga^{-1})$ to the chart $(B,\ZZ_n,\gy)$, it is sufficient to give an embedding $\map{\gl}{V_\ga \cap N}{B}$ such that $\gy \circ \gl = \id$.
It is easy to see pictorially that there are $n$ such embeddings (this is shown in Figure \ref{cone2}).
Therefore the union $\mc{V}$ of the above charts defines an orbifold atlas on $X$.
The atlas $\mc{V}$ is readily seen to be a refinement of $\mc{U}$, and so $\mc{V}$ and $\mc{U}$ define the same orbifold structure on $X$.

With the previous example in mind, we begin to see that each point $x \in X$ of an orbifold $(X,[\mc{U}])$ comes equipped with intrinsic algebraic data, namely the \emph{local group} $G_x$ of the point $x$.
In this example, the `cone point' $x \in X$ has local group $G_x \cong \ZZ_n$, while all points $y \in X \sm \{x\}$ have trivial local group.
This information distinguishes the orbifold $\mf{X}$ from the orbit space $\ZZ_n \sm B$.
We shall now make these intuitive notions rigourous.

\begin{lem}\label{satake1}
  Let $\mb{U}$ and $\mb{V}$ be two orbifold charts on a topological space $X$.
  Then every embedding $\map{\gl}{\mb{U}}{\mb{V}}$ induces an injective group homomorphism $\map{\gl_*}{G_U}{G_V}$.
\end{lem}

\begin{proofsketch}
  First note that for every group element $g \in G_U$, left multiplication by $g$ defines an embedding $\map{g}{\mb{U}}{\mb{U}}$.
  To complete the proof one goes on to show that for each $g \in G_U$ there is a unique $g^\prime \in G_V$ such that
  \begin{equation}\label{indhom}
    \gl \circ g = g^\prime \circ \gl.
  \end{equation}
  Setting $\gl_*(g) := g^\prime$ then defines the desired homomorphism.
  For a full proof, see \cite[Lemma 1]{iS57} for the case where the fixed point set $\wtd{U}^G$ has codimension $\leq n-2$ (where $n$ is the dimension of the chart $\mb{U}$)\footnote{In \cite{iS57}, this codimension condition is part of the definition of an orbifold chart/local uniformising system.}, or \cite[Proposition A.1]{iMdP97} for the general case.
\end{proofsketch}

The following lemma characterises the image of this homomorphism, and is needed to define local groups.

\begin{lem}\label{image}
  Let $\map{\gl}{\mb{U}}{\mb{V}}$ be an embedding as above, and suppose $h \in G_V$.
  Then $h \in \im \gl_*$ if and only if $\gl(\wtd{U}) \cap h\circ\gl(\wtd{U}) \neq \varnothing$, in which case $\gl(\wtd{U}) = h \circ \gl(\wtd{U})$.
\end{lem}

\begin{proof}
  (Adapted from \cite[Lemma 2]{iS57})
  If $h$ is in the image of $\gl_*$, then by \eqref{indhom} we have
  \begin{equation*}
    h \circ \gl(\wtd{U}) = \gl \circ \gl_*^{-1}(h)(\wtd{U}) = \gl(\wtd{U}),
  \end{equation*}
  from which the claim follows.
  Conversely, suppose that $h \notin \im \gl_*$ and $\gl(\wtd{U}) \cap h\circ\gl(\wtd{U})$ is nonempty, and let $\gl(x) \in \gl(\wtd{U}) \cap h\circ\gl(\wtd{U})$ be a point which is not fixed by $G_V$.
  Such a point exists provided $G_V$ is not trivial\footnote{If $G_V$ is trivial, then there is nothing to show.} since the fixed point set $\wtd{V}^{G_V}$ is nowhere dense in $\wtd{V}$ (by Corollary \ref{nowheredense}), and since $\gl(\wtd{U}) \cap h\circ\gl(\wtd{U})$ is open and by assumption nonempty in $\wtd{V}$.
  Since $\gl(x)$ is in $h \circ \gl(\wtd{U})$, there is a point $y \in \wtd{U}$ such that $\gl(x) = h\circ\gl(y)$.
  Thus
  \begin{equation*}
    \gf_V(\gl(x)) = \gf_V(h \circ \gl(y)) = \gf_V(\gl(y)),
  \end{equation*}
  and since $\gl$ is an embedding, this implies that $\gf_U(x) = \gf_U(y)$.
  Hence there exists $g \in G_U$ such that $y=gx$.
  Then
  \begin{equation*}
    \gl(x) = h \circ \gl(gx) = h\gl_*(g)\gl(x)
  \end{equation*}
  by \eqref{indhom}, and since $\gl(x)$ is not fixed by $G_V$, we have that $h\gl_*(g)$ is the identity in $G_V$.
  But then $h = \gl_*(g)^{-1} = \gl_*(g^{-1})$ is in the image of $\gl_*$, contradicting the assumption that $\gl(\wtd{U}) \cap h\circ\gl(\wtd{U})$ is nonempty.
\end{proof}

\begin{cor}
  Let $\mb{U}$ and $\mb{V}$ be two orbifold charts on a space $X$, let $\map{\gl}{\mb{U}}{\mb{V}}$ be an embedding, and suppose $x \in U \cap V$.
  Then for any choice of $x_U \in U$ and $x_V \in V$ lying above $x$ (in the sense that $\gf_U(x_U) = \gf_V(x_V) = x$), the isotropy groups $(G_U)_{x_U}$ and $(G_V)_{x_V}$ are isomorphic.
\end{cor}

\begin{proof}
  Choose $x_U \in U$ lying above $x$.
  By Lemma \ref{image}, $(G_V)_{\gl(x_U)}$ is in the image of $\gl_*$.
  Using the defining property \eqref{indhom} of $\gl_*$, it is then easy to show that $\gl_*$ maps $(G_U)_{x_U}$ isomorphically onto $(G_V)_{\gl(x_U)}$.
  For arbitrary $x_V \in \wtd{V}$ lying above $x$, since $x_V$ and $\gl(x_U)$ are in the same orbit of $G_V$, the isotropy groups $(G_V)_{x_V}$ and $(G_V)_{\gl(x_U)}$ are conjugate in $G_V$ and thus isomorphic.
\end{proof}

Now consider two orbifold charts $\mb{U}$ and $\mb{V}$ and a point $x \in U \cap V$ as in the corollary.
If we do not assume the existence of an embedding $\mb{U} \to \mb{V}$, then local compatibility guarantees that the conclusion of the corollary still holds.

\begin{dfn}
  \begin{enumerate}[(1)]
  \item
    Let $\mf{X}$ be an effective orbifold and suppose $x \in |\mf{X}|$.
    Choose an orbifold chart $\mb{U}$ in any atlas in the orbifold structure on $\mf{X}$, with $x \in U$, and choose any $x_U \in \wtd{U}$ lying above $x$.
    Then the \emph{local group} $G_x$ of $x$ is defined to be the isomorphism class of the group $(G_U)_{x_U}$.
  \item
    Let $\mf{X}$ be an effective orbifold.
    The \emph{singular set} of $\mf{X}$ is the set
    \begin{equation*}
      \gS(\mf{X}) := \{x \in |\mf{X}| \mid G_x \neq \mb{1}\}
    \end{equation*}
    of points with nontrivial local group.
    The points in $\gS(\mf{X})$ are called \emph{singular points}.
  \end{enumerate}
\end{dfn}

The above discussion shows that these definitions are independent of all choices made therein.
In particular, they are independent of the choice of atlas on $\mf{X}$.

In Appendix \ref{gpactions}, we show two important facts about smooth actions of finite groups on manifolds: they can be viewed in neighbourhoods of fixed points as orthogonal actions, and consequently (in the case of effective actions) their fixed point sets are nowhere dense.
Since effective orbifolds are locally given by smooth effective actions of finite groups on Euclidean spaces, we can apply these results locally, resulting in the following two propositions.

\begin{prop}
  The singular set $\gS(\mf{X})$ of an effective orbifold $\mf{X}$ is nowhere dense in $|\mf{X}|$.
\end{prop}

\begin{proof}
  This follows directly from Corollary \ref{nowheredense}.
\end{proof}

\begin{prop}\label{orth}
  Let $\mf{X}$ be an effective orbifold.
  Then there exists an orbifold atlas $\mc{U}_{\operatorname{orth}}$ on $\mf{X}$ such that every chart in $\mc{U}_{\operatorname{orth}}$ is of the form
  \begin{equation}\label{orthchart}
    (\wtd{U},G_x,\gf)
  \end{equation}
  where
  \begin{itemize}
  \item $x$ is a point in $|\mf{X}|$,
  \item $G_x$ is (isomorphic to) the local group of $x$,
  \item $\wtd{U}$ contains the origin, which is mapped to $x$ by $\gf$, and
  \item $G_x$ acts orthogonally on $\wtd{U}$.
  \end{itemize}
  Furthermore, there exists such a chart in $\mc{U}_{\operatorname{orth}}$ for each $x \in |\mf{X}|$. We call such an atlas an \emph{orthogonal atlas}.
\end{prop}

\begin{proof}
  Let $\mc{U}$ be an orbifold atlas on $\mf{X}$, and consider a chart $\mb{U} \in \mc{U}$ and a point $x \in \wtd{U}$.
  The isotropy group $(G_U)_x$ acts on $\wtd{U}$ and fixes $x$, so by the orthogonalisation lemma \ref{orthogonalisation} there exists a coordinate neighbourhood $\wtd{V}$ centred at $x$ upon which $(G_U)_x$ acts orthogonally.
  Denote the coordinate map $\map{\gy}{N_x}{\wtd{V}}$, where $N_x$ is a small open neighbourhood of the origin in some Euclidean space.
  The triple
  \begin{equation*}
    \big(N_x,(G_U)_x,\gf_U \circ \gy\big)
  \end{equation*}
  then defines an orbifold chart on $|\mf{X}|$.
  To construct the atlas $\mc{U}_\text{orth}$, take the set of all such charts, noting that we may need to take arbitrarily small neighbourhoods $\wtd{V}$ of each $x$ to guarantee local compatibility.
\end{proof}

The atlas we constructed on $\ZZ_n \sm B$ following Definition \ref{eff} is an example of an orthogonal atlas.
The existence of orthogonal atlases shows that we need not consider arbitrary finite group actions in local charts in order to deal with orbifolds.
However, in concrete examples such as $\ZZ_n \sm B$, the `simplest' atlases are not necessarily orthogonal.

\begin{rmk}
  If $\mf{X}$ is an orbifold with empty singular set, then an orbifold atlas for $\mf{X}$ yields a manifold atlas on $|\mf{X}|$, and equivalence of orbifold atlases corresponds to equivalence of manifold atlases.
  Likewise, any manifold $X$ can be made into an orbifold $\mf{X}$ with $|\mf{X}| = X$ by appending the trivial action of the trivial group to each coordinate chart.
  Therefore we consider orbifolds with empty singular sets to be manifolds, and vice versa.
\end{rmk}

\section{Further examples}\label{examples}

So far, the only explicit example of an orbifold that we have considered is that given by the action of a finite group on an open subset of Euclidean space.
Orbifolds are defined to be locally of this form, but we have not yet considered orbifolds which are different globally.

The most important example of an effective orbifold is given by a compact Lie group $G$ acting smoothly, effectively, and almost freely on a smooth manifold $M$.
The key result needed to construct an orbifold from such an action is the slice theorem \ref{slice}.

\begin{prop}\label{quotient}
  Let $G$ and $M$ be as above.
  Then there is a natural orbifold structure on the orbit space $G \sm M$.
\end{prop}

\begin{proofsketch}
  Let $x$ be a point of $M$.
  By the slice theorem \ref{slice}, there exists a $G$-neighbourhood $N$ of the orbit $G(x)$ which is $G$-equivariantly diffeomorphic to the twisted product $G \times_{G_x} U$, where $U$ is a coordinate $G_x$-neighbourhood of $x$.
  The triple $(U,G_x,\gf)$ is then an orbifold chart on $G \sm M$, where $U$ has been identified with an open subset of $\RR^n$ and where $\map{\gf}{U}{G \sm M}$ is the quotient map.
  Note that $G_x \sm U$ is homeomorphic to $G \sm N$, since $G \sm N$ is homeomorphic to $G \sm (G \times_{G_x} U)$, which can be identified with $G_x \sm U$.
  As in the proof of Proposition \ref{orth}, the union of all such charts will define an orbifold atlas on $G \sm M$, again noting that we may have to take arbitrarily small neighbourhoods $N$ to guarantee local compatibility.
\end{proofsketch}

The orbit space $G \sm M$ along with the above orbifold structure defines an orbifold: all that remains to be shown is that $G \sm M$ is paracompact and Hausdorff, and this follows once we assume that the manifold $M$ is paracompact and Hausdorff.\footnote{If $G$ is not assumed to be compact, then we have to assume that $G$ acts properly on $M$. If we do not assume this, $G \sm M$ may not be Hausdorff.}
An orbifold constructed in this way is called an \emph{effective quotient orbifold}, and may be denoted $G \sm M$ when no confusion with the orbit space $G \sm M$ is likely to occur.
This construction yields a wide variety of interesting examples.

\begin{example}\label{wp} \textbf{(Weighted projective spaces.)}
  Consider an odd-dimensional sphere
  \begin{equation*}
    S^{2n+1} = \left\{ z \in \CC^{n+1} : \norm{z} = 1\right\} \subset \CC^{n+1}
  \end{equation*}
  where $z = (z_0,\ldots,z_n)$ and $\norm{z} = \sum_{j=0}^n |z_j|^2$.
  Let $(a_0,\ldots,a_n)$ be an $(n+1)$-tuple of integers, and let the circle group $S^1$ act on $S^{2n+1}$ by
  \begin{equation*}
    \gl(z_0,\ldots,z_n) := (\gl^{a_0}z_0,\ldots,\gl^{a_n}z_n).
  \end{equation*}
  The action of $S^1$ on $S^{2n+1}$ is effective if and only if the integers $a_0,\ldots,a_n$ are coprime (in the sense that $\gcd(a_0,\ldots,a_n) = 1$), for if $b \neq 1$ divides each $a_j$, then the $b$-th root of unity $\gz_b$ acts as the identity on every element of $S^{2n+1}$.
  The resulting orbifold, called a \emph{weighted projective space}, is denoted $\WP(a_0,\ldots,a_n)$.
  
  As an explicit example, consider $\WP(2,3)$.
  In this case we have an action of $S^1$ on $S^3$ given by
  \begin{equation*}
    \gl(z_0,z_1) = (\gl^2 z_0, \gl^3 z_1).
  \end{equation*}
  Let us compute the singular set $\gS(\WP(2,3))$.
  Suppose $z=(z_0,z_1) \in S^3$, $\gl \in S^1$, and $\gl z = z$.
  Then we have
  \begin{equation*}
    z_0 = \gl^2 z_0 \quad \text{and} \quad z_1 = \gl^3 z_1.
  \end{equation*}
  If both $z_0$ and $z_1$ are nonzero, then upon dividing the above equalities by $z_0$ and $z_1$ respectively we see that $1 = \gl^2 = \gl^3$.
  The only complex number $\gl \in S^1$ for which this holds is the identity, so that such a point $(z_0 \neq 0, z_1 \neq 0)$ has trivial isotropy group.
  If $z_0 = 0$, then since $\norm{z} = 1$, $z_1$ must lie in $S^1$ (embedded in $S^3$ as the subspace $\{z \in S^3 \mid z_1 = 0\}$).
  All such points lie in the orbit of $(0,1)$.
  Now $1 = \gl^2$ if and only if $\gl = \pm 1$, so that the isotropy group of $(0,1)$ is isomorphic to the cyclic group $\ZZ_2$.
  Likewise, if $z_1 = 0$ then $z_0 = 1$; all such points lie in the orbit of $(1,0)$, and the isotropy group of this point is isomorphic to $\ZZ_3$.
\end{example}

A nice collection of similar examples can be found in \cite[\textsection 1.2]{ALR07}.

\begin{example}
  We need not restrict ourselves to actions of compact groups on manifolds; indeed, if a Lie group $G$ acts smoothly, effectively, and \emph{properly} on a manifold $M$, then the orbit space $G \sm M$ can be made into an orbifold exactly as in Proposition \ref{quotient}.
  In particular, if $G$ is also assumed to be discrete, we obtain what is called a \emph{good} orbifold; if $G$ is finite this is called an \emph{effective global quotient}.
  Note that for a proper action of a discrete group the orbit space $G \sm M$ is naturally a smooth manifold \cite[\textsection 0.4, Example 4.8]{mdC92}, but that it contains no `orbifold data': for example, for all $n \in \NN$ the orbit spaces $\ZZ_n \sm B$ from Example \ref{disc} are diffeomorphic, but the \emph{orbifolds} $\ZZ_n \sm B$ are distinct.
\end{example}

\begin{figure}[t]
  \centering
  \includegraphics[keepaspectratio=true,scale=1.8]{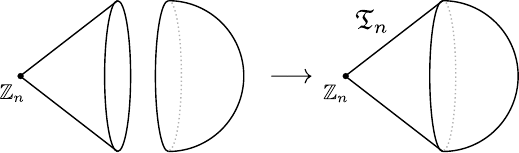}
  \caption{Construction of the teardrop orbifold $\mf{T}_n$.}
  \label{tdfig}
\end{figure}

\begin{example}\label{teardrop}
  Consider the construction of an orbifold atlas on $\ZZ_n \sm B$ following Remark \ref{ineff}.
  We can identify the orbit space $X$ with a cone and smoothly glue a hemisphere onto the open end, as in Figure \ref{tdfig}.
  The result is an orbifold with orbit space homeomorphic to $S^2$, but with a lone singular point with local group $\ZZ_n$.
  This orbifold is called a \emph{teardrop}, and we denote it by $\mf{T}_n$.
  It can be shown that teardrops are not effective global quotients \cite[Example 1.15]{ALR07}.
  
  Similarly, we can glue two different `orbifold cones' together.
  This yields another orbifold with orbit space homeomorphic to $S^2$, this time with two singular points with local groups $\ZZ_n$ and $\ZZ_m$.
  These orbifolds are not effective global quotients unless $n=m$ \cite[Example 1.16]{ALR07}.
\end{example}

\begin{rmk}\label{bigcharts}
  Since an effective global quotient is an orbifold, it is possible to relax our definition \eqref{chart} of orbifold charts by allowing $\wtd{U}$ to be a connected manifold, rather than a connected subset of $\RR^n$.
  This definition can be found in the literature, for example in \cite[\textsection 14.1]{jD11}.
  We will occasionally construct orbifold atlases using this characterisation.
\end{rmk}

\section{Smooth maps}

Defining smooth maps between smooth manifolds is quite easy: a map $\map{f}{M}{N}$ between manifolds is said to be smooth if, in all local coordinates, it is smooth as a map from $\RR^{\dim M}$ to $\RR^{\dim N}$.
We encounter a problem when we try to generalise this condition to orbifolds.
If we have two orbifolds $\mf{X}$, $\mf{Y}$ and a continuous map $\map{f}{|\mf{X}|}{|\mf{Y}|}$, we can say that $f$ is smooth if it can locally be lifted to a $G_U$-$G_V$-equivariant smooth map $\wtd{U} \to \wtd{V}$ for orbifold charts $\mb{U}$, $\mb{V}$ on $\mf{X}$, $\mf{Y}$.
However, there may be more than one possible choice of lifting, and we are faced with the question of deciding when two lifts should be considered to define the same map.

Satake's original definition \cite[\textsection 2]{iS57} of `smooth map' seems to be reasonable, but does not behave well under composition or with respect to orbifold vector bundles.
More satisfactory definitions have been given, for example in \cite[Definition 4.1.5]{CR02}.
Another reasonable way to define smooth maps of orbifolds is via \emph{groupoids}; we will take this up in Chapter \ref{gpoids}.
For now, we shall stick with Satake's definition and present the problems as they occur.

\begin{dfn}\label{smoothmaps}
  Let $\mf{X}$ and $\mf{Y}$ be two orbifolds, and let $\mc{U}$ and $\mc{V}$ be atlases on $\mf{X}$ and $\mf{Y}$ respectively.
  A \emph{smooth map} $\map{f}{\mf{X}}{\mf{Y}}$, defined with respect to the above atlases, is given by a map $\map{\gm}{\mc{U}}{\mc{V}}$ and a collection of smooth maps
  \begin{equation}\label{smmap}
    \left\{\map{f_{\mb{U}}}{\wtd{U}}{\wtd{\gm(U)}}\right\}_{\mb{U} \in \mc{U}},
  \end{equation}
  where we write
  \begin{equation*}
    \gm(\mb{U}) =  \left(\wtd{\gm(U)},G_{\gm(U)},\gf_{\gm(U)}\right),
  \end{equation*}
  such that for any embedding $\map{\gl}{\mb{U}}{\mb{V}}$ there exists an embedding $\map{\gl^\prime}{\gm(\mb{U})}{\gm(\mb{V})}$ such that
  \begin{equation*}
    \gl^\prime \circ f_\mb{U} = f_\mb{V} \circ \gl.
  \end{equation*}
\end{dfn}

\begin{rmk}
  Note that in \eqref{smmap}, $\wtd{\gm(U)}$ is defined to be the Euclidean space associated to the orbifold chart $\gm(\mb{U})$.
  It is not $\gm(\wtd{U})$---indeed, this is not even defined.
  This ambiguity could be considered a flaw of our notational system.
\end{rmk}

We can then go on to define refinements of maps defined with respect to pairs of atlases, and use this to define an equivalence relation on such locally-defined maps.
A smooth map $\mf{X} \to \mf{Y}$ is then an equivalence class of locally-defined maps.
A locally-defined map $\map{f}{\mf{X}}{\mf{Y}}$ induces a continuous map from $|\mf{X}|$ to $|\mf{Y}|$, which we denote by $|f|$; this induced map depends only on the equivalence class of $f$, and so we can speak of the map $|f|$ induced by a smooth map $\map{f}{\mf{X}}{\mf{Y}}$.

\begin{example}\label{tdmap}
  Let $\mf{T}_n$ be a teardrop orbifold as in Example \ref{teardrop}.
  Identify $|\mf{T}_n|$ with the sphere $S^2$ in such a way that the singular point $x$ of $\mf{T}_n$ corresponds to the `north pole' (i.e. any predetermined point) of $S^2$.
  Take an orbifold atlas for $\mf{T}_n$ such that there is a single chart centred at $x$ of the form $(B,\ZZ_n,\gf)$ and such that the remaining charts comprise a maximal manifold atlas for $S^2 \sm \{x\}$ (as in the construction following Remark \ref{ineff}).
  We can construct a smooth map $\map{f}{\mf{T}_n}{\mf{T}_n}$ which is the identity map when viewed in the chart at the singular point $x$, but which is rotation by $2\gp/n$ radians in the `north/south axis' when viewed in all other charts.
  The induced map $\map{|f|}{S^2}{S^2}$ will be rotation by $2\gp/n$ in this axis.
  To construct this map, write the atlas for $\mf{T}_n$ as
  \begin{equation*}
    \mc{U} = \{(B, \ZZ_n,\gf)\} \cup \mc{U}^\prime,
  \end{equation*}
  where $\mc{U}^\prime$ corresponds to the maximal manifold atlas on $S^2 \sm \{x\}$.
  The map $\map{\gm}{\mc{U}}{\mc{U}}$ sends $\mb{B}$ to itself and acts on $\mc{U}^\prime$ by `rotating the charts', i.e. by sending a coordinate chart $\mb{U} = (\wtd{U},\mb{1},\gf_U)$ to the chart $\gm(\mb{U}) := (\wtd{U},\mb{1},r \circ \gf_U)$ where $\map{r}{S^2 \sm \{x\}}{S^2 \sm \{x\}}$ is the rotation map.
  The chart $\gm(\mb{U})$ is in $\mc{U}^\prime$ due to maximality of $\mc{U}^\prime$ and smoothness of $r$.
  Define $f_\mb{V}$ to be the identity map for each $\mb{V} \in \mc{U}$.
  Then the condition on embeddings is automatically satisfied, and so $\{f_\mb{V}\}$ defines a smooth map $\mf{T}_n \to \mf{T}_n$.
\end{example}

The previous example is quite subtle: we defined a map from an orbifold to itself such that the maps $\{f_{\mb{V}}\}$ were all identity maps, but such that the induced map $|f|$ was not the identity.
This shows the role played by the the map $\gm$ between orbifold atlases.
This example also shows that maps between orbifolds are difficult to define in terms of atlases.
Note that if we attempted to define this map from $\mf{T}_n$ to $\mf{T}_m$ with $n \neq m$, then we would run into problems.
Indeed, if $\gl$ is an embedding from the chart $(B,\ZZ_n,\gf)$ to itself, then there would have to be a corresponding embedding $\gl^\prime$ from $(B,\ZZ_m,\gf^\prime)$ to itself (where $(B,\ZZ_m,\gf^\prime)$ is the corresponding chart of $\mf{T}_m$) such that $\id_B \circ \gl = \gl^\prime \circ \id_B$.
In particular, the maps $\gl$ and $\gl^\prime$ must be equal as maps from $B$ to itself.
Taking $\gl$ to be rotation by $2\gp/k$ for some $0 \leq k < n$, we get that rotation by $2\gp/k$ defines an embedding from $(B,\ZZ_m,\gf^\prime)$ to itself for all such $k$.
This is only true if $2\gp/k$ is an integer multiple of $2\gp/m$ for $0 \leq k < n$, and this implies that $n$ divides $m$.
If this is the case, then this construction does indeed define a smooth map $\mf{T}_n \to \mf{T}_m$; otherwise the construction fails.

\begin{rmk}
  We will not go to the trouble of defining compositions of smooth maps, since we will ultimately redefine smooth maps in terms of groupoids.
  Of course, compositions can be defined using these atlas-based definitions (full details are in \cite[\textsection 6.2]{aP10}), but to do so would require a more in-depth look at orbifold atlases.
\end{rmk}

Even though we have not defined compositions of smooth maps, we are still able to discuss diffeomorphisms using the definition in \cite[Definition 1.4]{ALR07}.\footnote{Note that in \cite{ALR07}, the equivalence relation between smooth maps defined in terms of atlases is not considered. Indeed, this is not necessary where diffeomorphisms are concerned.}

\begin{dfn}
  Let $\mf{X}$ and $\mf{Y}$ be orbifolds.
  A smooth map $\map{f}{\mf{X}}{\mf{Y}}$ is a \emph{diffeomorphism} if there exists a second smooth map $\map{g}{\mf{Y}}{\mf{X}}$ such that the induced maps $|f|$ and $|g|$ are mutual inverses.
\end{dfn}

We immediately see that the map $\map{f}{\mf{T}_n}{\mf{T}_n}$ constructed in Example \ref{tdmap} is an example of a diffeomorphism.

\section{Vector and principal bundles}\label{vapb}

Recall that a (real) \emph{vector bundle} over a topological space $X$ is, roughly speaking, given by a collection of real vector spaces
\begin{equation*}
  E = \bigsqcup_{x \in X} E_x
\end{equation*}
indexed by the points of $X$.
We have a natural projection map $\map{\gp}{E}{X}$, sending a vector $v \in E_x$ to the point $x$, and the total space $E$ is topologised in such a way that for each $x \in X$, there is a neighbourhood $U$ in $X$ such that the preimage
\begin{equation*}
  \gp^{-1}(U) = \bigsqcup_{x \in U} E_x
\end{equation*}
is homeomorphic to the product $U \times \RR^n$.
If $X$ is connected, then this $n$ is constant, and called the \emph{rank} of the vector bundle $E \to X$.
We must also ensure that the natural vector space topology on each $E_x$ agrees with that induced from $E$.
For a proper exposition of the basic theory of vector bundles, the reader should consult \cite[Chapter 1]{mA67}.

One can construct a vector bundle of rank $n$ over a space $X$ by taking an open cover $\{U_\ga\}_{\ga \in A}$ of $X$ and `gluing together' the products $U_\ga \times \RR^n$ along the intersections $U_\ga \cap U_\gb$ in a certain way.
This method will be used to construct vector bundles over orbifolds, so we'll briefly describe how it works for topological spaces.
For each pair $\ga,\gb$ of indices, we need a \emph{transition function}
\begin{equation*}
  \map{g_{\ga\gb}}{U_\ga \cap U_\gb}{\GL(n,\RR)}.
\end{equation*}
The collection of all such transition functions encodes how we need to glue together the products $U_\ga \times \RR^n$: if $x$ is in the intersection $U_\ga \cap U_\gb$, then a point $(x,v) \in U_\ga \times \RR^n$ will be glued to the point $(x,g_{\gb\ga}(x)(v)) \in U_\gb \times \RR^n$.
Now suppose that $x$ lies in the triple intersection $U_\ga \cap U_\gb \cap U_\gg$, and consider a point $(x,v) \in U_\ga \times \RR^n$.
We have defined two different ways of gluing this point to a point of $U_\gg \times \RR^n$: by the above definition, $(x,v)$ will be glued to $(x,g_{\gg\ga}(x)(v))$; however, $(x,v)$ will also be glued to $(x,g_{\gb\ga}(x)(v)) \in U_\gb \times \RR^n$, and this point will be glued to $(x,g_{\gg\gb}(x)g_{\gb\ga}(x)(v)) \in U_\gg \times \RR^n$.
For our construction to work, gluing must be transitive.
This leads to the \emph{cocycle condition}
\begin{equation}\label{cocycle}
  g_{\gg\gb}(x)g_{\gb\ga}(x) = g_{\gg\ga}(x) \qquad \text{for all $x \in U_\ga \cap U_\gb \cap U_\gg$}
\end{equation}
which must hold for all nonempty triple intersections $U_\ga \cap U_\gb \cap U_\gg$.
The above gluing procedure will then define a vector bundle over $X$.
The full details of this construction can be found in \cite[\textsection 3]{nS51}.

As with smooth maps of orbifolds, the classical definition \cite[\textsection 3]{iS57} of an abstract vector bundle over an orbifold is inconvenient at best.\footnote{At worst, the classical definitions of orbifold vector bundles and smooth maps are incompatible, as under these definitions the pullback of a vector bundle by a smooth map can fail to be a vector bundle (see \cite[\textsection 2.4]{ALR07}).}
In Section \ref{VBs} we will provide a satisfactory definition within the groupoid framework.
For now we will give an intuitive description of how orbifold vector bundles should behave, along with a few specific examples.

Let $\mf{X}$ be an orbifold with an atlas $\mc{U}$.
A vector bundle over $\mf{X}$ should locally be given by a product bundle; that is, for each chart $\mb{U} \in \mc{U}$ we should have a product bundle $\wtd{U} \times \RR^n \to \wtd{U}$, and furthermore $G_U$ should act on $\wtd{U} \times \RR^n$ in such a way that the projection $\map{\pr_1}{\wtd{U} \times \RR^n}{\wtd{U}}$ is $G_U$-equivariant.
If $\mb{U},\mb{V} \in \mc{U}$ are charts and $\map{\gl}{\mb{U}}{\mb{V}}$ is an embedding, then there should be a $G_U$-$G_V$-equivariant bundle map $\map{\gl^\prime}{\wtd{U} \times \RR^n}{\wtd{V} \times \RR^n}$ (with $G_V$ acting on $\wtd{V} \times \RR^n$ as above).
Such a map must be of the form $\gl^\prime(x,v) = (\gl(x),g_\gl(x)(v))$, with $\map{g_\gl}{\wtd{U}}{\GL_n(\RR)}$.
This suggests that we can construct vector bundles over orbifolds by a generalisation of the above `gluing' construction for topological spaces.
We shall carry out this construction and show that it defines a `base orbifold' $\mf{E}$ and a smooth map $\mf{E} \to \mf{X}$.

As we saw in the previous paragraph, the transition functions for an orbifold vector bundle should be indexed not by intersections $U_\ga \cap U_\gb$, but by embeddings.
As such, for every embedding $\map{\gl}{\mb{U}}{\mb{V}}$ between two charts in $\mc{U}$, suppose we are given a transition function
\begin{equation}\label{orbitf}
  \map{g_\gl}{\wtd{U}}{\GL(n,\RR)}.
\end{equation}
We assume that these transition functions satisfy the following modified version of the cocycle condition \eqref{cocycle}: for every pair $\map{\gl}{\mb{U}}{\mb{V}}$, $\map{\gm}{\mb{V}}{\mb{W}}$ of composable embeddings, we have
\begin{equation}\label{orbicocycle}
  g_{\gm}(\gl(x))g_{\gl}(x) = g_{\gm \circ \gl}(x) \qquad \text{for all $x \in \wtd{U}$.}
\end{equation}
Using these transition functions we construct an orbifold atlas on the space
\begin{equation*}
  E := \left( \bigsqcup_{\mb{U} \in \mc{U}} \wtd{U} \times \RR^n \right)\bigg/\sim,
\end{equation*}
where $\sim$ is the equivalence relation identifying a pair $(x,v) \in \wtd{U} \times \RR^n$ with $(\gl(x),g_{\gl}(x)(v)) \in \wtd{V} \times \RR^n$ whenever $\gl$ is an embedding from $\mb{U}$ to $\mb{V}$.
This space can be shown to be paracompact and Hausdorff.
The atlas is constructed as follows.
For each chart $\mb{U} \in \mc{U}$, define a chart $\mb{U}^*$ on $E$ by
\begin{equation*}
  \mb{U}^* = \left( \wtd{U}^* \times \RR^n, G_U, \gf_U^* \right)
\end{equation*}
where $\gg \in G_U$ acts (effectively) on $\wtd{U}^* := \wtd{U} \times \RR^n$ by $\gg(x,v) = (\gg(x),g_\gg(x)(v))$, recalling that $\gg$ is an embedding from $\mb{U}$ to itself.
Here $\gf_U^*$ is the composition
\begin{equation*}
  \gf_U^* \colon \wtd{U} \times \RR^n \to \bigsqcup_{\mb{U} \in \mc{U}} \wtd{U} \times \RR^n \to E.
\end{equation*}
Evidently $\gf_U^*$ induces a homeomorphism of $G_U \sm \wtd{U}^*$ onto $\gf_U^*(\wtd{U}^*)$.
Let $\mc{U}^*$ be the collection of all such charts $\mb{U}^*$.
It is easy to see that $\mb{U}$ is an orbifold atlas on $E$: $\mc{U}^*$ obviously covers $E$, and local compatibility is inherited from $\mc{U}$, since each embedding $\map{\gl}{\mb{U}}{\mb{V}}$ induces an embedding $\mb{U}^* \to \mb{V}^*$ which sends $(x,v) \in \wtd{U} \times \RR^n$ to $(\gl(x),g_\gl(x)(v))$.
The pair $(E,\mc{U}^*)$ therefore defines an effective orbifold $\mf{E}$.

For $\mf{E}$ to legitimately claim to be a bundle over $\mf{X}$, we need a smooth map $\map{\gp}{\mf{E}}{\mf{X}}$.
Such a map is easily constructed using the atlases $\mc{U}$ and $\mc{U}^*$ used above.
For each chart $\mb{U}^* \in \mc{U}^*$, $\gp_{U^*}$ maps $\mb{U}^*$ into $\mb{U}$ via the first projection
\begin{equation*}
  \map{\pr_1}{\wtd{U} \times \RR^n}{\wtd{U}}.
\end{equation*}
The reader is invited to check that this defines a smooth orbifold map.

Topologically, the result of this construction is a `vector bundle' $\map{\gp}{E}{X}$, where the fibre of $\gp$ over a non-singular point of $\mf{X}$ is a vector space.
However, the fibre over a singular point $x \in \gS(\mf{X})$ is of the form
\begin{equation*}
  \gp^{-1}(x) \cong G_x \sm \RR^n.
\end{equation*}
To see this, observe that in a local linear chart $(\wtd{U},G_x,\gf)$ the bundle $\mf{E}$ is of the form $(\wtd{U} \times \RR^n, G_x, \gf)$ with $G_x$ fixing $0$; thus $G_x$ acts on $\{0\} \times \RR^n$.
The fibre $\gp^{-1}(x)$ is given by the quotient space of this action.

\begin{example}\label{oftb}\textbf{(Tangent bundle.)}
  Let $\mf{X}$ be an orbifold with atlas $\mc{U}$.
  The \emph{tangent bundle} of $\mf{X}$, denoted $T\mf{X}$, is the orbifold vector bundle constructed using the set of transition functions
  \begin{equation*}
    g_\gl(x) = J_\gl(x),
  \end{equation*}
  where $\map{\gl}{\mb{U}}{\mb{V}}$ is an embedding, $x \in \wtd{U}$, and $J_\gl(x)$ is the Jacobian matrix of the embedding $\map{\gl}{\wtd{U}}{\wtd{V}}$ at the point $x$.
  The cocycle condition in this case follows from the chain rule for differentiation.

  For a concrete example of a tangent bundle, consider a `double cone' $\mf{C}$ from Example \ref{teardrop}.
  Recall that $|\mf{C}|$ is homeomorphic to $S^2$, and that $\mf{C}$ has two nonsingular points $x$ and $y$ with local groups $\ZZ_n$ and $\ZZ_m$ respectively.
  At the nonsingular points, the tangent bundle merely looks like the tangent bundle to a sphere.
  However, at the singular point $x$, $T_x \mf{C}$ looks like the quotient of $\RR^2$ by the rotation action of $\ZZ_n$---that is, $T_x \mf{C}$ is (diffeomorphic to) the `cone orbifold' $\ZZ_n \sm B$ from Example \ref{disc}.
  Likewise, $T_y \mf{C}$ is the cone orbifold $\ZZ_m \sm B$.
\end{example}

As with the tangent bundle, we can construct all the tensor bundles $T^r_s \mf{X}$ (see \cite[Example 5.4]{KN63} for the definition of these tensor bundles over manifolds) used in differential geometry.
This suggests that the usual techniques of differential geometry carry over to orbifolds---we take up this task in Chapter \ref{dg}.
Of course, to do this properly we need a rigorous notion of vector bundle, rather than a vague description.

To finish this section we shall sketch the construction of the \emph{orthonormal frame bundle} $\Or \mf{X}$ of an orbifold $\mf{X}$.
This is used to prove the following theorem.

\begin{thm}\label{presentability}
  Every effective orbifold is diffeomorphic to an effective quotient orbifold.
\end{thm}

We recall, again from Steenrod \cite[\textsection 8.1]{nS51}, the construction of the \emph{principal bundle} associated to a given vector bundle.
We will apply this construction to an orbifold vector bundle directly, rather than restating the construction for ordinary vector bundles.
Suppose we are given a set of transition functions \eqref{orbitf} for a vector bundle over an orbifold $\mf{X}$ with atlas $\mc{U}$.
We can replace the group $\GL(n,\RR)$ with a Lie subgroup $G < \GL(n,\RR)$ such that the images of the transition functions $g_\gl$ are all contained in $G$ (this is called a \emph{reduction of structure group} from $\GL(n,\RR)$ to $G$).
Carrying out the construction of the atlas $\mc{U}^*$, but replacing $\wtd{U} \times \RR^n$ with $\wtd{U} \times G$ and letting $\gg \in G_U$ act on $\wtd{U} \times G$ by
\begin{equation*}
  \gg(x,v) := (\gg(x),g_\gg(x) v),
\end{equation*}
we obtain a \emph{principal bundle} over the orbifold $\mf{X}$ with structure group $G$.\footnote{Recall from Remark \ref{bigcharts} that we can define an orbifold locally as a quotient of a \emph{manifold} by a finite group.}
The technicalities of this construction are essentially the same as those of the construction of an orbifold vector bundle.

When $M$ is a smooth $n$-dimensional manifold, the orthonormal frame bundle $\Or M$ is constructed by reducing the structure group of the tangent bundle $TM$ from $\GL(n,\RR)$ to the orthogonal group $O(n)$.
Such a reduction is equivalent to equipping $M$ with a Riemannian metric.
The same process works for an $n$-dimensional orbifold $\mf{X}$.
In this case, a Riemannian metric on $\mf{X}$ is given by a $G_U$-invariant Riemannian metric $g_U$ on each $T\wtd{U}$ for each orbifold chart $\mb{U}$ in a given atlas for $\mf{X}$, such that the metrics $g_U$ are compatible with embeddings between charts.
The existence of Riemannian metrics is proved in essentially the same way as for smooth manifolds, using the existence of partitions of unity; we provide a proof in Section \ref{metrics}.
Hence we can construct an orthonormal frame bundle $\Or \mf{X}$ to an orbifold $\mf{X}$.

\begin{prop}
  For each effective orbifold $\mf{X}$, the frame bundle $\Or \mf{X}$ is a manifold.
\end{prop}

\begin{proof}
  Let $\mc{U}$ be an orthogonal atlas for $\mf{X}$.
  Then the orbifold $\Or \mf{X}$ is given by orbifold charts of the form
  \begin{equation*}
    \mb{U}^* = \left(\wtd{U} \times O(n), G_U, \gf_U \circ \pr_1\right)
  \end{equation*}
  for each chart $\mb{U} \in \mc{U}$.
  Here $g \in G_U$ acts on $\wtd{U} \times O(n)$ by
  \begin{equation*}
    g(x,v) = (gx,Dg(v))
  \end{equation*}
  where $Dg(v)$ denotes the componentwise differential action of $g$ on a frame $v$.
  Now suppose a point $(x,v) \in \wtd{U} \times O(n)$ is fixed by $g \in G_U$.
  Then $g$ must be in $(G_U)_x$, and we get that $Dg(v) = v$.
  Since $v$ is a frame, this implies that $\map{Dg}{T_x(\wtd{U})}{T_x(\wtd{U})}$ is the identity map.
  Since the atlas $\mc{U}$ is orthogonal, $g$ acts linearly on $\wtd{U}$, so $Dg = \id_{T_x(\wtd{U})}$ implies that $g = 1$.
  (Here we have used that $G_U$ acts effectively on $\wtd{U} \times O(n)$.)
  Hence the singular set $\gS(\Or \mf{X})$ is empty, and therefore $\Or \mf{X}$ is a manifold.
\end{proof}

Now observe that the (compact) orthogonal group $O(n)$ acts smoothly, effectively, and almost freely on the manifold $\Or \mf{X}$ by fibrewise right translation.
The orbit space of this action is homeomorphic to $|\mf{X}|$, and indeed the effective quotient orbifold $O(n) \sm \Or \mf{X}$ is diffeomorphic to the orbifold $\mf{X}$.
This completes the sketch of the proof of Theorem \ref{presentability}.\footnote{A more complete proof is available in \cite[\textsection 1.3]{ALR07}.}

\section{Ineffective orbifolds}\label{ineffective}

Consider a compact Lie group $G$ acting smoothly and almost freely on a manifold $M$, as in Proposition \ref{quotient} but without the assumption of effectiveness.
The orbit space $G \sm M$ \emph{should} have the structure of an orbifold, but it does not, since we have demanded that local groups act effectively in local charts.
If we relax this condition and allow ineffective local actions, we get what is called an \emph{ineffective orbifold}.
These do occur naturally; for example, the weighted projective spaces $\WP(a_1,\ldots,a_n)$ (see Example \ref{wp}) are ineffective when the integers $a_1,\ldots, a_n$ are not coprime.
Another example is given by the \emph{moduli stack of elliptic curves} \cite{rH11,bF10}: this is the ineffective orbifold given by the quotient of the upper half-plane
\begin{equation*}
  \mc{H} := \{z \in \CC \mid \Im(z) > 0\}
\end{equation*}
by the discrete group $\SL(2,\ZZ)$ of $2 \times 2$ square matrices with integer entries, where $\SL(2,\ZZ)$ acts properly on $\mc{H}$ by M\"obius transformations
\begin{equation*}
  \twmat{a}{b}{c}{d}(z) := \frac{az + b}{cz + d}.
\end{equation*}
The matrix $-I$ acts trivially on every point on $\mc{H}$, so $\SL(2,\ZZ) \sm \mc{H}$ is an ineffective orbifold.

One can treat ineffective orbifolds in terms of orbifold atlases, as we have done with effective orbifolds in this chapter, but this approach is not perfect.
Foundational results such as Lemma \ref{satake1} no longer hold, and the singular set of an ineffective orbifold $\mf{X}$ can be far from nowhere dense in $|\mf{X}|$; indeed, in the case of an `ineffective global quotient' $G \sm M$, the singular set is the whole orbit space.
The orbifold atlas definition of ineffective orbifolds can be found in \cite[\textsection 4.1]{CR02}.
It is similar to the definitions made in Section \ref{effectives}, except for the following caveat: an embedding between charts $\mb{U} \to \mb{V}$ must also come with a prescribed injective group homomorphism $G_U \to G_V$, in order to get around the failure of Lemma \ref{satake1}.
This homomorphism must restrict to an isomorphism from $(G_U)_0$ to $(G_V)_0$, where $(G_U)_0 < G_U$ is the subgroup of elements of $G_U$ which act trivially on $\wtd{U}$.
We will not use this definition, since for our purposes it is simpler and more gratifying to define ineffective orbifolds in terms of groupoids.
This is the goal of the next chapter.

For $G$ and $M$ as above, the ineffective orbifold $G \sm M$ is relatively easy to understand.
Since $G_0$ is a subgroup of each isotropy group $G_x$, and since by assumption each $G_x$ is finite, $G_0$ is itself finite.
The quotient
\begin{equation*}
  G_\text{eff} := G/G_0
\end{equation*}
acts effectively on $M$, and so the orbit space $G_\text{eff} \sm M$ has a natural effective orbifold structure.
This orbit space is equal to $G \sm M$, so we have a natural way of turning $G \sm M$ into an effective orbifold.
We can think of the \emph{ineffective} orbifold $G \sm M$ as an orbifold with a `generic singularity' given by $G_0$.

An orbifold which is diffeomorphic (of course, the definition of diffeomorphism must first be extended to ineffective orbifolds) to such an ineffective global quotient is called \emph{presentable}.
The discussion above shows that presentable orbifolds can be thought of as a straightforward modification of effective orbifolds.
As for non-presentable orbifolds, we have the following interesting conjecture.

\begin{conj}\label{conjecture}(\cite[Conjecture 1.55]{ALR07})
  Every orbifold is presentable.
\end{conj}

We have shown (Theorem \ref{presentability}) that every effective orbifold is presentable, and so far the only examples of ineffective orbifolds we have given are already in the form of a quotient $G \sm M$.
Partial results towards this conjecture have been obtained; in \cite[Theorem 5.5]{HM04} it is shown that every connected orbifold can be written as a quotient $G \sm \mf{M}$, where $G$ is a compact Lie group and $\mf{M}$ is an \emph{purely ineffective} orbifold, i.e. an orbifold for which every local group $G_x$ acts trivially on some orbifold chart about the point $x$.
If this conjecture holds true, then we are in a good position to deal with ineffective orbifolds in general.
If not, then there exist ineffective orbifolds which can't be described in terms of a `generic singularity' as above; these orbifolds could globally behave quite differently to effective orbifolds.
This situation is not well understood, and we will restrict our attention to presentable orbifolds when necessary.

\begin{figure}[t]
  \centering
  \includegraphics[keepaspectratio=true,scale=1.8]{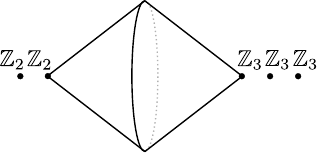}
  \caption{The inertia orbifold of a double cone with local groups $\ZZ_2$ and $\ZZ_3$}
  \label{inertfig}
\end{figure}

To conclude this chapter, we shall give a description of an ineffective orbifold which will be of use in what follows.
Let $\mf{X}$ be an effective orbifold with atlas $\mc{U}$.
The \emph{inertia orbifold} of $\mf{X}$, denoted $\wtd{\mf{X}}$, is an ineffective orbifold defined in terms of the atlas $\mc{U}$ as follows.
The space $\wtd{X} = |\mf{\wtd{X}}|$ is given by the set
\begin{equation*}
  \wtd{X} := \{(x,(g)) \mid x \in X, (g) \subset G_x\}
\end{equation*}
where $(g)$ denotes the conjugacy class of $g$ in the local group $G_x$.
For each chart $\mb{U} \in \mc{U}$ and for each conjugacy class $(g) \subset G_U$, define\footnote{Note that we are now stretching our notational system beyond its limits!} $\wtd{\mb{U}}_{(g)}$ to be the (ineffective) orbifold chart
\begin{equation*}
  \wtd{\mb{U}}_{(g)} := \left( \wtd{U}^{(g)}, Z(g), \wtd{\gf}_{(g)} \right)
\end{equation*}
where $\wtd{U}^{(g)}$ is the submanifold of $\wtd{U}$ consisting of all points fixed by $(g)$, $Z(g)$ is the centraliser of $(g)$ in $G_U$, and
\begin{equation*}
  \wtd{\gf}_{(g)}(x) := (\gf_U(x),(g)).
\end{equation*}
Note that $(g)$ is a conjugacy class in $G_U$ rather that $G_x$, but we can easily identify $(g)$ with a conjugacy class in $G_x$: since $(g)$ fixes $x$, it can be viewed as a conjugacy class in $(G_U)_x$, which can be identified with $G_x$.
Also note that if $g \in G_U$ is not the identity and if $\wtd{U}^{(g)}$ is nonempty, then $g \in Z(g)$ acts trivially on $\wtd{U}^{(g)}$, and so the action of $Z(g)$ on $\wtd{U}^{(g)}$ is not effective.
The collection of all such $\wtd{\mb{U}}_{(g)}$ turns out to be an ineffective orbifold atlas, defining an ineffective orbifold $\wtd{\mf{X}}$.
(Of course, if $\mf{X}$ is a manifold, then $\wtd{\mf{X}} = \mf{X}$ is effective.)

See Figure \ref{inertfig} for an illustration of the inertia orbifold of a `double cone' with $\ZZ_2$ and $\ZZ_3$ singularities.
Observe that the inertia orbifold is not connected, consists of components of differing dimension, and that in this case the local groups $\ZZ_2$ and $\ZZ_3$ act trivially in charts about the corresponding singular points.

\blanknonumber
\chapter{Groupoids}\label{gpoids}

The key ingredient of this thesis is a result of Moerdijk and Pronk \cite[Theorem 4.1, $4 \Rightarrow 1$]{iMdP97}, which gives a correspondence between effective orbifold groupoids and effective orbifolds.
Hence the first goal of this chapter is to develop enough Lie groupoid theory to state and prove this result.
The second goal is to give a rigorous definition of smooth maps between orbifolds.
To achieve this goal, we need to discuss the category-theoretic concept of a \emph{category of fractions} and the corresponding \emph{calculus of fractions}.
Finally, we will use this theory to define vector and principal bundles over groupoids, and to see how they behave with respect to Morita equivalence.
This leads to a satisfactory definition of bundles over orbifolds, as well as an analogue of topological $K$-theory for groupoids and orbifolds.

\section{Lie groupoids and examples}

\begin{dfn}
  A \emph{groupoid}\index{groupoid} is a small category in which every arrow is invertible.
\end{dfn}

By virtue of being a small category, such a groupoid $\mc{G}$ posesses an object set $G_0 := \Obj(\mc{G})$ and an arrow set $G_1 := \Mor(\mc{G})$.
It also comes equipped with \emph{source} and \emph{target} maps $\map{s,t}{G_1}{G_0}$, so that we can write an arrow $\gg \in G_1$ as
\begin{equation*}
  s(\gg) \stackrel{\gg}{\lra} t(\gg).
\end{equation*}
To represent the composition of arrows we have a \emph{composition} (or \emph{multiplication}) map $\map{m}{\fp{G_1}{t}{s}{G_1}}{G_1}$ defined on the fibred product $\fp{G_1}{s}{t}{G_1}$; this maps a pair of composable\footnote{By \emph{composable} we mean that $t(\gg) = s(\gd)$.} arrows $(\gg,\gd)$ to the composition $\gd \circ \gg$.
This is best represented diagrammatically as
\begin{equation*}
  x \stackrel{\gg}{\to} y \stackrel{\gd}{\to} z \quad \stackrel{m}{\longmapsto} \quad x \stackrel{\gd \circ \gg}{\lra} z.
\end{equation*}
Finally, we have a \emph{unit} (or \emph{identity}) map $\map{u}{G_0}{G_1}$ mapping an object $x \in G_0$ to the identity arrow $\map{\id_x}{x}{x}$, and an \emph{inverse} map $\map{i}{G_1}{G_1}$ mapping $\map{\gg}{x}{y}$ to $\map{\gg^{-1}}{y}{x}$.
These maps satisfy a collection of identities representing the axioms defining a category; for example, given an arrow $\gg \in G_1$, the fact that $\gg^{-1}$ is a left inverse for $\gg$ can be written as
\begin{equation*}
  m(\gg,i(\gg)) = u(s(\gg)).
\end{equation*}
Of course, we could define a groupoid $\mc{G} = (G_0,G_1,s,t,m,u,i)$ as being given by these structural maps and identities, and then remark that a groupoid defines a small category, but our definition is easier to remember.

\begin{dfn}
  A \emph{Lie groupoid}\index{groupoid!Lie} is a groupoid $\mc{G}$ such that $G_0$ and $G_1$ are smooth manifolds, $s$ is a submersion, and all remaining structural maps are smooth.
\end{dfn}

By Lemma \ref{fibredsub}, demanding that $s$ is a submersion ensures that the fibred product $G_1 \; {}_s {\times}_t \; G_1$ is a smooth manifold, so that we're able to discuss whether or not the multiplication $m$ is smooth.
Since the inverse map $\map{i}{G_1}{G_1}$ is a diffeomorphism, it follows that $t$ is also a submersion.

\begin{example}\label{actiongpoid}\textbf{(Action groupoids.)}\index{action groupoid}
  Let $M$ be a smooth manifold, and let $G$ be a Lie group acting smoothly on $M$ (from the left), with the action represented by a map $\map{\gQ}{G \times M}{M}$.
  Then we define the \emph{action groupoid} $G \ltimes M$ as follows.
  Set $(G \ltimes M)_0 := M$ and $(G \ltimes M)_1 := G \times M$, and let the source and target maps $\map{s,t}{G \times X}{X}$ be the projection $\pr_2$ and the action $\gQ$ respectively.
  Thus the arrows in $G \ltimes M$ are of the form
\begin{equation*}
    x \stackrel{(g,x)}{\lra} gx.
  \end{equation*}
  The other structural maps come from the group $G$ in the obvious fashion; for example, the inverse map $i$ sends an arrow $(g,x)$ to $(g^{-1},gx)$.
  The projection $\pr_2$ is a submersion and the remaining structural maps are evidently smooth, so $G \ltimes M$ is a Lie groupoid.
  In this way the notion of a Lie groupoid is a generalisation of that of a smooth group action.
  Taking $M$ to be a single point and allowing $G$ to act trivially on $M$, we obtain an action groupoid with a single object, which we denote by $G$.
\end{example}

\begin{example}\label{inducedgpoids}\textbf{(Induced groupoids.)}\index{induced groupoid}
  (See \cite[Section 5.3]{MM03})
  Let $\mc{G}$ be a Lie groupoid, $M$ a smooth manifold, and $\map{f}{M}{G_0}$ a smooth map.
  Under certain conditions we can construct a Lie groupoid $f^* \mc{G}$, called the \emph{induced groupoid}, in which $(f^* \mc{G})_0 = M$ and such that the set of arrows from $x \in M$ to $y \in M$ is identified with the set of arrows from $f(x)$ to $f(y)$ in $G_1$, with multiplication inherited from the multiplication in $\mc{G}$.
  Consider the diagram
  \begin{equation}\label{induced}
    \xymatrix{
      (f^* \mc{G})_1 \ar[d]_{\pr_1^\prime} \ar[rr]^{\pr_2^\prime} & & M \ar[d]^f \\
      \fp{M}{f}{s}{G_1} \ar[d]_{\pr_1} \ar[r]^(.66){\pr_2} & G_1 \ar[d]^s \ar[r]^t & G_0 \\
      M \ar[r]_f & G_0 & 
    }
  \end{equation}
  in which both squares are fibred products, i.e.
  \begin{equation*}
    (f^* \mc{G})_1 = \fp{(\fp{M}{f}{s}{G_1})}{t \circ \pr_2}{f}{M}.
  \end{equation*}
  Since $s$ is a submersion, the fibred product $\fp{M}{f}{s}{G_1}$ is a smooth manifold.
  If we further assume that either $f$ or $t \circ \pr_2$ is a submersion, then $(f^* \mc{G})_1$ is a smooth manifold, and the source and target maps ($\pr_1 \circ \pr_1^\prime$ and $\pr_2^\prime$ respectively) are smooth.
  It is not hard to see that that the multiplication and inverse maps coming from $\mc{G}$ are smooth, though this is rather unwieldy to write down.
  Thus the induced groupoid $f^* \mc{G}$ is defined provided that either $f$ or $t \circ \pr_2$ is a submersion.

  For later reference, we note that from diagram \ref{induced} it follows that the diagram
  \begin{equation}\label{foreshadowing}
    \xymatrix@C=50pt{
      (f^* \mc{G})_1 \ar[r]^{\pr_2 \circ \pr_1^\prime} \ar[d]_{(s,t)} & G_1 \ar[d]^{(s,t)} \\
      M^2 \ar[r]_{f^2} & G_0^2
    }
  \end{equation}
  is a fibred product.
\end{example}

\begin{example}\textbf{(Restrictions.)}
  Let $\mc{G}$ be a Lie groupoid, and suppose $U$ is an open submanifold of $G_0$.
  The inclusion $\map{\gi}{U}{G_0}$ is then a submersion, and we can define the \emph{restriction} of $\mc{G}$ to $U$ to be the induced groupoid
  \begin{equation*}
    \mc{G}|_U := \gi^* \mc{G}.
  \end{equation*}
  As a category, $\mc{G}|_U$ is the full subcategory of $\mc{G}$ determined by  $U \subset \Obj(\mc{G})$.
\end{example}

\begin{example}\label{manifoldgpoids}\textbf{(Manifolds as groupoids.)}
  Let $M$ be a smooth manifold.
  There are two natural ways of viewing $M$ as a Lie groupoid.
  First, we could consider a Lie groupoid with object and arrow space both equal to $M$, with $s=t=i=u=\id_M$, and with trivial multiplication.
  This yields the \emph{unit groupoid}\index{unit groupoid} associated to $M$, which we denote simply by $M$.

  Alternatively, let $U = \{(U_j,\gf_j)\}_{j \in J}$ be an atlas for the manifold $M$, where each $\map{\gf_j}{U_j \subset \RR^n}{X}$ is a coordinate map.
  We shall define a Lie groupoid $M[U]$ associated to this atlas, called a \emph{manifold atlas groupoid}, as follows.
  Let $U_0$ be the disjoint union
  \begin{equation*}
    U_0 := \bigsqcup_{j \in J} U_j
  \end{equation*}
  of the local Euclidan coordinate charts, and let $\map{f}{U_0}{M}$ be the disjoint union of the coordinate maps $\gf_j$.
  Since $f$ is a submersion, we can consider the Lie groupoid $f^* M$ induced from the unit groupoid associated to $M$.
  Thus two objects in $f^* M$ are connected by an arrow if and only if they represent the same point of the manifold $M$.
  We define $M[U] := f^* M$.
\end{example}

\begin{example}\label{localisation}\textbf{(Localisation of a groupoid over an open cover.)}
  We can generalise the idea behind the manifold atlas groupoid of the previous example to general Lie groupoids.
  Let $\mc{G}$ be a Lie groupoid, and suppose $U = \{U_j\}_{j \in J}$ is an open cover of $G_0$.
  Then we can construct a groupoid $\mc{G}[U]$ with object set
  \begin{equation*}
    U_0 := \bigsqcup_{j \in J} U_j
  \end{equation*}
  and with arrow set induced by the disjoint union of the inclusions $U_j \to G_0$.
  The groupoid $\mc{G}[U]$ is called the \emph{localisation of $\mc{G}$ over $U$}.
  Evidently the manifold atlas groupoid $M[U]$ above is equal to the localisation of the unit groupoid $M$ over the open cover $U$.
\end{example}

\begin{example}\label{orbifoldgpoids}\textbf{(Orbifolds as groupoids.)}
  Let $\mf{X}$ be an orbifold with $X := |\mf{X}|$, and let $\mc{U}$ be an orbifold atlas for $\mf{X}$. 
  We can construct an orbifold atlas groupoid $X[\mc{U}]$ analogous to the manifold atlas groupoid of Example \ref{manifoldgpoids}; however, we cannot use the technique of Example \ref{localisation}---there is no global `smooth space' to localise over.
  The groupoid $X[\mc{U}]$ must be constructed directly.\footnote{This is done in \cite[Theorem 4.1.1]{MP99}. Our approach is the same, but we provide more detail.}

  Let the object space $X[\mc{U}]_0$ be the disjoint union
  \begin{equation*}
    X[\mc{U}]_0 := \bigsqcup_{\mb{U} \in \mc{U}} \wtd{U}.
  \end{equation*}
  The objects of $X[\mc{U}]_0$ can thus be written in the form $(x,\wtd{U})$, where $\mb{U}$ is a chart in $\mc{U}$ and $x$ is a point in $\wtd{U}$.
  To construct the orbit space $X$ from the atlas $\mc{U}$, we identify two points $(x,\wtd{U})$ and $(y,\wtd{V})$ if there is a chart $\mb{W} \in \mc{U}$ and two embeddings
  \begin{equation}\label{xar}
    \xymatrix{
      \mb{U} & \mb{W} \ar[l]_\gl \ar[r]^\gm & \mb{V}
    }
  \end{equation}
  along with a point $z \in \wtd{W}$ such that $\gl(z) = x$ and $\gm(z) = y$.
  Hence an arrow $(x,\wtd{U}) \to (y,\wtd{V})$ should be given by a triple $(\gl,z,\gm)$ with $\gl$, $z$, and $\gm$ as above.
  To define such a triple we need the following information: a chart $\mb{W}$, two embeddings $\gl$ and $\gm$ out of $\mb{W}$, and a point $z \in \wtd{W}$.
  Thus we define the space
  \begin{equation}\label{arrowset}
    A := \bigsqcup_{\mb{W} \in \mc{U}} \bigsqcup_{(\gl,\gm) \in E(\mb{W})^2} \wtd{W},
  \end{equation}
  where $E(\mb{W})$ is the set of embeddings from $\mb{W}$ into any other chart of $\mc{U}$.
  An element of $A$ is given by a point $z \in \wtd{W}$ for some $\mb{W} \in \mc{U}$ along with two specified embeddings out of $\mb{W}$, so that $A$ can be identified with the set of triples $(\gl,z,\gm)$ as above.
  As a disjoint union of smooth manifolds, $A$ is itself a smooth manifold.

  The space $A$ is a decent candidate for the arrow set of the groupoid $X[\mc{U}]$, but it is too big.
  If we have a two points $(x,\wtd{U})$ and $(y,\wtd{V})$ in $X[\mc{U}]_0$, a diagram of embeddings
  \begin{equation*}
    \xymatrix{
      & \mb{W}^\prime \ar[d]^\gn & \\
      \mb{U} & \mb{W} \ar[l]_\gl \ar[r]^\gm & \mb{V},
    }
  \end{equation*}
  and a point $z^\prime \in \wtd{W}^\prime$ with $\gl(\gn(z^\prime)) = x$ and $\gm(\gn(z^\prime)) = y$, then we obtain a second triple $(\gl \circ \gn,z^\prime,\gm \circ \gn)$ which ought to define the same arrow as $(\gl,\gn(z),\gm)$.
  This defines an equivalence relation $\sim$ on $A$, and we define $X[\mc{U}]_1$ to be the quotient space $A/\sim$.

  \begin{lem}
    The space $X[\mc{U}]_1$ has a canonical smooth structure.
  \end{lem}

  \begin{proof}
    Let $\map{\gp}{A}{A/\sim}$ be the quotient map.
    Then it suffices to show that for each point $x \in A$ there is a neighbourhood $x \in V \subset A$ such that $\gp$ maps $V$ homeomorphically onto $\gp(V)$ \cite[\textsection 1.5.3]{cT11}.
    Consider an element $(\wtd{W},(\gl,\gm))$ of the disjoint union \eqref{arrowset}.
    Each point of $A$ has such an element as an open neighbourhood, and clearly the equivalence relation $\sim$ is trivial on $(\wtd{W},(\gl,\gm))$, so it suffices to show that $\gp((\wtd{W},(\gl,\gm)))$ is open.
    By the definition of the quotient topology, this is equivalent to showing that the saturation $S$ of $(\wtd{W},(\gl,\gm))$ is open in $A$.

    Suppose $(\gl^\prime,z^\prime,\gm^\prime)$ is in $S$.
    Then there is a point $z \in \wtd{W}$ such that
    \begin{equation*}
      (\gl^\prime,z^\prime,\gm^\prime) \sim (\gl,z,\gm).
    \end{equation*}
    This means either that we have a diagram
    \begin{equation*}
      \xymatrix{
        & \mb{W}^\prime \ar[dl]_{\gl^\prime} \ar[d]^{\gn^\prime} \ar[dr]^{\gm^\prime} & \\
      \mb{U} & \mb{W} \ar[l]_\gl \ar[r]^\gm & \mb{V}
      }
    \end{equation*}
    with $z^\prime \in \wtd{W}^\prime$ and $\gn^\prime(z^\prime) = z$, or a diagram
    \begin{equation*}
      \xymatrix{
        & \mb{W} \ar[dl]_\gl \ar[d]^\gn \ar[dr]^\gm & \\
        \mb{U} & \mb{W}^\prime \ar[l]_{\gl^\prime} \ar[r]^{\gm^\prime} & \mb{V}
      }
    \end{equation*}
    with $\gn(z) = z^\prime$.
    In either case, we claim that there is an open neighbourhood of $(\gl^\prime,z^\prime,\gm^\prime)$ contained in $S$; this will prove that $S$ is open in $A$.
    The proof is the same either way, so we will focus only on the first case.
    Since $\map{\gn^\prime}{\wtd{W}^\prime}{\wtd{W}}$ is an embedding with $\gn^\prime(z^\prime) = z$, there exists a neighbourhood $N$ of $z$ in $\wtd{W}$which is contained in $\gn^\prime(\wtd{W}^\prime)$.
    The set
    \begin{equation*}
      \{ (\gl^\prime,{\gn^\prime}^{-1}(x),\gm^\prime) \mid x \in N\}
    \end{equation*}
    is then an open neighbourhood of $(\gl^\prime,z^\prime,\gm^\prime)$ in $A$ (noting that ${\gn^\prime}^{-1}$ is well-defined on $N$).
    Furthermore, this set is contained in $S$, since if $x$ is in $N$, then $(\gl^\prime,{\gn^\prime}^{-1}(x),\gm^\prime)$ is equivalent to $(\gl,x,\gm) \in (\wtd{W},(\gl,\gm))$.
    Completing the proof in the second case then shows that $S$ is open, and ultimately proves that $X[\mc{U}]_1$ has a canonical smooth structure.
  \end{proof}

  For a triple $(\gl,z,\gm)$ in $A$ (as in \eqref{xar}), write the corresponding arrow in $X[\mc{U}]$ as $[\gl,z,\gm]$.
  The structural maps $s,t,u,i$ for the groupoid $X[\mc{U}]$ are defined by
  \begin{align*}
    s([\gl,z,\gm]) &:= (\gl(z),\wtd{U}) \\
    t([\gl,z,\gm]) &:= (\gm(z),\wtd{V}) \\
    u((x,\wtd{U})) &:= [\id_{\mb{U}},x,\id_{\mb{U}}] \\
    i([\gl,z,\gm]) &:= [\gm,z,\gl];
  \end{align*}
  it is easy to show that these are all well-defined and smooth, and that $s$ is a submersion.
  To define the multiplication map $m$, consider a diagram of embeddings
  \begin{equation*}
    \xymatrix{
      \mb{U} & \mb{X} \ar[l]_\gl \ar[r]^\gm & \mb{V} & \mb{X}^\prime \ar[l]_{\gl^\prime} \ar[r]^{\gm^\prime} & \mb{W}
    }
  \end{equation*}
  along with two points $x \in \wtd{X}$, $x^\prime \in \wtd{X}^\prime$ such that $\gm(x) = \gl^\prime(x)$, so that the arrows $[\gl,x,\gm]$, $[\gl^\prime,x^\prime,\gm^\prime]$ are composable.
  Then since
  \begin{equation*}
    \gf_X(x) = \gf_V(\gm(x)) = \gf_V(\gm^\prime(x^\prime)) = \gf_{X^\prime}(x^\prime),
  \end{equation*}
  the point $\gf_X(x)$ is in the intersection $\wtd{X} \cap \wtd{X}^\prime$.
  Hence by local compatibility there is a chart $\mb{Y} \in \mc{U}$, two embeddings
  \begin{equation*}
    \xymatrix{
      \mb{X} & \mb{Y} \ar[l]_\gs \ar[r]^{\gs^\prime} & \mb{X}^\prime,
    }
  \end{equation*}
  and a point $y \in \wtd{Y}$ such that $\gl(\gs(y)) = \gl(x)$ and $\gm^\prime(\gs^\prime(y)) = \gm^\prime(x)$.
  We define the composition $[\gl^\prime,x^\prime,\gm^\prime]\circ[\gl,x,\gm]$ to be the arrow $[\gl \circ \gs,y,\gm^\prime \circ \gs^\prime]$.
  Again, one can show that this is well-defined and smooth.
  
  This is manifestly more complicated than defining an orbifold in terms of an orbifold atlas; the power of the groupoid approach will come from the ability to define orbifolds \emph{without using atlases}.
\end{example}

Given a Lie groupoid $\mc{G}$ and points $x,y \in G_0$, we can consider the subspace $\mc{G}(x,y) \subset G_1$ of arrows $x \to y$, endowed with the subspace topology.
Taking $x = y$ yields the \emph{isotropy group} $\mc{G}_x := \mc{G}(x,x)$, which has a group structure inherited from the multiplication $m$.\footnote{In fact, $\mc{G}_x$ is a Lie group, and the arrow spaces $\mc{G}(x,y)$ are closed submanifolds of $G_1$; a proof can be found in \cite[Theorem 5.4]{MM03}. Ultimately we'll only consider the case where each $\mc{G}_x$ is finite, so we won't need this information.}
The \emph{orbit}\index{orbit} of $x$ is defined to be the set
\begin{equation*}
  \mc{G}(x) := \{t(\gg) \mid \gg \in G_1, s(\gg) = x\} = t(s^{-1}(x)) \subset G_0
\end{equation*}
along with the subspace topology.
The \emph{orbit space}\index{orbit space!of Lie groupoid} of $\mc{G}$ is the set
\begin{equation*}
  |\mc{G}| := \{\mc{G}(x) \mid x \in G_0\}
\end{equation*}
with the weak topology coming from the projection map $\map{\gp}{G_0}{|\mc{G}|}$ sending a point $x$ to the orbit $\mc{G}(x)$.
When $\mc{G} = G \ltimes X$ is an action groupoid, these notions evidently coincide with the usual action-theoretic notions of isotropy group, orbit, and orbit space.
For a groupoid $M[U]$ coming from a smooth manifold $(M,U)$ as in Example \ref{manifoldgpoids}, all isotropy groups are trivial, and the orbit space is homeomorphic to the manifold $M$.
For an orbifold atlas groupoid $X[\mc{U}]$ as in Example \ref{orbifoldgpoids}, the isotropy group of a point $x \in \wtd{U}$ is given by the isotropy group $(G_U)_x$ in the local group $G_U$, and the orbit space is homeomorphic to $X$.

Note that if $\map{\gg}{x}{y}$ is an arrow in $G_1$, then $\gg$ defines an isomorphism from $\mc{G}_x$ to $\mc{G}_y$ which sends $\map{\gd}{x}{x}$ to $\map{\gg\circ \gd\circ \gg^{-1}}{y}{y}$.
Thus if $x$ is a point in the orbit space $|\mc{G}|$, we can speak of the isotropy group $\mc{G}_x$ as an isomorphism class of groups, just as we can speak of the local group of a point of an orbifold.

Often we like to think of the orbit space as providing the `topology' of $\mc{G}$.
As such, we say the groupoid $\mc{G}$ is \emph{compact}\index{compact} if $|\mc{G}|$ is compact, and \emph{connected}\index{connected} if $|\mc{G}|$ is connected.
For a manifold (or orbifold) atlas groupoid $M[U]$ as above, compactness and connectedness coincide with compactness and connectedness of the manifold (or orbifold) $M$.
Note however that in general $M[U]_0$ is neither compact nor connected.

\section{Homomorphisms and equivalences}

\begin{dfn}
  Let $\mc{G}$ and $\mc{H}$ be Lie groupoids.
  A \emph{homomorphism}\index{homomorphism!of Lie groupoids} $\map{\gf}{\mc{G}}{\mc{H}}$ of Lie groupoids is a functor between the categories $\mc{G}$ and $\mc{H}$ such that the object and morphism maps $\map{\gf_0}{G_0}{H_0}$ and $\map{\gf_1}{G_1}{H_1}$ are smooth.
\end{dfn}

Evidently the composition of two Lie groupoid homomorphisms remains a Lie groupoid homomorphism, so we can form the category $\cat{LGpoid}$ of Lie groupoids and homomorphisms.\footnote{Associativity follows directly from associativity of functor composition, and the identity functor is obviously a Lie groupoid homomorphism.}
We then say two Lie groupoids are \emph{isomorphic} if they are isomorphic in this category.
A morphism $\map{\gf}{\mc{G}}{\mc{H}}$ in this category induces a map from $|\mc{G}|$ to $|\mc{H}|$ which we denote $|\gf|$; it is easy to see that this map is continuous.
Indeed, for an open subset $U \subset H_0$ we have
\begin{equation*}
  |\gf|^{-1}(\gp_H(U)) = \gp_G(\gf_0^{-1}(U)),
\end{equation*}
where $\gp_G$ and $\gp_H$ are the projections of $G_0$ and $H_0$ onto their respective orbit spaces.
Since the topologies on $|\mc{G}|$ and $|\mc{H}|$ are the weak topologies coming from these projections, this shows that $|\gf|$ is continuous.
This correspondence can be interpreted as defining a forgetful functor $\map{| \cdot |}{\cat{LGpoid}}{\cat{Top}}$.

We briefly recall the definition of a natural transformation.
Let $\cat{C}$ and $\cat{D}$ be categories, and suppose that $F$ and $G$ are two functors from $\cat{C}$ to $\cat{D}$.
A \emph{natural transformation} $\gt$ from $F$ to $G$ is a correspondence which associates to each object $x \in \Obj(\cat{C})$ a morphism $\gt(x)$ from $F(x)$ to $G(x)$ such that for all morphisms $\map{\ga}{x}{y}$ in $\Mor(\cat{C})$, the diagram
\begin{equation}
  \xymatrix{
    F(x) \ar[d]_{\gt(x)} \ar[r]^{F(\ga)} & F(y) \ar[d]^{\gt(y)} \\
    G(x) \ar[r]_{G(\ga)} & G(y)
  }
\end{equation}
commutes.
When $\cat{C}$ and $\cat{D}$ are small categories, $\gt$ can be thought of as a map from $\Obj(\cat{C})$ to $\Mor(\cat{D})$ such that the diagram above commutes for all $\ga \in \Mor(\cat{C})$.

\begin{dfn}
  Let $\map{\gf,\gy}{\mc{G}}{\mc{H}}$ be two Lie groupoid homomorphisms.
  A \emph{natural transformation} from the homomorphism $\gf$ to the homomorphism $\gy$ is given by a natural transformation $\gt$ from $\gf$ to $\gy$ (viewed as functors) such that the associated map $\map{\gt}{G_0}{H_1}$ is smooth.
\end{dfn}

\begin{example}\textbf{(Homomorphisms between manifolds and Lie groups.)}
  Let $M$ and $N$ be smooth manifolds viewed as unit groupoids.
  Then every Lie groupoid homomorphism $M \to N$ can be identified with a smooth map between the manifolds $M$ and $N$, and vice versa.
  Similarly, if $G$ and $H$ are Lie groups viewed as Lie groupoids with one object, then a Lie groupoid homomorphism $G \to H$ is the same thing as a smooth group homomorphism between the Lie groups $G$ and $H$.
  Therefore Lie groupoid homomorphisms generalise both smooth maps and smooth group homomorphisms.

  Let $\map{\gf,\gy}{M}{N}$ be Lie groupoid homomorphisms (or equivalently, smooth maps) and suppose $\map{\gt}{\gf}{\gy}$ is a natural transformation.
  Then for each $x \in M$, $\gt(x)$ must be an arrow from $\gf(x)$ to $\gy(x)$.
  But since the only arrows in $N$ are identity arrows, we must have that $\gf(x) = \gy(x)$, and so $\gt$ is the identity natural transformation.

  Alternatively, consider two Lie groupoid homomorphisms $\map{\gf,\gy}{G}{H}$ (or equivalently, smooth group homomorphisms), and suppose $\map{\gt}{\gf}{\gy}$ is a natural transformation.
  Then $\gt$ is completely determined by a choice of group element $\gt \in H$, and we must have
  \begin{equation*}
    \gt\gf(g) = \gy(g)\gt
  \end{equation*}
  for all $g \in G$.
  Thus the natural transformations $\gf \to \gy$ can be identified with the elements $\gt \in H$ such that
  \begin{equation*}
    \gf = \gt^{-1} \gy \gt.
  \end{equation*}
\end{example}

Because every arrow in a groupoid is invertible, with the map sending an arrow to its inverse being smooth, every natural transformation between Lie groupoid homomorphisms is invertible with smooth inverse.
In light of this fact, we will refer to natural transformations as \emph{natural isomorphisms}.
If $\gf$ and $\gy$ are Lie groupoid homomorphisms, then we write $\gf \sim \gy$ to mean that $\gf$ and $\gy$ are naturally isomorphic.
It is easy to see that two naturally isomorphic Lie groupoid homomorphisms induce the same maps of orbit spaces.
We would like to consider naturally isomorphic homomorphisms as defining the same map.
To this end, we define $\cat{LGpoid}^\prime$ to be the quotient category of $\cat{LGpoid}$ obtained by identifying two homomorphisms if they are related by a natural transformation.
By the discussion above, this is a well-defined equivalence relation.

\begin{rmk}\label{gspaces}
  Since our intentions are not yet clear, we should take a moment to explain how we plan on using groupoids to represent orbifolds.
  We are working towards a theory of `generalised smooth spaces'; orbifolds are particular examples of such spaces.
  These spaces are represented by Lie groupoids (with non-isomorphic groupoids potentially representing the same space), and maps between these spaces are essentially represented by Lie groupoid homomorphisms (with distinct homomorphisms potentially representing the same map).
  For example, consider a connected manifold $M$ with two distinct atlases $U$ and $V$, and consider a second manifold $N$ with an atlas $W$.
  The `smooth space' $M$ is represented by the unit groupoid $M$ as well as the atlas groupoids $M[U]$ and $M[V]$, but in general no two of these groupoids are isomorphic.
  Indeed, for the groupoid $M$ to be isomorphic to $M[U]$, the manifold $M$ must be diffeomorphic to the disjoint union of the coordinate charts in $U$, and this is only possible if $M$ is covered by a single chart.
  A similar problem arises when we try to define a map from $M$ to $N$; we can do this by defining a Lie groupoid homomorphism from either $M$, $M[U]$, or $M[V]$ to either $N$ or $N[W]$.
  The following sections will explain how to deal with these issues.
\end{rmk}

\begin{dfn}\label{equivs}
  Let $\map{\ge}{\mc{G}}{\mc{H}}$ be a homomorphism of Lie groupoids.
  \begin{enumerate}[(1)]
  \item
    We call $\ge$ a \emph{weak equivalence}\index{equivalence!weak} if
    \begin{enumerate}[(i)]
    \item \label{we1}
      the composition
      \begin{equation*}
      \fp{G_0}{\ge_0}{s}{H_1} \stackrel{\pr_2}{\lra} H_1 \stackrel{t}{\lra} H_0
    \end{equation*}
    is a surjective submersion, and
  \item \label{we2}
    the square
    \begin{equation}\label{we2fp}
      \xymatrix@C=50pt{
        G_1 \ar[r]^{\ge_1} \ar[d]_{(s,t)} & H_1 \ar[d]^{(s,t)} \\
        G_0 \times G_0 \ar[r]^{\ge_0 \times \ge_0} & H_0 \times H_0 \\
      }
    \end{equation}
    is a fibred product of manifolds.
  \end{enumerate}
\item
  We call $\ge$ a \emph{strong equivalence}\index{equivalence!strong} if it represents an isomorphism in $\cat{LGpoid}^\prime$, i.e. if there exists a second homomorphism $\map{\gh}{\mc{H}}{\mc{G}}$ and two natural isomorphisms $\gh \circ \ge \sim \id_{\mc{G}}$ and $\ge \circ \gh \sim \id_{\mc{H}}$.
  Such a homomorphism $\gh$ is called a \emph{weak inverse} of $\gf$.
\end{enumerate}
\end{dfn}

The set of weak (resp. strong) equivalences is closed under composition, and every strong equivalence is a weak equivalence.\footnote{The first statement is obvious in the case of strong equivalences, and for weak equivalences this is \cite[Proposition 5.12(iii)]{MM03}. The second statement is \cite[Proposition 5.11]{MM03}.}
Lie groupoid isomorphisms are automatically strong equivalences; in particular, identity homomorphisms are strong equivalences.

\begin{example}
  Let $\mc{G}$ be a Lie groupoid and $M$ a smooth manifold, and suppose $\map{f}{M}{G_0}$ is a smooth map such that the composition $\map{t \circ \pr_2}{\fp{M}{f}{s}{G_1}}{G_0}$ is a surjective submersion.
  Then, recalling diagrams \eqref{induced} and \eqref{foreshadowing}, the homomorphism $\map{F}{f^* \mc{G}}{\mc{G}}$ with $F_0 = f$ and
  \begin{equation*}
    \map{F_1 = \pr_2 \circ \pr_1^\prime}{(f^* \mc{G})_1}{G_1}
  \end{equation*}
  is a weak equivalence (with notation as in diagram \eqref{induced}) .
  In particular, if $U$ is an atlas on $M$, then this homomorphism gives a weak equivalence from $M[U]$ to the unit groupoid $M$.
  Similarly, if $U$ and $V$ are atlases on $M$ and $V$ refines $U$, then each refinement $V \to U$ induces a weak equivalence $M[V] \to M[U]$.
  This remains true if we replace manifolds and manifold atlases with orbifolds and orbifold atlases.
\end{example}

One can show that if $\gf$ and $\gy$ are naturally isomorphic homomorphisms, then $\gf$ is a weak (resp. strong) equivalence if and only if $\gy$ is also a weak (resp. strong) equivalence.
Thus the notion of weak equivalence descends from $\cat{LGpoid}$ to the quotient category $\cat{LGpoid}^\prime$, and so we can talk of weak equivalences in $\cat{LGpoid}^\prime$.\footnote{Strong equivalences in $\cat{LGpoid}$ descend to isomorphisms in $\cat{LGpoid}^\prime$, so there is no need to talk of strong equivalences in $\cat{LGpoid}^\prime$.}

\begin{rmk}
  One can show (as in the remarks following \cite[Definition 1.42]{ALR07}) that a weak equivalence $\map{\gf}{\mc{G}}{\mc{H}}$ between two Lie groupoids is a weak equivalence of categories.
  Every weak equivalence of categories has a weak inverse (this follows from the axiom of choice), but such a weak inverse is not necessarily a Lie groupoid homomorphism, as the defining maps will usually fail to be smooth.
  Likewise, a weak equivalence of categories $\mc{G} \to \mc{H}$ will generally not be a weak equivalence of Lie groupoids.
\end{rmk}

\section{Morita equivalence}

Recall from Remark \ref{gspaces} that the notion of Lie groupoid isomorphism is too strong for our purposes.
The notion of weak equivalence is more appropriate: we wish to consider two Lie groupoids $\mc{G}$ and $\mc{H}$ as representing the same `generalised smooth space' if there is a weak equivalence between them.
Speaking more categorically, our `category of generalised smooth spaces' should be constructed from the category $\cat{LGpoid}^\prime$ in such a way that weak equivalences in $\cat{LGpoid}^\prime$ induce isomorphisms in this new category.

This suggests a more general problem: given a category $\cat{C}$ and a collection of morphisms $\gS \subset \Mor(\cat{C})$, does there exist a category $\cat{C}[\gS^{-1}]$ such that, in an appropriate sense, $\cat{C}[\gS^{-1}]$ contains $\cat{C}$ and such that all morphisms in $\gS$ become invertible in $\cat{C}[\gS^{-1}]$?
 The \emph{calculus of fractions} offers conditions on $\gS$ for which such a category exists and is well-behaved, along with a description of the morphisms in this category.
 We give a brief overview of the definitions and theorems (without proofs) that we will require.
 Along with proofs, these can be found in \cite[Chapter 1]{GZ67}.

 \begin{dfn}\cite[1.1]{GZ67}
   Let $\cat{C}$ and $\cat{D}$ be categories.
   \begin{enumerate}[(1)]
   \item
     Let $\map{F}{\cat{C}}{\cat{D}}$ be a functor and suppose $\ge$ is a morphism in $\cat{C}$.
     We say that $F$ \emph{makes $\ge$ invertible} if $F(\ge)$ is invertible in $\cat{D}$.
   \item
     Let $\gS \subset \Mor(\cat{C})$.
     A \emph{category of fractions} for $\gS$ is defined by a category $\cat{C}[\gS^{-1}]$ and a functor $\map{P_\gS}{\cat{C}}{\cat{C}[\gS^{-1}]}$ such that
     \begin{enumerate}[(i)]
     \item
       $P_\gS$ makes all morphisms in $\gS$ invertible, and
     \item
       the pair $(\cat{C}[\gS^{-1}],P_\gS)$ is universal among functors satisfying condition (i), i.e. if a functor $\map{F}{\cat{C}}{\cat{D}}$ makes all morphisms in $\gS$ invertible, then there exists a unique functor $\bar{F}$ making a commutative diagram
       \begin{equation*}
         \xymatrix{
           \cat{C} \ar[d]_{P_\gS} \ar[r]^F & \cat{D} \\
           \cat{C}[\gS^{-1}] \ar@{.>}[ur]_{\bar{F}}. & \\
         }
       \end{equation*}
     \end{enumerate}
   \end{enumerate}
 \end{dfn}

\begin{dfn}\label{rcof}\cite[dual of 2.2]{GZ67}
  Let $\cat{C}$ be a category and $\gS \subset \Mor(\cat{C})$.
  We say that $\gS$ \emph{admits a right calculus of fractions} if
  \begin{enumerate}[(i)]
  \item
    $\gS$ contains all identity morphisms in $\cat{C}$,
  \item
    $\gS$ is closed under composition of composable morphisms,
  \item
    for each diagram
    \begin{equation*}
      \xymatrix{
        X \ar[r]^{\ge} & I & Y \ar[l]_{\gf}
      }
    \end{equation*}
    in $\cat{C}$ with $\ge \in \gS$, there exists a commutative square
    \begin{equation*}
      \xymatrix{
        Y^\prime \ar@{.>}[d]_{\gf^\prime} \ar@{.>}[r]^{\ge^\prime} & Y \ar[d]^\gf \\
        X \ar[r]_\ge & I
      }
    \end{equation*}
    with $\ge^\prime \in \gS$, and
  \item
    for each commutative diagram
    \begin{equation*}
      \xymatrix{
        X \ar@/^/[r]^\gf \ar@/_/[r]_\gy & Y \ar[r]^\ge & Y^\prime \\
      }
    \end{equation*}
    in $\cat{C}$ with $\ge \in \gS$, there exists a morphism $\gh \in \gS$ making a commutative diagram
    \begin{equation*}
      \xymatrix{
        X^\prime \ar@{.>}[r]^\gh & X \ar@/^/[r]^\gf \ar@/_/[r]_\gy & Y \ar[r]^\ge & Y^\prime. \\
      }
    \end{equation*}
  \end{enumerate}
\end{dfn}

If $\gS$ and $\cat{C}$ are as in the above definition, then there exists a well-behaved category of fractions $(\cat{C}[\gS^{-1}],P_\gS)$ for $\gS$.
It is important that we describe the objects and morphisms in this category.
The objects of $\cat{C}[\gS^{-1}]$ are simply the objects of $\cat{C}$,
\begin{equation*}
  \Obj(\cat{C}[\gS^{-1}]) = \Obj(\cat{C}).
\end{equation*}
The morphisms $\cat{C}[\gS^{-1}](X,Y)$ between two objects $X,Y$ of $\cat{C}$ can be described as follows.
Consider the collection $H(X,Y)$ of diagrams in $\cat{C}$ of the form
\begin{equation*}
  \xymatrix{
    X & X^\prime \ar[l]_\ge \ar[r]^\gf & Y
  }
\end{equation*} 
with $\ge \in \gS$.
We shall call these diagrams \emph{generalised morphisms}.
A generalised morphism as above will be denoted by the pair $(\ge,\gf)$.
We declare two such generalised morphisms $(\ge,\gf)$ and $(\gh,\gy)$ to be equivalent if there is a commutative diagram in $\cat{C}$ of the form
\begin{equation*}
  \xymatrix{
    & X^\prime \ar[dl]_\ge \ar[dr]^\gf & \\
    X & X^{\prime\prime\prime} \ar@{.>}[u]^\ga \ar@{.>}[d]_\gb & Y \\
    & X^{\prime\prime} \ar[ul]^\gh \ar[ur]_\gy & \\
  }
\end{equation*}
with $\ge\circ\ga = \gh\circ\gb \in \gS$.
The conditions in Definition \ref{rcof} ensure that this defines an equivalence relation.
By the construction of $\cat{C}[\gS^{-1}]$ in \cite[\textsection I.1 and \textsection I.2]{GZ67}, $\cat{C}[\gS^{-1}](X,Y)$ can be naturally identified by the quotient of $H(X,Y)$ under this equivalence relation.
We denote the equivalence class of a generalised morphism $(\ge,\gf)$ by $[\ge,\gf]$.

To see how composition of morphisms in $\cat{C}[\gS^{-1}]$ works, consider two generalised morphisms given by a diagram
\begin{equation*}
  \xymatrix{
    X & X^\prime \ar[l]_\ge \ar[r]^\gf & Y & Y^\prime \ar[l]_\gh \ar[r]^\gy & Z
  }
\end{equation*}
in $\cat{C}$ with $\ge,\gh \in \gS$.
In order to compose the morphisms $[\ge,\gf]$ and $[\gh,\gy]$, notice that condition (iii) in Definition \ref{rcof} ensures the existence of a commutative diagram
\begin{equation*}
  \xymatrix@R=7pt{
    & & X^{\prime\prime} \ar@{.>}@/_/[dl]_{\gh^\prime} \ar@{.>}@/^/[dr]^{\gf^\prime} & & \\
    X & X^\prime \ar[l]_\ge \ar[r]^\gf & Y & Y^\prime \ar[l]_\gh \ar[r]^\gy & Z
  }
\end{equation*}
with $\gh^\prime \in \gS$.
Since $\gS$ is closed under composition by condition (ii), we end up with an element $(\ge \circ \gh^\prime,\gy \circ \gf^\prime) \in H(X,Z)$, whose equivalence class we define to be the composition $[\gh,\gy] \circ [\ge,\gf]$.
One can easily show that this class is independent of the choice of $(\gh^\prime,\gf^\prime)$.

We return to the category $\cat{LGpoid}^\prime$ of Lie groupoids and homomorphims modulo natural isomorphisms.

\begin{thm}
  Let $\gS \subset \Mor(\cat{LGpoid}^\prime)$ be the set of weak equivalences in the quotient category $\cat{LGpoid}^\prime$.
  Then $\gS$ admits a right calculus of fractions.
\end{thm}

\begin{proofsketch}
  We already know that conditions (i) and (ii) in Definition \ref{rcof} hold, as remarked after Definition \ref{equivs}.
  Conditions (iii) and (iv) follow from the existence of weak fibred products of groupoids \cite[Proposition 5.12(iv)]{MM03}; see \cite[Theorem 3.31]{MM00} for a complete proof.
\end{proofsketch}

Thus we are able to work in the category of fractions $\cat{LGpoid}^\prime[\gS^{-1}]$ in which all weak equivalences become isomorphisms.
A morphism $\mc{G} \to \mc{H}$ in this category is an equivalence class of diagrams in $\cat{LGpoid}^\prime$ of the form
\begin{equation*}
  \xymatrix{
    \mc{G} & \mc{K} \ar[l]_\ge \ar[r]^\gf & \mc{H}
  }
\end{equation*}
where $\ge$ is a weak equivalence; such a diagram is called a \emph{generalised map}\index{generalised map}, and denoted $\map{(\ge,\gf)}{\mc{G}}{\mc{H}}$.
The equivalence class of such a generalised map is denoted by $[\ge,\gf]$.

\begin{dfn}
  A \emph{Morita equivalence}\index{equivalence!Morita} between Lie groupoids $\mc{G}$ and $\mc{H}$ is a generalised map $\map{(\ge,\gf)}{\mc{G}}{\mc{H}}$ in which $\gf$ is a weak equivalence.
  We say $\mc{G}$ and $\mc{H}$ are \emph{Morita equivalent} if there exists a Morita equivalence $\mc{G} \to \mc{H}$.
  If $\ge$ and $\gf$ are strong equivalences, we call $(\ge,\gf)$ a \emph{strong Morita equivalence}
\end{dfn}

\begin{prop}\label{meisom}
  Two Lie groupoids are Morita equivalent if and only if they are isomorphic in $\cat{LGpoid}^\prime[\gS^{-1}]$.
\end{prop}

For a proof of this fact, see \cite[Corollary 4.23]{nL01}.\footnote{Note that Landsman's $\cat{LG}^\prime[S^{-1}]$ is our $\cat{LGpoid}^\prime[\gS^{-1}]$, and that the definition of Morita equivalence used is the `$G$-$H$-bibundle' definition, which is equivalent to our definition (see \cite[Corollary 3.21]{MM00}).}
This tells us that the category $\cat{LGpoid}^\prime[\gS^{-1}]$ is a decent candidate for our `category of generalised smooth spaces'.
In order to study these spaces, we must consider properties of groupoids that are invariant under Morita equivalence.
For example, since a weak equivalence $\mc{G} \to \mc{H}$ induces a homeomorphism of orbit spaces $|\mc{G}| \to |\mc{H}|$, the homeomorphism class of the orbit space is Morita invariant.
Of course, general Lie groupoids are hard to deal with, so we will not spend too much time considering the category $\cat{LGpoid}^\prime[\gS^{-1}]$ for its own sake.
Our goal is to study orbifolds, so study orbifolds we shall.

\section{Orbifold groupoids}

We now restrict our attention to a particular class of Lie groupoids, namely \emph{orbifold groupoids}.
These will turn out to represent orbifolds---both effective and ineffective.
First we must consider \'etale groupoids.

\begin{dfn}
  We say that a Lie groupoid $\mc{G}$ is \emph{\'etale} if $s$ and $t$ are local diffeomorphisms.
\end{dfn}

If $\mc{G}$ is \'etale and $\dim G_0$ is defined (i.e. the connected components of $G_0$ have equal dimension), then $G_1$ and $G_0$ must have equal dimension, so we can define the \emph{dimension}\index{dimension} of $\mc{G}$ to be the number
\begin{equation*}
  \dim \mc{G} := \dim G_1 = \dim G_0.
\end{equation*}
In fact, since $s$ is assumed to be a submersion, if $\dim G_0$ is defined then $\mc{G}$ is \'etale if and only if $\dim G_1 = \dim G_0$.

\begin{prop}\label{ofetale}
  Let $\mf{X}$ be an $n$-dimensional orbifold with atlas $\mc{U}$ and underlying space $X$.
  Then the orbifold atlas groupoid $X[\mc{U}]$ is \'etale.
\end{prop}

\begin{proof}
  The object space $X[\mc{U}]_0$ is a disjoint union of $n$-dimensional manifolds, so to show that $X[\mc{U}]$ is \'etale it is sufficient to show that $\dim X[\mc{U}]_1 = n$.
  Recall from the construction of $X[\mc{U}]$ in Example \ref{orbifoldgpoids} that $X[\mc{U}]_1$ is the quotient of the $n$-dimensional manifold $A$ (defined in \eqref{arrowset}) by an equivalence relation $\sim$.
  The smooth structure on $A/\sim$ is such that the quotient map $\map{\gp}{A}{A/\sim}$ is a local diffeomorphism, and so
  \begin{equation*}
    \dim X[\mc{U}]_1 = \dim A = n,
  \end{equation*}
  and we are done.
\end{proof}

\begin{rmk}
  Let $\mc{G} = G \ltimes X$ be an action groupoid.
  Then since
  \begin{equation*}
    \dim G_1 = \dim G + \dim X = \dim G + \dim G_0,
  \end{equation*}
  the groupoid $\mc{G}$ is \'etale if and only if the group $G$ is discrete.
  In particular, if $G$ is finite, then $G \ltimes X$ is \'etale.
\end{rmk}

Unfortunately the \'etale property is not Morita invariant.
It is proven in \cite[Proposition 5.20]{MM03} that a Lie groupoid $\mc{G}$ is Morita equivalent to an \'etale groupoid if and only if it is a \emph{foliation groupoid}; that is, if for each $x \in G_0$ the isotropy group $\mc{G}_x$ is discrete.
Note that \'etale groupoids are automatically foliation groupoids: all isotropy groups $\mc{G}_x = (s,t)^{-1}(x)$ are discrete since $s$ and $t$ are local diffeomorphisms.

\begin{rmk}
  An action groupoid $G \ltimes X$ is a foliation groupoid precisely when the action of $G$ on $X$ has discrete stabilisers.
  In particular, if $G$ acts almost freely on $X$, then $G \ltimes X$ is a foliation groupoid.
  In the case where $G$ itself is non-discrete, $G \ltimes X$ is not \'etale, though it will be Morita equivalent to an \'etale groupoid.
\end{rmk}

\begin{rmk}\label{smorita}
  Although the \'etale property is not Morita invariant, Morita equivalence behaves quite nicely when restricted to \'etale groupoids.
  Indeed, it is possible to show that a Morita equivalence $(\ge,\gh)$ between two \'etale groupoids is equivalent to a strong Morita equivalence $(\ge^\prime,\gh^\prime)$, in the sense that $(\ge,\gh)$ and $(\ge^\prime,\gh^\prime)$ represent the same morphism in $\cat{LGpoid}^\prime[\gS^{-1}]$.\footnote{To show this, one must first show that if $\ge$ is a weak equivalence of \'etale groupoids then $P_\gS(\ge)$ is equivalent to a strong Morita equivalence. This is most conveniently done by using the alternative definitions of Morita equivalence and weak equivalence appearing in \cite[\textsection 1.3]{aH01} and \cite[Definition 5.2]{mW}, noting that `Morita equivalence' in the second reference is the same as our 'weak equivalence'.}
\end{rmk}

\begin{rmk}\label{germ}
  Let $\mc{G}$ be an \'etale groupoid, take two possibly equal points $x,y \in G_0$, and let $\Diff(x,y)$ denote the set of germs of local diffeomorphisms from $x$ to $y$.\footnote{A local diffeomorphism from $x$ to $y$ is given by a neighbourhood $U$ of $x$, a neighbourhood $V$ of $y$, and a diffeomorphism $\map{f}{U}{V}$ with $f(x) = y$. Two such local diffeomorphisms $\map{f}{U}{V}$, $\map{g}{U^\prime}{V^\prime}$ are said to be equivalent if there exists a neighbourhood $W$ of $x$, with $W \subset U \cap U^\prime$, such that the restrictions of $f$ and $g$ to $W$ are equal. A germ of a local diffeomorphism from $x$ to $y$ is an equivalence class of local diffeomorphisms as above.}
  Suppose $\map{\gg}{x}{y}$ is an arrow in $\mc{G}$.
  Then there exists a neighbourhood $U_\gg$ of $\gg$ in $G_1$ and neighbourhoods $V_x,V_y \subset G_0$ of $x$ and $y$ respectively such that the source and target maps $s,t$ map $U_\gg$ diffeomorphically onto $V_x$ and $V_y$ respectively.
  Through these diffeomorphisms, $\gg$ defines a diffeomorphism $V_x \to V_y$, and hence a germ $\wtd{\gg} \in \Diff(x,y)$.
  This germ is independent of the choice of neighbourhood $U_\gg$, so we have defined a natural map $\mc{G}(x,y) \to \Diff(x,y)$.
  Note that $\Diff(x,x)$ is a group under composition, and that the above map defines a group homomorphism
  \begin{equation*}
    \mc{G}_x \to \Diff(x,x).
  \end{equation*}
\end{rmk}

\begin{dfn}\label{effgpoid}
  An \'etale groupoid $\mc{G}$ is said to be \emph{effective} if for each $x \in G_0$ the group homomorphism $\mc{G}_x \to \Diff(x,x)$ is injective.
\end{dfn}

So the effective \'etale groupoids are those for which the only arrows which act trivially as germs of diffeomorphisms are the identity arrows.

\begin{prop}
  Let $X[\mc{U}]$ be an orbifold atlas groupoid as in Proposition \ref{ofetale}.
  Then $X[\mc{U}]$ is effective.
\end{prop}

\begin{proof}
  Let $(x,\wtd{U})$ be a point in $X[\mc{U}]_0$.
  The isotropy group $X[\mc{U}]_{(x,\wtd{U})}$ consists of all arrows of the form $[\gl,z,\gm]$, where $z$ is in $\wtd{W}$ for some chart $\mb{W} \in \mc{U}$ and $\map{\gl,\gm}{\mb{W}}{\mb{U}}$ are embeddings such that $\gl(z) = \gm(z) = x$.
  Using Lemmas \ref{satake1} and \ref{image},\footnote{In particular, we use the corollary that every embedding $\mb{W} \to \mb{U}$ is of the form $g \circ \gl$ for some $g \in G_W$. See the remark following \cite[Lemma 2]{iS57}.}, we find that $\gm = g \circ \gl$ for some $g \in (G_U)_{x}$.
  From the diagram
  \begin{equation*}
    \xymatrix{
      & \mb{W} \ar[dl]_\gl \ar[d]_\gl \ar[dr]^\gm & \\
      \mb{U} & \mb{U} \ar[l]^{1} \ar[r]_g & \mb{U},
    }
  \end{equation*}
  we find that $[\gl,z,\gm] = [1,x,g]$.
  Two arrows $[1,x,g]$ and $[1,x,h]$ of this form (i.e. with $g,h \in (G_U)_{x}$) are equal if and only if $g=h$, and one can show that
  \begin{equation*}
    [1,x,h] \circ [1,x,g] = [1,x,hg].
  \end{equation*}
  Therefore $X[\mc{U}]_{(x,\wtd{U})}$ is isomorphic to the local group $(G_U)_x$.
  The action of $g \in (G_U)_x$ on a neighbourhood of $(x,\wtd{U})$ as a local diffeomorphism coincides with the action of $G_U$ on $\wtd{U}$, which is effective since $\mf{X}$ is an effective orbifold.
  Therefore the \'etale groupoid $X[\mc{U}]$ is effective.
\end{proof}
 
\begin{dfn}
  We say that a Lie groupoid $\mc{G}$ is \emph{proper} if the if the map
  \begin{equation*}
    \map{(s,t)}{G_1}{G_0 \times G_0}
  \end{equation*}
  is proper.
\end{dfn}

As a direct consequence of this definition, every isotropy group $\mc{G}_x$ ($x \in G_0$) of a proper Lie groupoid is compact.
Properness is preserved under weak equivalence; this is a straightforward topological lemma proved in \cite[Lemma 5.25 and Proposition 5.26]{MM03}.

It can be shown (see \cite[Lemma 3.14]{mT10}) that orbifold atlas groupoids $X[\mc{U}]$ are proper.
Recalling Proposition \ref{ofetale}, we make the following definition.

\begin{dfn}
  A Lie groupoid is called an \emph{orbifold groupoid} if it is both proper and \'etale.
\end{dfn}

An effective orbifold $\mf{X}$ with atlas $\mc{U}$ and underlying space $X$ determines an effective orbifold groupoid $X[\mc{U}]$.
The converse is also true, as we will now show.
The following proofs are adapted from \cite[Theorem 4.1, $4 \Rightarrow 1$]{iMdP97}.

\begin{thm}\label{localaction}
  Let $\mc{G}$ be an orbifold groupoid, and let $x \in G_0$ be a point of $\mc{G}$.
  Then there exist arbitrarily small neighbourhoods $N_x$ of $x$ such that the isotropy group $\mc{G}_x$ acts on $x$; this action corresponds to the natural action of $\mc{G}_x$ near $x$ as local diffeomorphisms.
  Further, for each of these neighbourhoods, the restriction $\mc{G}|_{N_x}$ is isomorphic to the action groupoid $\mc{G}_x \ltimes N_x$
\end{thm}

\begin{proof}
  As in Remark \ref{germ}, for each $\gg \in \mc{G}_x$ choose a neighbourhood $W_\gg$ of $\gg$ in $G_1$ upon which $s$ and $t$ restrict to diffeomorphisms.
  The isotropy group $\mc{G}_x$ is finite, so these neighbourhoods can be chosen to be pairwise disjoint.
  Our first candidate for the neighbourhood $N_x$ is the set
  \begin{equation}\label{ux}
    U_x := \bigcap_{\gg \in \mc{G}_x} s(W_\gg);
  \end{equation}
  each $s(W_\gg)$ is an open neighbourhood of $x$, thus so is $U_x$.
  However, $\mc{G}_x$ may not act on $U_x$, so we need to choose a smaller neighbourhood.

  The next step is to find a neighbourhood $V_x \subset U_x$ of $x$ such that
  \begin{equation}\label{forlem}
    (V_x \times V_x) \cap (s,t) \left( G_1 \sm \bigcup_{\gg \in \mc{G}_x} W_\gg\right) = \varnothing;
  \end{equation}
  such a neighbourhood has the property that if an arrow $\gd \in G_1$ is such that $s(\gd), t(\gd) \in V_x$, then $\gd$ is in $W_\gg$ for some $\gg \in \mc{G}_x$.
  For now we shall assume such a neighbourhood exists; a proof that this is so is given in Lemma \ref{techlem}.
  We note that properness of $\mc{G}$ is used in the construction of $V_x$.

  The neighbourhood $V_x$ is not necessarily $\mc{G}_x$-stable, so we shall shrink $V_x$ to ensure that this is so.
  For each $\gg \in \mc{G}_x$ let $\wtd{\gg}$ denote the diffeomorphism $s(W_\gg) \to t(W_\gg)$ induced by $\gg$.
  Since $V_x \subset s(W_\gg)$ for each $\gg \in \mc{G}_x$ (by \ref{ux}), each $\wtd{\gg}$ is defined on $V_x$, and so we can define
  \begin{equation*}
    N_x := \{y \in V_x \mid \text{$\wtd{\gg}(y) \in V_x$ for all $\gg \in \mc{G}_x$}\} = \bigcap_{\gg \in \mc{G}_x} \wtd{\gg}^{-1}(V_x \cap t(W_\gg)) \cap V_x.
  \end{equation*}
  The second equality shows that $N_x$ is open, and since $x$ is fixed by each $\gg \in \mc{G}_x$, $x$ is in $N_x$.
  If $y \in N_x$ and $\gg \in \mc{G}_x$, then for all $\gd \in \mc{G}_x$
  \begin{equation*}
    \wtd{\gd}(\wtd{\gg}(y)) = \wtd{\gd\gg}(y) \in N_x,
  \end{equation*}
  so $\mc{G}_x$ acts on $N_x$.
  Clearly we can make $N_x$ arbitrarily small.
  
  Now we must show that the restriction $\mc{G}|_{N_x}$ is isomorphic to the action groupoid $\mc{G}_x \ltimes N_x$.
  For each $\gg \in \mc{G}_x$, define
  \begin{equation*}
    O_\gg := W_\gg \cap s^{-1}(N_x).
  \end{equation*}
  If $y \in N_x$ and $\gg \in \mc{G}_x$, then $\wtd{\gg}(y) \in N_x$ (by the definition of $N_x$).
  By the definitions of $\wtd{\gg}$ and $W_\gg$, this shows that
  \begin{equation*}
    O_\gg = W_\gg \cap (s,t)^{-1}(N_x \times N_x).
  \end{equation*}
  By the definition of $V_x$, if $\gd \in (s,t)^{-1}(N_x \times N_x)$ then $\gd$ is in $W_\gg$ for some $\gg \in \mc{G}_x$.
  Consequently
  \begin{equation*}
    (\mc{G}|_{N_x})_1 = (s,t)^{-1}(N_x \times N_x) = \bigcup_{\gg \in \mc{G}_x} O_\gg.
  \end{equation*}
  Since the $W_\gg$ are disjoint, so are the $O_\gg$, and hence the above union is a disjoint union.
  It follows (with a bit of thought) that $\mc{G}|_{N_x}$ is isomorphic to $\mc{G}_x \ltimes N_x$.
  (Note that $(\mc{G}_x \ltimes N_x)_1 = \mc{G}_x \times N_x$, and since $\mc{G}_x$ is finite, this is diffeomorphic to a disjoint union of $|\mc{G}_x|$ copies of $N_x$.)
\end{proof}

\begin{lem}\label{techlem}
  Using the notation of the proof above, there exists a neighbourhood $V_x \subset U_x$ of $x$ such that \eqref{forlem} holds.
\end{lem}

\begin{proof}
  Since $(s,t)$ is a proper map into a compactly generated Hausdorff space, $(s,t)$ is a closed map.
  Thus $(s,t)(G_1 \sm \cup_{\gg \in \mc{G}_x} W_\gg)$ is closed.
  Manifolds are regular spaces, so there exists a closed neighbourhood $R_x$ of the point $(x,x) \in G_0 \times G_0$ which avoids the closed set $(s,t)(G_1 \sm \cup_{\gg \in \mc{G}_x} W_\gg)$.
  Define the open neighbourhood $V_x$ of $x$ by
  \begin{equation*}
    V_x := \pr_1(R_x \cap U_x \times U_x) \cap \pr_2(R_x \cap U_x \times U_x).
  \end{equation*}
  Then since $V_x \times V_x \subset R_x$, we have that $V_x$ satisfies \eqref{forlem}.
\end{proof}

\begin{cor}\label{ofat}
  Let $\mc{G}$ be an effective orbifold groupoid.
  Then there is a canonically defined orbifold atlas on the orbit space $|\mc{G}|$.
\end{cor}

\begin{proof}
  Let $\mc{U}$ be the set of triples of the form
  \begin{equation*}
    (N_x, \mc{G}_x, \gp_x)
  \end{equation*}
  where $x \in G_0$, $N_x$ and $\mc{G}_x$ are as in Theorem \ref{localaction}, and $\gp_x$ is the restriction of the quotient map $G_0 \to |\mc{G}|$ to $N_x$.
  Note that for each $x$ we should take arbitrarily small neighbourhoods $N_x$.
  Since the action of $\mc{G}_x$ on $N_x$ is given by the arrows in $\mc{G}$, the map $\gp_x$ is automatically $\mc{G}_x$-invariant and induces an injective map $N_x \sm G_x \to |\mc{G}|$.
  Hence to show that each such triple is an orbifold chart, it suffices to show that $\gp_x$ is an open map; but this holds since $\gp_x$ is the quotient map of a smooth action of a finite group.

  To show that $\mc{U}$ is an orbifold atlas on $|\mc{G}|$, we are left with proving local compatibility.
  Suppose $(N_x,\mc{G}_x,\gp_x)$ and $(N_y,\mc{G}_y,\gp_y)$ are two charts in $\mc{U}$ such that there is a point $z \in \gp_x(N_x) \cap \gp_y(N_y)$.
  Then there is a point $z^\prime \in G_0$ and two arrows $\map{\gg_1}{z^\prime}{x}$, $\map{\gg_2}{z^\prime}{y}$ in $\mc{G}$.
  Choose a chart $(N_{z^\prime},\mc{G}_{z^\prime},\gp_{z^\prime})$ as above, along with neighbourhoods $W_{\gg_1}$, $W_{\gg_2}$ of $\gg_1$ and $\gg_2$ in $G_1$ such that the following conditions are satisfied:
  \begin{itemize}
  \item
    $s$ and $t$ restrict to local diffeomorphisms on $W_{\gg_1}$ and $W_{\gg_2}$,
  \item
    $s(W_{\gg_1}) = s(W_{\gg_2}) = N_{z^\prime}$, and
  \item
    $t(W_{\gg_1}) \subset N_x$ and $t(W_{\gg_2}) \subset N_y$.
  \end{itemize}
  Then the maps $\map{\wtd{\gg_1}}{N_{z^\prime}}{N_x}$ and $\map{\wtd{\gg_2}}{N_{z^\prime}}{N_y}$ induce embeddings
  \begin{equation*}
    \xymatrix{
      (N_x,\mc{G}_x,\gp_x) & \ar[l] (N_{z^\prime},\mc{G}_{z^\prime},\gp_{z^\prime}) \ar[r] & (N_y,\mc{G}_y,\gp_y).
    }
  \end{equation*}
  Now $\gp_{z^\prime}(z^\prime) = z$, and the conditions above ensure that $\gp_{z^\prime}(N_{z^\prime}) \subset \gp_x(N_x) \cap \gp_y(N_y)$, so this proves local compatibility.
  Therefore $\mc{U}$ is an orbifold chart on $|\mc{G}|$.
\end{proof}

In summary, effective orbifold groupoids have natural orbifold structures on their orbit spaces, and an orbifold atlas on a paracompact Hausdorff space determines an effective orbifold groupoid.
Furthermore, a refinement of orbifold atlases induces a weak equivalence of these orbifold groupoid, and so an equivalence class of atlases corresponds to a Morita equivalence class of effective orbifold groupoids.
With these correspondences in mind, we make the following definition.

\begin{dfn}
  An \emph{orbifold structure} on a paracompact Hausdorff space $X$ is given by an orbifold groupoid $\mc{G}$ and a homeomorphism
  \begin{equation}\label{st}
    |\mc{G}| \to X.
  \end{equation}
  An equivalence of orbifold structures $|\mc{G}| \to X$, $|\mc{H}| \to X$ is given by a Morita equivalence $\map{(\ge,\gh)}{\mc{G}}{\mc{H}}$ such that the diagram
  \begin{equation}\label{eq}
    \xymatrix@C=5pt{
      |\mc{G}| \ar[dr] \ar[rr]^{(\ge,\gh)} & & |\mc{H}| \ar[dl] \\
      & X & \\
    }
  \end{equation}
  commutes.
  An \emph{orbifold} $\mf{X}$ is a paracompact Hausdorff space $X$ along with an equivalence class of orbifold structures.
  A \emph{groupoid representation} of an orbifold $\mf{X}$ is given by a particular choice of orbifold structure \eqref{st} in this equivalence class; the groupoid $\mc{G}$ appearing in this structure is called a \emph{representing groupoid} for $\mf{X}$.
\end{dfn}

\begin{rmk}
  Since the orbit space of a Morita equivalence class of Lie groupoids is only defined up to canonical homeomorphism, to `fix' the topology of an orbifold we must give a particular homeomorphism as in \eqref{st}.
\end{rmk}

Note that we have not assumed that our orbifold groupoids are effective.
An \emph{effective} orbifold is one for which any representing groupoid is effective; an \emph{ineffective} orbifold is one which is not effective (recall that effectiveness is Morita invariant).
This allows us to consider a wider class of orbifolds than before, without having to concern ourselves with ineffective orbifold atlases.
In particular, if $G$ is a compact Lie group acting smoothly, almost freely, and \emph{not necessarily effectively} on a smooth manifold $M$, then the action groupoid $G \ltimes M$ (more specifically, a Morita equivalent \'etale groupoid) defines an ineffective orbifold structure on the orbit space $G \sm M = |G \ltimes M|$.
Therefore the examples given in section \ref{ineffective} take their rightful place as orbifolds.

\begin{dfn}
  An orbifold $\mf{X}$ is called \emph{presentable} if it has a representing groupoid $\mc{G}$ which is Morita equivalent to an action groupoid $G \ltimes M$, where $G$ is a compact Lie group acting smoothly and almost freely on a smooth manifold $M$.
\end{dfn}

\begin{rmk}
  Note that the action groupoid $G \ltimes M$ is not \'etale unless $G$ is discrete, so $G \ltimes M$ will usually not be a representing groupoid for an orbifold.
\end{rmk}

With this new machinery in place, we can give a satisfactory definition of a smooth map between orbifolds.

\begin{dfn}
  A \emph{smooth map} between two orbifolds $\mf{X},\mf{Y}$, with groupoid representations
  \begin{equation*}
    |\mc{G}| \to X, \qquad |\mc{H}| \to Y
  \end{equation*}
  respectively, is given by a morphism $\map{f}{\mc{G}}{\mc{H}}$ in $\cat{LGpoid}^\prime[\gS^{-1}]$.
  Such a morphism induces a unique map $\map{|f|}{X}{Y}$ such that the diagram
  \begin{equation*}
    \xymatrix{
      |\mc{G}| \ar[d] \ar[r]^f & |\mc{H}| \ar[d] \\
      X \ar[r]_{|f|} & Y
    }
  \end{equation*}
  commutes.
\end{dfn}

\begin{rmk}
  If we define a diffeomorphism as an invertible smooth map, then by Proposition \ref{meisom} two orbifolds are diffeomorphic if and only if they have Morita equivalent representing groupoids.
  This justifies our definition of presentability.
\end{rmk}

With $\mf{X}$ and $\mf{Y}$ as in the above definition, if
\begin{equation*}
  |\mc{G}^\prime| \to X, \qquad |\mc{H}^\prime| \to Y
\end{equation*}
are equivalent groupoid representations of $\mf{X}$ and $\mf{Y}$ respectively, then the morphism $f$ uniquely determines a morphism $\map{f^\prime}{\mc{G}^\prime}{\mc{H}^\prime}$ yielding a smooth map from $\mf{X}$ to $\mf{Y}$.
The morphisms $f$ and $f^\prime$ are considered to determine the same map from $\mf{X}$ to $\mf{Y}$.\footnote{To be more precise, we could define an equivalence relation on maps defined with respect to groupoid representations, and define a smooth map between orbifolds to be an equivalence class of these maps. At this point, this level of rigour seems unnecessary, and to fill in the gaps explicitly would make for dry reading.}

\begin{rmk}
  It can be shown that this definition of smooth map is equivalent to the `good maps' introduced in \cite[\textsection 4.4]{CR02} and the `strong maps' of \cite[\textsection 5]{iMdP97}. (See \cite[\textsection 2.4]{ALR07}.)
\end{rmk}

A \emph{smooth function} on an orbifold $\mf{X}$ with groupoid representation $|\mc{G}| \to X$ is given by a smooth map from $\mf{X}$ to the manifold $\RR$ (represented by the unit groupoid $\RR$).
The set of smooth functions $\mf{X} \to \RR$ is denoted by $C^\infty(\mf{X})$, and is naturally identified with the set of generalised maps from $\mc{G}$ to $\RR$.
Let $f \in C^\infty(\mf{X})$ be a smooth function on $\mf{X}$, given by a generalised map
\begin{equation*}
  \xymatrix{
    \mc{G} & \mc{G}^\prime \ar[l]_\ge \ar[r]^\gf & \RR. \\
  }
\end{equation*}
Since $\ge$ is a weak equivalence, the induced map $|\ge|$ is invertible, giving a map $|\gf| \circ |\ge|^{-1}$ from $|\mc{G}|$ to $\RR$.
Composing this map with the homeomorphism $X \to |\mc{G}|$ defines a map from $X$ to $\RR$, which we denote by $|f|$.
This map allows us to discuss `pointwise' behavior of the map $f$: we define
\begin{equation*}
  f(x) := |f|(x)
\end{equation*}
for all $x \in X$, and we define the \emph{support} of $f$ to be the set
\begin{equation*}
  \supp f := \overline{\{x \in X : f(x) \neq 0\}}.
\end{equation*}
Hence we can speak of compactly supported real-valued functions on orbifolds.

\begin{rmk}
  Defining a smooth map between two orbifolds $\mf{X}$, $\mf{Y}$ presented by \'etale groupoids $\mc{G}$, $\mc{H}$ respectively is essentially the same as defining a generalised map
  \begin{equation*}
    \xymatrix{
      \mc{G} & \mc{K} \ar[l] \ar[r] & \mc{H}.
    }
  \end{equation*}
  Observe that the Lie groupoid $\mc{K}$ is weakly equivalent to $\mc{G}$, and hence Morita equivalent to $\mc{G}$.
  Therefore $\mc{K}$ is a foliation groupoid.
  However, we cannot assume that $\mc{K}$ is \'etale.\footnote{If we prohibited the use of foliation groupoids to represent orbifolds, then we would not be able to take advantage of the fact that presentable orbifolds can be represented by action groupoids $G \ltimes M$ for non-discrete groups $G$. Representing presentable orbifolds in this way is crucial when looking at their $K$-theory, so such representations shold be allowed.}
  This will cause us problems in the future.
\end{rmk}

To conclude this section, we note that the groupoid approach to orbifolds allows us to easily define restrictions of orbifolds, using the fact that the restriction $\mc{G}|_U$ of an orbifold groupoid $\mc{G}$ to an open set $U$ is an orbifold groupoid.

\begin{dfn}
  Let $\mf{X}$ be an orbifold with groupoid representation $\map{p}{|\mc{G}|}{X}$, and let $U \subset X$ be an open set.
  Writing $\map{\gp}{G_0}{|\mc{G}|}$ for the quotient map, we define the orbifold $\mf{X}|_U$ by the groupoid representation
  \begin{equation*}
    |\mc{G}|_{(p \circ \gp)^{-1}(U)} \to U.
  \end{equation*}
  where we view $|\mc{G}|_{(p \circ \gp)^{-1}(U)}$ as a subset of $|\mc{G}|$.
\end{dfn}

\section{$\mc{G}$-spaces and bundles}\label{gbundles}

Unlike the classical definition of a fibre bundle over an orbifold, the definition of a fibre bundle over a groupoid is fairly intuitive.
A fibre bundle over a Lie groupoid $\mc{G}$ is given by a fibre bundle $\map{\gp}{E}{G_0}$ (in the usual sense) along with a fibrewise action of the arrows of $\mc{G}$ on $E$.
That is, every arrow $\map{\gg}{x}{y}$ in $\mc{G}$ induces a map $\map{\gg}{E_x}{E_y}$ (linear in the case of a vector bundle).
These induced maps must respect the composition law of $\mc{G}$.
To make this notion precise, we have to discuss what it means for a Lie groupoid to act on a smooth manifold.

\begin{dfn}
  Let $\mc{G}$ be a Lie groupoid, $M$ a smooth manifold, and $\map{\ga}{M}{G_0}$ a smooth map.
  A \emph{left action of $\mc{G}$ on $M$ along $\ga$} is given by a smooth map
  \begin{equation*}
    \map{\gQ}{\fp{G_1}{s}{\ga}{M}}{M}
  \end{equation*}
  such that for all $(\gg^\prime,\gg,y) \in \fp{G_1}{s}{t}{\fp{G_1}{s}{\ga}{M}}$ (writing $\gg y := \gQ(\gg,y)$ when $s(\gg) = \ga(y)$),
  \begin{itemize}
  \item
    $\ga(\gg y) = t(\gg)$,
  \item
    $\id_{\ga(y)} = y$, and
  \item
    $\gg^\prime(\gg y) = (\gg^\prime \gg)y$.
  \end{itemize}
  The map $\ga$ is called the \emph{anchor} of the $\mc{G}$-action, and $M$ is referred to as a \emph{left $\mc{G}$-space}.
  In the same way we can define right $\mc{G}$-spaces.
\end{dfn}

Given such a left $\mc{G}$-space, we can define an action groupoid $\mc{G} \ltimes X$ analogous to that of Example \ref{actiongpoid}.
Observe that for a Lie group $G$ viewed as a groupoid with one object, this notion corresponds to the usual notion of $G$-space (as there is only one possible choice of anchor).

A \emph{$\mc{G}$-map} $\map{f}{X}{Y}$ between two $\mc{G}$-spaces (with anchors $\ga_X$ and $\ga_Y$ respectively) is a smooth map with $\ga_Y \circ f = \ga_X$ such that $f(\gg x) = \gg\cdot f(x)$ for all $(\gg,x) \in \fp{G_1}{s}{\ga_X}{X}$.

\begin{example}\label{gpoid-as-gspace}
  Let $\mc{G}$ be a Lie groupoid.
  Then $G_0$ can be viewed as a $\mc{G}$ space where the anchor is the identity map and $\map{\gg}{x}{y}$ acts by
  \begin{equation*}
    \gg x := y.
  \end{equation*}
\end{example}

\begin{dfn}
  Suppose that $\mc{G}$ is a Lie groupoid and $H$ is a Lie group.
  \begin{enumerate}[(1)]
  \item
    A \emph{vector $\mc{G}$-bundle} $E \to \mc{G}$ is given by a vector bundle $\map{\gp}{E}{G_0}$ and a fibrewise linear left action of $\mc{G}$ on $E$ along $\gp$.
  \item
    A (right) \emph{$H$-principal $\mc{G}$-bundle} $P \to \mc{G}$ is given by a (right) $H$-principal bundle $\map{\gp}{P}{G_0}$ and a $H$-equivariant left action of $\mc{G}$ on $P$ along $\gp$.
  \end{enumerate}
  The term \emph{$\mc{G}$-bundle} will be used generically to refer to either a vector or principal $\mc{G}$-bundle.
\end{dfn}

\begin{rmk}\label{actionbundles}
  When $\mc{G} = G \ltimes X$ is an action groupoid, vector $\mc{G}$-bundles are naturally identified with $G$-equivariant vector bundles over the $G$-space $X$.
  This observation is crucial in studying the $K$-theory of presentable orbifolds (see Section \ref{presentablek}).
\end{rmk}

\begin{dfn}
  Suppose that $\mc{G}$ is a Lie groupoid and that $E,F \to \mc{G}$ are vector (resp. principal) $\mc{G}$-bundles.
  \begin{enumerate}[(1)]
  \item
    A \emph{bundle homomorphism} from $E$ to $F$ is a $\mc{G}$-map $\map{f}{E}{F}$ which is also a vector (resp. principal) bundle homomorphism in the usual sense.
    Bundle isomorphisms are defined as usual to be invertible bundle homomorphisms.
  \item
    A \emph{section} of the $\mc{G}$-bundle $E$ is a $\mc{G}$-map $\map{\gs}{G_0}{E}$ which is also a section of $E$ in the usual sense, where $G_0$ is viewed as a $\mc{G}$-space as in Example \ref{gpoid-as-gspace}.
  \end{enumerate}
\end{dfn}

Consider a vector $\mc{G}$-bundle $\map{\gp}{E}{\mc{G}}$ and a Lie groupoid homomorphism $\map{\gf}{\mc{H}}{\mc{G}}$.
We would like to define a pullback $\mc{H}$-bundle $\gf^* E \to \mc{H}$, generalising the construction of the pullback of an ordinary vector bundle.
We shall outline this construction briefly.
The underlying vector bundle of $\gf^* E$ is the ordinary pullback
\begin{equation*}
  \map{\gp^*}{\gf_0^* E}{H_0}.
\end{equation*}
To make $\gf_0^* E$ into a $\mc{H}$-bundle, we need an action of $\mc{H}$ on $\gf_0^* E$ along $\gp^*$.
Suppose $\map{\gd}{x}{y}$ is an arrow in $\mc{H}$.
Then $\gf_1(\gd)$ is an arrow in $\mc{G}$ from $\gf_0(x)$ to $\gf_0(y)$, and hence induces a map from $E_{\gf_0(x)}$ to $E_{\gf_0(y)}$.
By definition of the pullback, $(\gf_0^* E)_x$ and $(\gf_0^* E)_y$ are identified with $E_{\gf_0(x)}$ and $E_{\gf_0(y)}$ respectively; the arrow $\gd$ induces a map $(\gf_0^* E)_x \to (\gf_0^* E)_y$ by way of this identification.
In this way $\gf_0^* E$ is made into a $\mc{H}$-space, thus defining the pullback $\mc{H}$-bundle $\gf^* E$.

The above construction shows that pullback bundles are defined in the category $\cat{LGpoid}$; the following lemma extends this fact to $\cat{LGpoid}^\prime$ (up to canonical isomorphism).

\begin{lem}\label{canism}
  Let $E$, $\mc{G}$, and $\mc{H}$ be as above, and suppose $\map{\gf,\gy}{\mc{H}}{\mc{G}}$ are naturally isomorphic Lie groupoid homomorphisms.
  Then the pullback $\mc{H}$-bundles $\gf^* E$ and $\gy^* E$ are canonically isomorphic.
\end{lem}

\begin{proofsketch}
  Let $\map{\gt}{\gf}{\gy}$ be a natural isomorphism.
  Then for each point $x \in H_0$, there is an arrow $\gt(x)$ in $G_1$ from $\gf_0(x)$ to $\gy_0(x)$.
  This arrow induces a map $E_{\gf_0(x)} \to E_{\gy_0(x)}$.
  Using the identifications $E_{\gf_0(x)} = (\gf_0^* E)_x$ and $E_{\gy_0(x)} = (\gy_0^* E)_x$, this defines a $\mc{H}$-bundle map $\gf_0^* E \to \gy_0^* E$.
  The inverse $\gt^{-1}$ provides an inverse to this $\mc{H}$-bundle map.
\end{proofsketch}

\begin{rmk}
  The above constructions work equally well for principal $\mc{G}$-bundles; we will mainly consider vector $\mc{G}$-bundles from now on.
\end{rmk}

Consider two vector $\mc{G}$-bundles $E,F \to \mc{G}$.
As vector bundles over the manifold $G_0$, we can consider the Whitney sum $E \oplus F \to G_0$ and tensor product $E \otimes F \to G_0$.
These are vector bundles with fibres of the form $(E \oplus F)_x \cong E_x \oplus F_x$ and $(E \otimes F)_x \cong E_x \otimes F_x$.
The $\mc{G}$-actions on $E$ and $F$ lead to a $\mc{G}$-actions on the bundles $E \oplus F$ and $E \otimes F$, endowing them with a vector $\mc{G}$-bundle structures.
This is done in the obvious way: for a vector $(v,w) \in E \oplus F$ and an arrow $\map{\gg}{\gp((v,w))}{y}$ in $\mc{G}$, we define
\begin{equation*}
  \gg (v,w) := (\gg v,\gg w) \in (E \oplus F)_y.
\end{equation*}
Since the actions of $\mc{G}$ on $E$ and $F$ are smooth, this action is also smooth.
We define the action of $\mc{G}$ on $E \otimes F$ analogously.
Thus the operations $\oplus$ and $\otimes$ are defined on the set of vector $\mc{G}$-bundles.
As one would expect, these operations behave well with respect to $\mc{G}$-bundle isomorphism.
To see this, let $E^\prime,F^\prime \to \mc{G}$ be two more vector $\mc{G}$-bundles, and suppose $\map{f}{E}{E^\prime}$ and $\map{g}{F}{F^\prime}$ are $\mc{G}$-bundle isomorphisms.
Then in particular $f$ and $g$ are isomorphisms of vector bundles over $G_0$, and the map
\begin{equation*}
  \map{f \oplus g}{E \oplus F}{E^\prime \oplus F^\prime}
\end{equation*}
sending a vector $(v,w)$ to $(f(v),g(w))$ is a vector bundle isomorphism (see \cite[\textsection 1.2]{mA67}).
Since $\mc{G}$ acts componentwise on sums and since $f$ and $g$ are $\mc{G}$-maps, $f \oplus g$ is also a $\mc{G}$-map, and hence gives a $\mc{G}$-bundle isomorphism $E \oplus F \to E^\prime \oplus F^\prime$.
Likewise, $E \otimes F$ and $E^\prime \otimes F^\prime$ are isomorphic as $\mc{G}$-bundles.
Therefore we can define the semiring $\Vect_\FF(\mc{G})$, whose elements are isomorphism classes of vector $\mc{G}$-bundles over the field $\FF$ (either $\RR$ or $\CC$), with the operations $\oplus$ and $\otimes$.\footnote{We have only shown that these operations are well-defined. The other axioms, such as associativity and distributivity, have been neglected. The proof that the remaining axiom holds is identical to the proof in the case of vector bundles over topological spaces.}
In Section \ref{gpoidk} we will use this semiring (with $\FF = \CC$) to construct the $K$-theory of the Lie groupoid $\mc{G}$.

If $\map{\gf}{\mc{H}}{\mc{G}}$ is a Lie groupoid homomorphism and if $E,F \to \mc{G}$ are two isomorphic $\mc{G}$-bundles, then the pullbacks $\gf^*E$ and $\gf^*F$ are isomorphic (for the same reason as in the case of vector bundles over manifolds).
In fact, if we view $\gf$ as a morphism in the quotient category $\cat{LGpoid}^\prime$, Lemma \ref{canism} shows that $\gf^*E$ and $\gf^*F$ remain isomorphic.
Hence a morphism $\map{\gf}{\mc{H}}{\mc{G}}$ in $\cat{LGpoid}^\prime$ induces a map
\begin{equation*}
  \map{\gf^*}{\Vect_\FF(\mc{G})}{\Vect_\FF(\mc{H})}.
\end{equation*}
This map is easily shown to be a semiring homomorphism.
We are more interested in the situation in the category $\cat{LGpoid}^\prime[\gS^{-1}]$; more precisely, we would like to pull back bundles by generalised maps.
As such, let $E \to \mc{G}$ be a vector $\mc{G}$-bundle, and suppose
\begin{equation}\label{gm}
  \xymatrix{
    \mc{H} & \mc{K} \ar[l]_\ge \ar[r]^\gf & \mc{G}
  }
\end{equation}
is a generalised map.
We can pull $E$ back to a vector bundle $\gf^* E$ over $\mc{K}$; to get a bundle over $\mc{H}$, we must then push $\gf^*E$ \emph{forward} by $\ge$.
Here we encounter a problem: if $\ge$ is a weak equivalence with no further assumptions, then it is not at all obvious how to proceed.
however, if $\ge$ is a \emph{strong} equivalence, then there exists a weak inverse $\ge^{-1}$ to $\ge$ (unique up to natural isomorphism), and we can pull $\gf^* E$ back by $\ge^{-1}$ to get a bundle over $\mc{H}$, which we denote $(\ge,\gf)^*E$.
Of course, this bundle is only defined up to canonical isomorphism.
Nevertheless we get a well-defined semiring homomorphism
\begin{equation*}
  \map{(\ge,\gf)^*}{\Vect_\FF(\mc{G})}{\Vect_\FF(\mc{H})}
\end{equation*}
under the assumption that $\ge$ is a strong equivalence.

We must then show that this homomorphism is well-defined on equivalence classes of generalised maps.
In attempting to prove this, we encounter another problem.
We must show that if we have a vector bundle $E \to \mc{G}$ and a diagram in $\cat{LGpoid}^\prime$ of the form
\begin{equation*}
  \xymatrix{
    & \mc{K} \ar[dl]_\ge \ar[dr]^\gf & \\
    \mc{H} & \mc{M} \ar[l]^\ga \ar[r]_\gb \ar[u]^a & \mc{G}
  }
\end{equation*}
such that $\ga$ is a weak equivalence, then $(\ge,\gf)^* E \cong (\ga,\gb)^* E$.
But without assuming that $\ga$ is a strong equivalence, this is meaningless.
However, if $\ga$ is assumed to be a strong equivalence, then we can take a weak inverse $\ga^{-1}$ to $\ga$ and write
\begin{equation*}
  (\ge^{-1})^* \gf^* E \cong (a \circ \ga^{-1})^* \gf^* E \cong (\ga^{-1})^* a^* \gf^* E \cong (\ga^{-1})^* \gb^* E,
\end{equation*}
so that $(\ge,\gf)^* E \cong (\ga,\gb)^* E$.
Hence we have a semiring homomorphism
\begin{equation*}
  \map{[\ge,\gf]^*}{\Vect_\FF(\mc{G})}{\Vect_\FF(\mc{H})}.
\end{equation*}
It follows that the isomorphism class of $\Vect_\FF(\mc{G})$ is a strong Morita invariant, in the sense that if $\mc{G}$ and $\mc{H}$ are strongly Morita equivalent, then $\Vect_\FF(\mc{G}) \cong \Vect_\FF(\mc{H})$

If in the generalised maps above all groupoids are assumed to be \'etale, then we can avoid the issue of having to assume our weak equivalences are strong equivalences.
Given a generalised map $(\ge,\gf)$ as in \eqref{gm}, since $\mc{H}$ and $\mc{M}$ are assumed to be \'etale, by Remark \ref{smorita} we can replace $\ge$ with a strong Morita equivalence, resulting in a diagram
\begin{equation*}
  \xymatrix@C=5pt@R=5pt{
    & \mc{K}^\prime \ar[dl]_{\ge^\prime} \ar[dr]^\gy & & & \\
    \mc{H} & & \mc{K} \ar[rr]^\gf & & \mc{G} \\
  }
\end{equation*}
in which $\ge^\prime$ is a strong equivalence.
Thus we are reduced to the case where we have a strong equivalence on the left, and we can define
\begin{equation*}
  (\ge,\gf)^*E := ((\ge^\prime)^{-1})^* \gy^* \gf^* E.
\end{equation*}
However, this doesn't solve the problem when only $\mc{H}$ and $\mc{G}$ are assumed to be \'etale, as is the case when we consider maps between orbifolds.

\section[Associated $\mc{G}$-bundles and H-S morphisms]{Associated $\mc{G}$-bundles and Hilsum-Skandalis morphisms}\label{HS}
In the previous section, we showed that (isomorphism classes of) vector $\mc{G}$-bundles were invariant under strong Morita equivalence.
To show this, we used the fact that bundles easily push forward under strong equivalences; it is not at all obvious that they can be pushed forward under \emph{weak} equivalences.
This prevents us from showing the Morita invariance of vector $\mc{G}$-bundles.
In fact, this Morita invariance \emph{does} hold, but it is not easily seen without an alternative characterisation of vector $\mc{G}$-bundles, and this alternative characterisation is not obvious without an alternative characterisation of generalised morphisms.
The required characterisation is that of Hilsum-Skandalis (see \cite[Chapter II]{jM96} and \cite{HS87}; we shall not state the definition, but will instead list one consequence of which we will make use.

\begin{prop}
  Let $H$ be a Lie group and $\mc{G}$ a Lie groupoid.
  Then the set of isomorphism classes of $H$-principal $\mc{G}$-bundles $P \to \mc{G}$ is in natural bijective correspondence with the set
  \begin{equation*}
    \cat{LGpoid}^\prime[\gS^{-1}](\mc{G},H)
  \end{equation*}
  of generalised morphisms from $\mc{G}$ to the Lie group $H$ (viewed as a groupoid with one object as in Example \ref{actiongpoid}).
\end{prop}

Using this fact, we can show that isomorphism classes of rank $n$ vector bundles over $\mc{G}$ are the same as generalised morphisms $\mc{G} \to GL(n,\RR)$.
Once this has been shown, if $\mc{H}$ is Morita equivalent to $\mc{G}$, then $\mc{H}$ and $\mc{G}$ are isomorphic in $\cat{LGpoid}^\prime[\gS^{-1}]$ (by Proposition \ref{meisom}) and so the morphism sets
\begin{equation*}
  \cat{LGpoid}^\prime[\gS^{-1}](\mc{G},\GL(n,\RR)) \approx \cat{LGpoid}^\prime[\gS^{-1}](\mc{H},\GL(n,\RR))
\end{equation*}
can be naturally identified.
This identification can be shown to respect Whitney sums and tensor products, and so it will follow that the isomorphism class of the semiring $\Vect_\FF(\mc{G})$ is Morita invariant.

To achieve this, we need to show that rank $n$ vector $\mc{G}$-bundles (up to isomorphism) are the same as $\GL(n,\RR)$-principal $\mc{G}$-bundles.
This is done by generalising the associated bundle construction to the case of Lie groupoids.
Recall that to a $\GL(n,\RR)$-principal bundle $\map{\gp}{P}{G_0}$ over $G_0$ we can associate a vector bundle $V \to G_0$, with the same transition functions as $P$, by setting
\begin{equation*}
  V := \GL(n,\RR) \sm (P \times \RR^n),
\end{equation*}
where $A \in \GL(n,\RR)$ acts on $P \times \RR^n$ by
\begin{equation*}
  A(p,v) := (p \cdot A, A^{-1} v).
\end{equation*}
If $P$ is a $\mc{G}$-space, then we can make $V$ into a $\mc{G}$-space by setting
\begin{equation*}
  \gg [p,v] = [\gg p,v]
\end{equation*}
whenever $\gg$ is an arrow in $\mc{G}$ with $s(\gg) = \gp(p)$.
Since the actions of $\GL(n,\RR)$ and $\mc{G}$ on $P$ are compatible, this action is well-defined.
Therefore the associated bundle construction carries over to the world of Lie groupoids, and as with the case of bundles over topological spaces, this shows that rank $n$ vector $\mc{G}$-bundles are given by $\GL(n,\RR)$-principal $\mc{G}$-bundles (and conversely; all up to isomorphism of course).
By the previous paragraph, we can conclude that the set of isomorphism classes of vector $\mc{G}$-bundles over a Lie groupoid is Morita invariant.

To convince the reader that the semiring $\Vect_\FF(\mc{G})$ is Morita invariant, we provide a sketch of the characterisation of the Whitney sum in terms of $\GL(n,\RR)$-principal $\mc{G}$-bundles.
Suppose we have two generalised maps
\begin{equation*}
  \mc{G} \to \GL(n,\RR), \quad \mc{G} \to \GL(m,\RR)
\end{equation*}
with the same mediating groupoid, i.e. a diagram of the form
\begin{equation*}
  \xymatrix@R=5pt{
    & & \GL(n,\RR) \\
    \mc{G} & \mc{H} \ar[l]_\ge \ar[ur]^(.36)\gf \ar[dr]_(.36)\gy & & \\
    & & \GL(m,\RR).
  }
\end{equation*}
Then we can construct a Lie groupoid homomorphism $\map{\gf \oplus \gy}{\mc{H}}{\GL(n+m,\RR)}$ by sending each arrow $h \in H_1$ to the block matrix
\begin{equation*}
  \gf_1(h) \oplus \gy_1(h) = \twmat{\gf_1(h)}{0}{0}{\gf_2(h)}.
\end{equation*}
This homomorphism corresponds to the construction of the Whitney sum of two vector bundles.
Likewise we can construct a homomorphism $\map{\gf \otimes \gy}{\mc{H}}{\GL(nm,\RR)}$ corresponding to the tensor product of vector bundles.\footnote{This is by no means a proof that $\Vect_\FF(\mc{G})$ is Morita invariant as a semiring. Unfortunately time constraints have prevented the inclusion of such a proof.}
  
\section{$K$-theory of Lie groupoids and orbifolds}\label{gpoidk}

Recall that the $K$-theory $K^0(M)$ of a compact manifold $M$ is defined to be the Grothendieck group of the commutative monoid $\Vect_\CC(M)$ (see \cite[Chapter II]{mA67}), with a ring structure coming from the multiplicative structure on $\Vect_\CC(M)$.
Likewise we can define the $K$-theory of a compact Lie groupoid $\mc{G}$: $\Vect_\CC(\mc{G})$ is a commutative monoid (forgetting the multiplicative structure), and we define $K^0(\mc{G})$ to be the Grothendieck group of this monoid.
The multiplication on $\Vect_\CC(\mc{G})$ then induces a ring structure on $K^0(\mc{G})$.
Of course, the ring $K^0(\mc{G})$ is only defined up to canonical isomorphism.
Since the semiring $\Vect_\CC(\mc{G})$ is Morita invariant, $K^0(\mc{G})$ is also Morita invariant.
Also, each generalised map $\map{(\ge,\gf)}{\mc{G}}{\mc{H}}$ between two compact Lie groupoids induces a ring homomorphism
\begin{equation*}
  \map{(\ge,\gf)^*}{K^0(\mc{H})}{K^0(\mc{G})}
\end{equation*}
thanks to the universal property of the Grothendieck group.

If $\mc{G}$ is a Lie groupoid (not necessarily compact), then we define a \emph{vector $\mc{G}$-bundle with compact support} to be a pair of $\mc{G}$-bundles $E^0,E^1 \to \mc{G}$ along with a bundle morphism $\map{\gs}{E^0}{E^1}$ such that the set
\begin{equation*}
  \{x \in G_0 \mid \text{$\map{\gs_x}{E_x^0}{E_x^1}$ is not an isomorphism}\}
\end{equation*}
defines a compact set in $|\mc{G}|$ \cite[Definition 2.4]{HW12}.
The \emph{compactly supported $K$-theory} ring $K_c^0(\mc{G})$ is then defined analogously to the $K$-theory ring in the compact case, and likewise every generalised map $\map{(\ge,\gf)}{\mc{G}}{\mc{H}}$ between two Lie groupoids induces a ring homomorphism $K_c^0(\mc{H}) \to K_c^0(\mc{G})$.
If $\mc{G}$ is compact then $K_c^0(\mc{G})$ and $K^0(\mc{G})$ are isomorphic, so in general (i.e. when $\mc{G}$ is not necessarily compact) we will write $K^0(\mc{G})$ instead of $K_c^0(\mc{G})$.

For a possibly non-compact Lie groupoid $\mc{G}$, let $\mc{G} \times \RR^n$ denote the groupoid with $(\mc{G} \times \RR^n)_0 := G_0 \times \RR^n$, $(\mc{G} \times \RR^n)_1 := G_1 \times \RR^n$ and with source and target maps $(s,\id)$ and $(t,\id)$ respectively.
The orbit space of $\mc{G} \times \RR^n$ is homeomorphic to $|\mc{G}| \times \RR^n$, so $\mc{G} \times \RR^n$ is not compact.
There is an obvious Lie groupoid homomorphism $\map{\gi}{\mc{G}}{\mc{G} \times \RR^n}$: define
\begin{equation*}
  \gi_0 := \id_{G_0} \times c^0, \quad \gi_1 = \id_{G_1} \times c^1
\end{equation*}
where $\map{c^0}{G_0}{\RR^n}$ (resp. $\map{c^1}{G_1}{\RR^n}$) sends each point in $G_0$ (resp. $G_1$) to the origin $0 \in \RR^n$.
The homomorphism $\gi$ (or indeed any generalised map equivalent to $\gi$) induces a ring homomorphism
\begin{equation*}
  \map{\gi^*}{K^0(\mc{G} \times \RR^n)}{K^0(\mc{G})}.
\end{equation*}
We define $K^{-n}(\mc{G}) := \ker \gi^*$ for each natural number $n$.
If $\mc{G}$ and $\mc{H}$ are Morita equivalent Lie groupoids, then $\mc{G} \times \RR^n$ and $\mc{H} \times \RR^n$ are Morita equivalent,\footnote{To see this, one needs to show that if $\map{\ge}{\mc{G}}{\mc{H}}$ is a weak equivalence, then so is the homomorphism $\map{(\ge,\id)}{\mc{G} \times \RR^n}{\mc{H} \times \RR^n}$, where $(\ge,\id)_0 = (\ge_0,\id_{\RR^n})$ and $(\ge,\id)_1 = (\ge_1,\id_{\RR^n})$. One can then construct a Morita equivalence between $\mc{G} \times \RR^n$ and $\mc{H} \times \RR^n$ by means of a Morita equivalence between $\mc{G}$ and $\mc{H}$.} so $K^{-n}(\mc{G})$ is Morita invariant for each $n$.
Using the above definitions, we can define the rings $K^{-n}(\mf{X})$ for any orbifold $\mf{X}$: we simply define $K^{-n}(\mf{X})$ to be $K^{-n}(\mc{G})$, where $\mc{G}$ is a representing groupoid for $\mf{X}$.
Clearly any smooth map $\map{f}{\mf{X}}{\mf{Y}}$ between two orbifolds induces a corresponding ring homomorphism $\map{f^*}{K^{-n}(\mf{Y})}{K^{-n}(\mf{X})}$ in $K$-theory.

At this point in the development of topological $K$-theory we would prove the Bott periodicity theorem, implying that $K^{-n} \cong K^{-n-2}$ for each $n$, and thus prompting us to treat $K$-theory as a $\ZZ_2$-graded theory.
However, the proof of Bott periodicity does not immediately generalise to the much wider scope of Lie groupoids.
In fact, it does not even generalise immediately to arbitrary orbifolds, although when we restrict our attention to presentable orbifold groupoids Bott periodicity will hold automatically.
We will leave this to Section \ref{presentablek}.

\begin{rmk}
  Bott periodicity does hold for proper Lie groupoids.
  This can be proven by identifying the $K$-theory of $\mc{G}$ with the $K$-theory of the reduced $C^*$-algebra $C_\text{red}^*(\mc{G})$ as in \cite[Chapter 2]{aC94} and using that Bott periodicity holds in this context \cite{HG04}.
\end{rmk}


\chapter[Geometry and $K$-theory]{Geometry and $K$-theory of orbifolds}\label{dg}

In this chapter, we look at the differential geometry and `topological' $K$-theory of orbifolds using the groupoid language developed in Chapter \ref{gpoids}.
From the viewpoint of differential geometry, orbifolds are not too different from manifolds: many geometric constructions on smooth manifolds can be generalised to orbifolds.
That this generalisation is possible is a consequence of the existence of partitions of unity subordinate to open covers of orbifolds.
As in the theory of smooth manifolds, this allows us to construct globally-defined objects from locally-defined objects.
This result is proven in Section \ref{po1s}, and the following sections detail some of its consequences.
In the first half of this chapter it becomes clear that the groupoid language does not tell us anything new about the differential geometry of orbifolds---this can be studied perfectly adequately using the language of orbifold atlases.
However, we will see in the second half that the groupoid approach is very useful in studying the $K$-theory of presentable orbifolds.
The inertia orbifold will be introduced in Section \ref{inof} (this is most easily done using groupoids) and then appear as the natural target for the `delocalised Chern character'.
This generalises the usual Chern character for smooth manifolds, and defines an isomorphism between the $K$-theory of a presentable orbifold (modulo torsion) and the de Rham cohomology of the corresponding inertia orbifold.
This isomorphism is the subject of Section \ref{dcci}.

\section{Orbifold bundles revisited}\label{VBs}

\begin{dfn}
  Let $\mf{X}$ be an orbifold with groupoid representation $|\mc{G}| \to X$.
  A \emph{vector} (resp. \emph{principal}) bundle over $\mf{X}$ is given by a vector (resp. principal) $\mc{G}$-bundle $E \to \mc{G}$.
\end{dfn}

Sections of bundles over orbifolds are defined by sections of the corresponding $\mc{G}$-bundles.

\begin{rmk}
  With $\mf{X}$ and $E$ as in the above definition, if $|\mc{G}^\prime| \to X$ is another orbifold presentation of $\mf{X}$ and
  \begin{equation*}
    \xymatrix{
      \mc{G} & \mc{M} \ar[l]_\ge \ar[r]^\gf & \mc{G}^\prime
    }
  \end{equation*}
  is a Morita equivalence defining an equivalence of orbifold structures, then $(\ge,\gf)$ and $E$ induce a well-defined \emph{isomorphism class} $(\ge,\gf)_*(E)$ of bundles over $\mc{G}^\prime$.
  If $(\ge,\gf)$ is a strong Morita equivalence, then $(\ge,\gf)_*(E) := \gf_* \ge* E$ is well-defined up to canonical isomorphism (depending on the choice of weak inverse for $\gf$ used to define $\gf_*$ as in Section \ref{gbundles}.
  Otherwise, we can use the techniques of Section \ref{HS} to define $(\ge,\gf)_* E$ as an isomorphism class of bundles.
  This shows that the notion of a \emph{bundle} over an orbifold---as opposed to the notion of an \emph{isomorphism class of bundles}---is not well-defined if we represent orbifolds with groupoids.
  This is an important flaw of the groupoid approach.
  To get around this, whenever we have to deal with an actual bundle rather than an isomorphism class, we will work with a specific groupoid representation and not consider equivalent presentations.
\end{rmk}

Many constructions in differential geometry---in particular, the notion of a connection on a vector bundle---rely on the existence of nonzero sections of such a vector bundle.
Of course, a vector bundle on a manifold has plenty of nonzero sections,\footnote{The reader is warned not to confuse `nonzero' sections with `nonvanishing' sections.} so we take this requirement for granted.
For a vector bundle over a Lie groupoid this is no longer the case, so we are led to restrict the class of vector bundles under consideration.
Recall that for an \'etale groupoid $\mc{G}$ and a point $x \in G_0$, there is a natural group homomorphism $\map{\gl_x}{\mc{G}_x}{\Diff(x,x)}$ (see Remark \ref{germ}).
If $E \to \mc{G}$ is a vector $\mc{G}$-bundle, then the isotropy group $\mc{G}_x$ also acts on the fibre $E_x$, thus defining a group homomorphism $\map{\gr_x}{\mc{G}_x}{\End{E_x}}$.

\begin{dfn}
  The vector $\mc{G}$-bundle $E \to \mc{G}$ is called \emph{good} if
  $\ker \gl_x \subset \ker \gr_x$
  for each $x \in G_0$.
  The $\mc{G}$-bundle $E$ is called \emph{bad} if it is not good.
\end{dfn}

For a good $\mc{G}$-bundle, if a loop $\map{g}{x}{x}$ acts trivially about $x$ (as a local diffeomorphism), then it acts trivially on the fibre $E_x$.
Note that if the \'etale groupoid $\mc{G}$ is effective, then every vector bundle over $\mc{G}$ is automatically good, since in this case $\ker \gl_x = \{\id_x\}$ for all $x$.

\begin{example}\label{badbundle}\textbf{(A bad bundle over an ineffective orbifold.)}
  Let $\ZZ_2$ act trivially on an arbitrary smooth manifold $M$, and let $\mf{X}$ denote the corresponding ineffective orbifold (with representing groupoid $\ZZ_2 \ltimes M$).
  Let $\ZZ_2$ act on the trivial bundle $M \times \RR$ by
  \begin{equation*}
    \bar{1}(x,v) := (x,-v)
  \end{equation*}
  where $\bar{1}$ denotes the non-identity element of the additive group $\ZZ_2$.
  Equipped with this $\ZZ_2$-action, $M \times \RR$ is an bad orbifold vector bundle over $\mf{X}$.
  Now suppose $\gs$ is a section of this bundle over $\mf{X}$.
  For all $x \in M = |\mf{X}|$, we must have
  \begin{equation*}
    \gs(x) = \gs(\bar{1}x) = \bar{1}\gs(x) = -\gs(x),
  \end{equation*}
  from which it follows that the only section of $M \times \RR$ over $\mf{X}$ is the zero section.
\end{example}

We refer to a vector bundle over an orbifold $\mf{X}$, given by a vector $\mc{G}$-bundle $E$ over a representing groupoid $\mc{G}$, as \emph{good} or \emph{bad} if $E$ is good or bad respectively.

\'Etale groupoids share some of the geometric properties of smooth manifolds.
In particular, the notion of tangent bundle carries over to \'etale groupoids.
Given an \'etale groupoid $\mc{G}$, we can consider the tangent bundle $TG_0$ as a $\mc{G}$-space in the following way.
Let $\map{\gg}{x}{y}$ be an arrow in $\mc{G}$.
Then $\gg$ determines the germ of a local diffeomorphism from $x$ to $y$; the differential of such a germ\footnote{To define the differential of the germ of a local diffeomorphism, simply observe that a representative $\map{f}{U}{V}$ of such a germ induces a diffeomorphism $\map{Df}{T_x G_0}{T_y G_0}$, and that this map is independent of our choice of representative by the definition of a germ.} gives an isomorphism of tangent spaces $T_x G_0 \to T_y G_0$.
Thus the action of $\mc{G}$ on $G_0$ lifts to an action upon $TG_0$, making $TG_0$ into a $\mc{G}$-space.
We denote the resulting vector $\mc{G}$-bundle by $T\mc{G}$; this is a good bundle.
For an orbifold $\mf{X}$ with representing groupoid $\mc{G}$, we define the tangent bundle $T\mf{X}$ to be the vector bundle over $\mf{X}$ given by $T\mc{G}$.
When $\mc{G}$ is an orbifold atlas groupoid, this bundle coincides with the tangent bundle constructed in Section \ref{oftb}.

A local diffeomorphism as above induces maps not only between tangent spaces $T_x G_0 \to T_y G_0$, but also between all tensor bundles $(T_s^r)_x G_0 \to (T_s^r)_y G_0$ (see \cite[Example 5.4]{KN63} and Section \ref{vapb}); thus to an \'etale groupoid $\mc{G}$ we can associate tensor bundles $T_s^r \mc{G}$.
Likewise we can define symmetric bundles $S^r T\mc{G}$ and exterior bundles $\wedge^r T\mc{G}$.
As with the tangent bundle, these are all good bundles. 
In particular, we can construct the exterior powers $\wedge^r T^* \mc{G}$ of the cotangent bundle of $\mc{G}$, which allow us to consider differential forms on the \'etale groupoid $\mc{G}$.
When $\mc{G}$ is a represnting groupoid for an orbifold $\mf{X}$, we obtain the corresponding notions for orbifolds.

\section{Partitions of unity}\label{po1s}

The existence of partitions of unity (see \cite[\textsection II.3]{sL99}) on a smooth manifold $M$ allows us to extend local constructions---such as Riemannian metrics, and connections on bundles---to global constructions.
They are also crucial in proving the exactness of the generalised Mayer-Vietoris sequence \cite[Proposition 8.5]{rBlT82}, from which follows the \v{C}ech-de Rham isomorphism \cite[Theorem 8.9]{rBlT82}.
The same is true for orbifolds: partitions of unity exist, they can be used to prove the existence of Riemannian metrics on orbifolds and connections on bundles, and from their existence follows a \v{C}ech-de Rham type isomorphism for orbifold de Rham cohomology.
We will define these concepts and prove the aforementioned results in the following sections.

\begin{dfn}\label{po1}
  Let $\mf{X}$ be an orbifold, and let $\mc{U} = \{U_\ga\}_{\ga \in A}$ be an open cover of $X = |\mf{X}|$.
  A \emph{partition of unity} subordinate to the cover $\mc{U}$ is given by a locally finite refinement $\{V_\gb\}_{\gb \in B} \to \mc{U}$ along with a collection of smooth non-negative functions $\{\gr_\gb\}_{\gb \in B} \subset C^\infty(\mf{X})$ such that
  \begin{enumerate}[(a)]
  \item
    for each $\gb$, $\supp \gr_\gb$ is compact and contained in $V_\gb$, and
  \item
    for each $x \in X$,
    \begin{equation*}
      \sum_{\gb} \gr_\gb(x) = 1.
    \end{equation*}
  \end{enumerate}
\end{dfn}

Of course, to speak of smooth functions in a truly well-defined way we must choose a groupoid representation of $\mf{X}$.
There is no harm in doing this.

\begin{thm}\label{po1exists}
  Let $\mf{X}$ be an orbifold with groupoid representation $|\mc{G}| \to X$, and let $\mc{U}$ be an open cover of $X$ as in Definition \ref{po1}.
  Then there exists a partition of unity subordinate to $\mc{U}$.
\end{thm}

\begin{proof}
  Use paracompactness of $X$ and Theorem \ref{localaction} to find a locally finite refinement $\map{\gl}{\{U_\ga\}_{\ga \in A}}{\mc{U}}$ of $\mc{U}$, with $U_\ga = \gp(\wtd{U}_\ga)$ where $\{\wtd{U}_\ga\}_{\ga \in A}$ is a refinement of the cover $\gp^{-1} \mc{U}$ of $G_0$ such that for each $\ga$,
  \begin{equation}\label{isisis}
    \mc{G}|_{\wtd{U}_\ga} \cong \mc{G}_{x(\ga)} \ltimes \wtd{U}_\ga
  \end{equation}
  for some $x(\ga) \in U_\ga$ (here $\gp$ denotes the map $\mc{G} \to |\mc{G}| \to X$).
  We can also find for each $\ga$ a relatively compact set $W_\ga$ with
  \begin{equation*}
    \overline{\wtd{W}_\ga} \subset \wtd{U}_\ga \qquad (\wtd{W}_\ga := \gp^{-1}(W_\ga))
  \end{equation*}
  such that $\{W_\ga\}_{\ga \in A}$ is a locally finite refinement of $\mc{U}$ (see \cite[Propositions II.3.1 and II.3.2]{sL99}).

  For each $\ga$, choose a compactly supported non-negative smooth function $\gy_\ga \in C^\infty(G_0)$ such that
  \begin{equation*}
    \text{$\gy_\ga \equiv 1$ on $\wtd{W}_\ga$}, \quad \text{$\gy_\ga \equiv 0$ on $\wtd{U}_\ga^c$.}
  \end{equation*}
  Now define functions $\gy_\ga^* \in C^\infty(G_0)$ by
  \begin{equation*}
    \gy_\ga^*(y) := \left\{ \begin{array}{ll} 0 & y \notin \wtd{U}_\ga, \\ \left| \mc{G}_{x(\ga)} \right|^{-1} \sum_\gg \gy_\ga(\gg y) & x \in \wtd{U}_\ga, \end{array} \right.
  \end{equation*}
  where the sum is over all $\gg \in G_1$ such that $s(\gg) = y$ and $t(\gg) \in \wtd{U}_\ga$.
  Since $\wtd{U}_\ga = \gp^{-1}(U_\ga)$, if $\gg$ is in $G_1$, then the source and target of $\gg$ are either both in $\wtd{U}_\ga$ or both not in $U_\ga$.
  Furthermore, by the isomorphism \eqref{isisis}, this sum is identified with the average
  \begin{equation*}
    \frac{1}{|\mc{G}_{x(\ga)}|} \sum_{\gg \in \mc{G}_{x(\ga)}} \gy_\ga(\gg y)
  \end{equation*}
  over the group $\mc{G}_{x(\ga)}$.
  It follows that $\gy_\ga^*$ is a well-defined smooth function in $C^\infty(\mc{G})$, and thus defines a smooth function on $\mf{X}$.

  Finally, as with the usual construction of a partition of unity on a manifold, for each $\ga$ define
  \begin{equation*}
    \gr_\ga := \frac{\gy_\ga^*}{\sum_\gh \gy_\gh^*}.
  \end{equation*}
  It is now easy to check that this collection of functions defines a partition of unity subordinate to $\mc{U}$.
\end{proof}

The above proof is adapted from \cite[Proposition 1.2]{yC90}.
Observe that this is essentially an `orbifold atlas' proof translated into groupoid language via Theorem \ref{localaction}; properties of groupoids were not used in any profound way.

\section{Differential forms and de Rham cohomology}

\begin{dfn}
  A \emph{differential $n$-form} $\go$ on an \'etale groupoid $\mc{G}$ is a section of the $n^\text{th}$ exterior power of the cotangent $\mc{G}$-bundle
  \begin{equation*}
    \go \in \gG\left(\mc{G}, \wedge^n T^* \mc{G}\right).
  \end{equation*}
  Equivalently, a differential $n$-form on $\mc{G}$ is given by a differential $n$-form $\go$ on $G_0$ such that for all arrows $\map{\gg}{x}{y}$ in $\mc{G}$,
  \begin{equation*}
    \go(y) = \gg \cdot \go(x).
  \end{equation*}
\end{dfn}

A differential form on an orbifold $\mf{X}$ is defined as a section of an exterior power of the cotangent bundle $T^* \mf{X}$, or equivalently as a differential form on a representing groupoid $\mc{G}$.
As with any concept involving sections of vector bundles, the concept of a differential form on an orbifold $\mf{X}$ is somewhat hazy without considering a particular groupoid representation of $\mf{X}$.
The space of differential $n$-forms over a groupoid $\mc{G}$ (resp. an orbifold $\mf{X}$) is denoted by $\gO^n(\mc{G})$ (resp. $\gO^n(\mf{X})$).

\begin{prop}
  Let $\mc{G}$ be an \'etale groupoid, and suppose $\go \in \gO^n(G_0)$ is a $\mc{G}$-equivariant $n$-form on $G_0$.
  Then $d\go \in \gO^{n+1}(G_0)$ is also $\mc{G}$-equivariant.
\end{prop}

\begin{proof}
  Suppose $\map{\gg}{x}{y}$ is an arrow in $\mc{G}$.
  Then $\gg$ extends to a local diffeomorphism $\map{\wtd{\gg}}{U_x}{U_y}$ between neighbourhoods of $x$ and $y$.
  Using $D\wtd{\gg}$ to denote the action of $\mc{G}$ on $\wedge^n T^* G_0$, $\mc{G}$-equivariance of $\go$ on the aforementioned neighbourhoods can be expressed by the commutative diagram
  \begin{equation*}
    \xymatrix{
      U_x \ar[d]_\go \ar[r]^{\wtd{\gg}} & U_y \ar[d]^\go \\
      \bigwedge^n T^* U_x \ar[r]_{D\wtd{\gg}} & \bigwedge^n T^* U_y.
    }
  \end{equation*}
  By naturality of the exterior derivative, this extends to a commutative diagram
  \begin{equation*}
    \xymatrix{
      U_x \ar[d]_\go \ar[r]^{\wtd{\gg}} & U_y \ar[d]^\go \\
      \bigwedge^n T^* U_x \ar[d]_d \ar[r]_{D\wtd{\gg}} & \bigwedge^n \ar[d]^d T^* U_y \\
      \bigwedge^{n+1} T^* U_x \ar[r]_{D\wtd{\gg}} & \bigwedge^{n+1} T^* U_y
    }
  \end{equation*}
  which expresses the local $\mc{G}$-equivariance of $d\go$.
  This implies global $\mc{G}$-equivariance as desired.
\end{proof}

\begin{cor}
  The exterior derivative is defined on $\gO^\bullet(\mc{G}) := \oplus_{n \geq 0} \gO^n(\mc{G})$.
\end{cor}

Likewise the wedge product $\wedge$ is defined on $\gO^\bullet(\mc{G})$.
Thus we are able to define the de Rham complex $(\gO^\bullet(\mc{G}),d)$ and cohomology rings $H_\text{dR}^\bullet(\mc{G})$ of an \'etale groupoid, as well as the corresponding complex $(\gO^\bullet(\mf{X}),d)$ and cohomology rings $H_\text{dR}^\bullet(\mf{X})$ of an orbifold.
However, $H_\text{dR}^\bullet(\mf{X})$ is not a very strong `orbifold invariant', as expressed by the following theorem.

\begin{thm}\label{cdrs} \textbf{(\v{C}ech-de Rham-Satake)}
  Let $\mf{X}$ be a presentable orbifold.
  Then the de Rham cohomology ring of $\mf{X}$ is isomorphic to the singular cohomology ring of $|\mf{X}|$ with real coefficients:
  \begin{equation}\label{cdrse}
    H_\text{dR}^\bullet(\mf{X}) \cong H^\bullet(|\mf{X}|;\RR).
  \end{equation}
\end{thm}

To prove that there is a \emph{group} isomorphism  of the form \eqref{cdrse} we shall consider the sheaf of locally-defined differential forms over an orbifold $\mf{X}$.\footnote{For a nice exposition of sheaf theory, see \cite[Chapter II]{rW08}.}
Let $\map{f}{|\mc{G}|}{X}$ be a groupoid representation of $\mf{X}$, and consider the correspondence $\gO^p$ sending an open set $U \subset X$ to the group
\begin{equation*}
  \gO^p\left(\mc{G}|_{(f \circ \gp)^{-1}(U)}\right)
\end{equation*}
of differential forms on the restriction of $\mc{G}$ to $(f \circ \gp)^{-1}(U)$.

\begin{prop}
  The correspondence $\gO^p$ above defines a fine sheaf of Abelian groups on $X$.
\end{prop}

\begin{proof}
  The correspondence is well-defined---the restriction of an orbifold groupoid $\mc{G}$ to an open set is \'etale by Theorem \ref{localaction}, so we can speak of differential forms over restrictions.
  For two open sets $U \subset W \subset X$, we define the restriction map $\map{r_U^W}{\gO^p(W)}{\gO^p(U)}$ to be the usual restriction of differential forms
  \begin{equation*}
    \gO^p(\wtd{W}) \to \gO^p(\wtd{U}),
  \end{equation*}
  writing $\wtd{W} := (f \circ \gp)^{-1}(W)$ and $\wtd{U} := (f \circ \gp)^{-1}(U)$.
  If $\go$ is in $\gO^p(W)$, then $\go$ can be viewed as a $\mc{G}|_{\wtd{W}}$-equivariant differential form on $\wtd{W}$.
  Since the arrows of $\mc{G}|_{\wtd{U}}$ are contained in $\mc{G}|_{\wtd{W}}$, the restriction of $\go$ to $\wtd{U}$ is $\mc{G}|_{\wtd{U}}$-equivariant, so these restriction maps are well-defined.
  It follows immediately that if $U \subset V \subset W \subset X$, then $r_U^V \circ r_V^W = r_U^W$.
  Thus $\gO^p$ is a presheaf.

  The two sheaf axioms \cite[Definition 1.2]{rW08} follow easily.
  To see this, let $\{U_i\}$ be an open cover of a subset $U \subset X$.
  If $\go,\gh \in \gO^p(U)$ agree when restricted to any $U_i$, then they must be equal as differential forms on $G_0$, hence they are equal in $\gO^p(U)$, proving the first axiom.
  As for the second axiom, if we have a form $\go_i \in \gO^p(U_i)$ for each $i$ and if for every nonempty intersection $U_i \cap U_j$ we have that
  \begin{equation*}
    r_{U_i \cap U_j}^{U_i}(\go_i) = r_{U_i \cap U_j}^{U_j}(\go_j),
  \end{equation*}
  then as differential forms on $G_0$ the $\go_i$ `patch together' to define a form
  \begin{equation*}
    \go \in \gO^p(\wtd{U}) \qquad (\wtd{U} := (\gf \circ \gp)^{-1}(U)).
  \end{equation*}
  It remains to show that $\go$ is $\mc{G}|_{\wtd{U}}$-equivariant.
  Suppose $\map{\gg}{x}{y}$ is an arrow in $\mc{G}|_{\wtd{U}}$.
  Then $f \circ \gp(x) = f \circ \gp(y)$, and there is some $U_i$ containing this point.
  Since $\gg$ is an arrow in $\mc{G}|_{\wtd{U}_i}$, the form $\go_i$ is $\gg$-equivariant, and so
  \begin{equation*}
    \go(x) = \go_i(x) = \go_i(y) = \go(y).
  \end{equation*}
  This shows that $\gO^p$ is a sheaf.

  Finally, to show that $\gO^p$ is fine, we simply mimic the proof that the sheaf of differential $p$-forms on a manifold is fine, using the existence of partitions of unity (Theorem \ref{po1exists}).
  Given a locally finite open cover $\{V_\ga\}$ of $X$, take a partition of unity $\{\gr_\ga\}$ subordinate to this cover and define sheaf homomorphisms $\map{\gh_\ga}{\gO^p}{\gO^p}$ by sending a form $\go \in \gO^p(U)$ to the form
  \begin{equation*}
    x \mapsto \gh_\ga(x) \go(x) \in \gO^p(U).
  \end{equation*}
\end{proof}

Consider now the exact sheaf sequence
\begin{equation*}
  0 \to \gO^0 \stackrel{d}{\to} \gO^1 \stackrel{d}{\to} \cdots \to \gO^p \stackrel{d}{\to} \cdots.
\end{equation*}
The kernel of $\map{d}{\gO^0}{\gO^1}$ is precisely the sheaf $\underline{\RR}$ of constant functions $X \to \RR$.\footnote{Viewing a form $\go \in \gO^0(U)$ as a $\mc{G}|_{\wtd{U}}$-equivariant $0$-form (i.e. a $\mc{G}|_{\wtd{U}}$-equivariant function) on $\wtd{U}$, we have that $d\go = 0$ if and only if $\go$ is constant. Constant functions on $\wtd{U}$ are automatically $\mc{G}|_{\wtd{U}}$-equivariant and obviously are no different to constant functions on $U$.}
Since the sheaves $\gO^p$ are fine, we have a fine resolution
\begin{equation*}
  0 \to \underline{\RR} \to \gO^\bullet.
\end{equation*}
The abstract de Rham theorem \cite[Theorem 3.13]{rW08} hence provides an isomorphism
\begin{equation*}
  H^\bullet_\text{dR}(\mf{X}) \cong H^\bullet_\text{\v{C}}(X;\underline{\RR})
\end{equation*}
of the de Rham cohomology groups of $\mf{X}$ with the \v{C}ech cohomology groups of $X$ with coefficients in $\underline{\RR}$.
If $X$ has the homotopy type of a CW-complex, then we have a natural isomorphism
\begin{equation*}
  H^\bullet_\text{\v{C}}(X;\underline{\RR}) \cong H^\bullet(X;\RR)
\end{equation*}
between \v{C}ech and singular cohomology \cite[Theorem 15.8]{rBlT82}.\footnote{Singular cohomology in \cite{rBlT82} has coefficients in $\ZZ$; we obtain our isomorphism by replacing $\ZZ$ with $\RR$ throughout, or alternatively by the universal coefficient theorem \cite[Theorem 15.14]{rBlT82}.}
This yields a group isomorphism of the form \eqref{cdrse} whenever $X$ has the homotopy type of a CW-complex.
In particular, if $\mf{X}$ is presentable, then $X$ is homeomorphic to $G \sm M$ for a smooth almost-free action of a compact Lie group $G$ on a smooth manifold $M$, and such orbit spaces have canonical CW-complex structures (see \cite[Proposition 1.6]{tM71} and \cite[Corollary 7.2]{sI83}).
To prove that \eqref{cdrse} is a ring isomorphism, one must use the `double complex' approach (see \cite[Theorem 14.28]{rBlT82}) rather than the sheaf-theoretic approach.
We will not prove this here, as the proof would take us too far off track.

\begin{rmk}
  The proof of the group isomorphism \eqref{cdrse} given here is a modernised version of that sketched (or at least indicated) by Satake for effective orbifolds in \cite[Theorem 1]{iS56}.
\end{rmk}

\begin{rmk}
  As usual, we can define differential $n$-forms with complex coefficients as sections of the bundle $\wedge^n T^* \mc{G} \otimes \CC$.
  All the results presented above then carry over to complex coefficients.
  In particular, the \v{C}ech-de Rham-Satake theorem then states that
  \begin{equation*}
    H_\text{dR}^\bullet(\mf{X};\CC) \cong H^\bullet(|\mf{X}|;\CC)
  \end{equation*}
  for a presentable orbifold $\mf{X}$.
  We can also define compactly supported de Rham cohomology, where the support of a differential form $\go \in \gO^n(\mc{G})$ is defined to be the closure of the set
  \begin{equation*}
    \gp\left(\left\{ x \in G_0 : \go(x) \neq 0\right\}\right) \subset |\mc{G}|
  \end{equation*}
  where $\map{\gp}{G_0}{|\mc{G}|}$ is the quotient map.
  The \v{C}ech-de Rham-Satake theorem holds for compactly supported cohomology (replacing simplicial cohomology with compactly supported simplicial cohomology).
  Consequently we have a covariant Mayer-Vietoris sequence: if $\mf{X}$ is a presentable orbifold with underlying space $X$,  and if $\{U,V\}$ is an open cover of $X$, then there is an exact sequence
  \begin{equation*}
    \xymatrix{
      \cdots \to H_\text{dR}^p(\mf{X}|_{U \cap V}) \to H_\text{dR}^p(\mf{X}|_U) \oplus H_\text{dR}^p(\mf{X}|_V) \to H_\text{dR}^p(\mf{X}) \to H_\text{dR}^{p+1}(\mf{X}|_{U \cap V}) \to \cdots
    }
  \end{equation*}
  given by the corresponding covariant Mayer-Vietoris exact sequence in compactly supported simplicial cohomology.
\end{rmk}

\section{Riemannian metrics}\label{metrics}

\begin{dfn}
  A \emph{Riemannian metric} on an \'etale groupoid $\mc{G}$ is a positive-definite section
  \begin{equation*}
    g \in \gG(\mc{G},\gS^2 T^* \mc{G})
  \end{equation*}
  of the $2^\text{nd}$ symmetric power of the cotangent bundle of $\mc{G}$.
  Equivalently, a Riemannian metric on $\mc{G}$ is given by a Riemannian metric on the smooth manifold $G_0$ which is $\mc{G}$-invariant in the sense that if $\map{\gg}{x}{y}$ is an arrow in $\mc{G}$, then
  \begin{equation*}
    g_y(\gg v,\gg w) = g_x(v,w) \qquad (v,w \in T_x G_0),
  \end{equation*}
  recalling that $TG_0$ is naturally a $\mc{G}$-space.
\end{dfn}

For two vector fields $\gx,\gn \in \gG(\mc{G},T\mc{G})$ on an \'etale groupoid equipped with a Riemannian metric $g$, we let $\langle \gx,\gn \rangle \in C^\infty(\mc{G})$ denote the smooth function
\begin{equation*}
  x \mapsto g_x(\gx_x,\gn_x).
\end{equation*}
As with differential forms, a Riemannian metric on an orbifold $\mf{X}$ is defined to be a Riemannian metric on a representing groupoid $\mc{G}$.

\begin{thm}\label{metricsexist}
  Let $\mf{X}$ be an orbifold with groupoid representation $|\mc{G}| \to X$.
  Then there exists a Riemannian metric on $\mf{X}$.
\end{thm}

\begin{proof}
  Using paracompactness of $X$ and Theorem \ref{localaction}, take a locally finite open cover $\{U_\ga\}_{\ga \in A}$ of $X$ such that for each $\ga$
  \begin{equation*}
    \mc{G}|_{\wtd{U}_\ga} \cong \mc{G}_{x(\ga)} \ltimes \wtd{U}_\ga
  \end{equation*}
  using the notation from the proof of Theorem \ref{po1exists}.
  On each $\wtd{U}_\ga$ there exists a $\mc{G}|_{\wtd{U}_\ga}$-invariant Riemannian metric $\wtd{g}_\ga$.
  To see this, for a fixed $\ga$ take any Riemannian metric $g_\ga$ on $\wtd{U}_\ga$ and define for $y \in \wtd{U}_\ga$ and $v,w \in T_y \wtd{U}_\ga$
  \begin{equation*}
    \wtd{g}_\ga(v,w) := \frac{1}{|\mc{G}_{x(\ga)}|} \sum_{\gg \in \mc{G}_{x(\ga)}} g_\ga(\gg v, \gg w).
  \end{equation*}
  Just as in Theorem \ref{po1exists}, this averaging procedure defines a $\mc{G}|_{\wtd{U}_\ga}$-invariant Riemannian metric on $\wtd{U}_\ga$.

  Now we can use a partition of unity to piece together these locally-defined metrics in order to obtain a globally-defined one.
  Let $\{\gr_\ga\}$ be a partition of unity subordinate to the cover $\{U_\ga\}$, and define
  \begin{equation*}
    g_x := \sum_{\ga} \gr_\ga(x) \wtd{g}_\ga.
  \end{equation*}
  Exactly as in the smooth manifold case (see \cite[Proposition III.1.4]{KN63}), this defines a Riemannian metric on the orbifold $\mf{X}$.
\end{proof}

\begin{rmk}\label{herm}
  By generalising the definition of a Riemannian metric, we can consider metrics on any vector bundle (real or complex\footnote{By a metric over a complex vector bundle we mean a Hermitian metric.}) over a groupoid $\mc{G}$.
  The proof of Theorem \ref{metricsexist} carries over to this more general context, telling us that there exist metrics on any vector bundle over an orbifold.
\end{rmk}

\section{Connections and curvature}

\begin{dfn}
  Let $E$ be a vector bundle over an \'etale groupoid $\mc{G}$.
  A \emph{connection} on $E$ is an $\RR$-linear\footnote{Recall our convention of assuming all vector bundles are real vector bundles unless otherwise stated.} map
  \begin{equation*}
    \map{\nabla}{\gG(\mc{G},E)}{\gO^1(\mc{G},E) := \gG(\mc{G},\wedge^1 T^* M \otimes E)}
  \end{equation*}
  such that the Leibniz rule
  \begin{equation}\label{leib}
    \nabla(\gf \gx) = d\gf \otimes \gx + \gf\nabla\gx
  \end{equation}
  holds for all functions $\gf \in C^\infty(\mc{G})$ and sections $\gx \in \gG(\mc{G},E)$.
  Given such a connection $\nabla$ and a vector field $\gh \in \gG(\mc{G},T\mc{G})$, we define
  \begin{equation*}
    \map{\nabla_\gh}{\gG(\mc{G},E)}{\gG(\mc{G},E)}
  \end{equation*}
  to be the operator $\gi(X) \nabla$, where $\map{\gi(\gh)}{\gO^1(\mc{G},E)}{\gG(\mc{G},E)}$ is contraction by $\gh$.\footnote{See \cite[Definitions 1.6 and 1.14]{BGV92}}
\end{dfn}

Connections on orbifold bundles are defined in the now-apparent way.

\begin{example}\textbf{(Levi-Civita connection.)}
  Let $\mf{X}$ be an orbifold with representing groupoid $\mc{G}$ and with a Riemannian metric $g$ defined on $T\mc{G}$.
  Then, in particular, $G_0$ is a Riemannian manifold, and so we have the \emph{Levi-Civita connection} $\nabla^\text{LC}$ defined on the tangent bundle $TG_0$.
  This is the unique connection on $TG_0$ such that
  \begin{align*}
    d\langle \gx,\gh \rangle &= \langle \nabla^\text{LC} \gx,\gh \rangle + \langle \gx, \nabla^\text{LC} \gh \rangle \quad \text{and} \\
    [\gx,\gh] &= \nabla_\gx^\text{LC} \gh - \nabla_\gh^\text{LC} \gx
  \end{align*}
  for all vector fields $\gx,\gh \in \gG(G_0,TG_0)$.
  Recall once more the action of an arrow $\map{\gg}{x}{y}$ in $\mc{G}$ as a local diffeomorphism from $x$ to $y$.
  By the definition of a Riemannian metric on an \'etale groupoid, $\gg$ in fact acts as a local isometry.
  Since the Levi-Civita connection is determined by local properties of the metric on $G_0$, it follows that $\nabla^\text{LC}$ is $\mc{G}$-invariant.
  Therefore $\nabla^\text{LC}$ defines a connection on $T\mc{G}$, and in particular on $T\mf{X}$.
\end{example}

Every vector bundle over a smooth manifold admits a connection; the proof (see \cite[Proposition 3.3.4]{lN07}) is a standard partition of unity argument.
The same is true for vector bundles over orbifolds.

\begin{prop}\label{connections}
  Let $\mf{X}$ be an orbifold with groupoid representation $|\mc{G}|\to X$, and suppose $E \to \mc{G}$ is a $\mc{G}$-vector bundle.
  Then there exists a connection on $E$.
\end{prop}

\begin{proof}
  As usual, take a locally finite open cover $\{U_\ga\}$ of $X$ with
  \begin{equation*}
    \mc{G}|_{\wtd{U}_\ga} \cong \mc{G}_{x(\ga)} \ltimes \wtd{U}_\ga,
  \end{equation*}
  again using notation from the proof of Theorem \ref{po1exists}, and take a partition of unity $\{\gr_\ga\}$ subordinate to this cover.
  Over each $\wtd{U}_\ga$, the vector $\mc{G}$-bundle $E$ takes the form of a $\mc{G}_{x(\ga)}$-equivariant vector bundle over $\wtd{U}_\ga$.
  Fixing $\ga$, let $\nabla_\ga$ be any connection on $E_\ga$, and define a new $\mc{G}_{x(\ga)}$-equivariant connection $\wtd{\nabla}_\ga$ on $E_\ga$ by
  \begin{equation*}
    \wtd{\nabla}_\ga(\gx) := \frac{1}{|\mc{G}_x(\ga)|} \sum_{\gg \in \mc{G}_{x(\ga)}} \nabla_\ga (\gg \gx) \quad (\gx \in \gG(\wtd{U}_\ga,E_\ga)),
  \end{equation*}
  where $\gg \in \mc{G}_{x(\ga)}$ acts on $\gG(\wtd{U}_\ga,E_\ga)$ by
  \begin{equation*}
    \gg\gx(x) = \gg(\gx(\gg^{-1} x)) \quad (x \in \wtd{U}_\ga).
  \end{equation*}
  It is straightforward to check that $\wtd{\nabla}_\ga$ defines a connection on the vector $\mc{G}|_{\wtd{U}_\ga}$-bundle $E_\ga$ (noting that the Leibniz rule \eqref{leib} only holds for $\mc{G}|_{\wtd{U}_\ga}$-equivariant sections).

  Now define a connection $\nabla$ on the $\mc{G}$-bundle $E$ by
  \begin{equation*}
    \nabla(\gx) := \sum_\ga \gr_\ga \wtd{\nabla}_\ga (\gx|_{\wtd{U}_\ga}) \quad (\gx \in \gG(\mc{G},E)).
  \end{equation*}
  As a locally finite sum of smooth sections, $\nabla(\gx)$ is a well-defined smooth section of $E$ over $G_0$, whose $\mc{G}$-invariance follows immediately.
  Since the $\wtd{\nabla}_\ga$ are connections in the local $\mc{G}|_{\wtd{U}_\ga}$-bundles $E_\ga$, it follows that $\nabla$ is a connection.
\end{proof}

Exactly as in the smooth manifold case, we get the following corollary.

\begin{cor}
  The set of connections on a vector bundle $E$ over an orbifold $\mf{X}$ is a non-empty affine space modelled on the space $\gO^1(\End(E))$ of $\End(E)$-valued $1$-forms.
\end{cor}

\begin{rmk}
  Our proof of Proposition \ref{connections} is essentially the same as that given in \cite[Lemma 4.3.2]{CR02}.
  However, the proof is much simpler when presented in the language of groupoids.
\end{rmk}

\begin{rmk}
  Although we did not restrict the definition of a connection to good vector bundles, the concept is not paticularly useful for bad vector bundles---as shown in Example \ref{badbundle}---without modification (such a modification is given in \cite{cS07}, with application to Chern-Weil theory).
\end{rmk}

\begin{rmk}
  One should ask whether or not connections exist on (vector or principal) bundles over non-orbifold \'etale groupoids.
  It turns out that not all principal bundles over such groupoids admit connections.
  This is explored in \cite[\textsection 3.3]{LTX07} (a specific counterexample is given in Example 3.12).
\end{rmk}

Let $E \to M$ be a vector bundle over a smooth bundle, and suppose $\nabla$ is a connection on $E$.
Recall that the \emph{curvature} of $\nabla$ is an $\End(E)$-valued $2$-form
\begin{equation*}
  F_\nabla \in \gO^2(\End(E))
\end{equation*}
defined by
\begin{equation}\label{curv}
  F_\nabla(\gx,\gh) := [\nabla_\gx,\nabla_\gh] - \nabla_{[\gx,\gh]}
\end{equation}
for vector fields $\gx,\gh \in \gG(M,TM)$ \cite[\textsection 1.1]{BGV92}.
By using \eqref{curv} we can define the curvature of a connection on a $\mc{G}$-bundle $E \to \mc{G}$ over an \'etale groupoid, provided that the Lie bracket of two $\mc{G}$-equivariant vector fields is itself $\mc{G}$-equivariant.
We shall show that this is the case.
Let $\gx,\gh \in \gG(G_0,TG_0)$ be two $\mc{G}$-equivariant vector fields, and write the Lie bracket as
\begin{equation*}
  [\gx,\gh] = \left. \frac{d}{dt} \right|_{t=0} \gf_t \gh,
\end{equation*}
where $\{\gf_t\}$ is the family of diffeomorphisms of $G_0$ given by the flow of $\gx$.
To show that $[\gx,\gh]$ is $\mc{G}$-equivariant, it thus suffices to show for any arrow $\map{\gg}{x}{y}$ in $\mc{G}$ that
\begin{equation}\label{floweq}
  \gf_t(y) = \wtd{\gg} \cdot \gf_t(x)
\end{equation}
for sufficiently small $t>0$, where $\wtd{\gg}$ denotes the action of $\gg$ as a local diffeomorphism from $x$ to $y$.
But since $\wtd{\gg}$ is a local diffeomorphism and $D\wtd{\gg}(\gx_x) = \gx_y$ (as $\gx$ is a $\mc{G}$-equivariant vector field), the equality \eqref{floweq} holds for small $t$.
Therefore equation \eqref{curv} can be used to define the curvature of a connection $\nabla$ on a vector $\mc{G}$-bundle $E$.

\section{Chern-Weil theory}\label{ccs}

Let $E \to M$ be a vector bundle of rank $r$ over the field $\FF$ (here $\FF$ is either $\RR$ or $\CC$) over a smooth manifold $M$, and suppose $\nabla$ is a connection on $E$.
The \emph{Chern-Weil construction} (see \cite[\textsection III.3]{rW08}) associates to each invariant $k$-linear form\footnote{A $k$-linear form over $M_{r \times r}(\FF)$ is a multilinear map
  \begin{equation*}
    \map{\gf}{\underbrace{M_{r \times r}(\FF) \times \cdots \times M_{r \times r}(\FF)}_\text{$k$ times}}{\FF};
  \end{equation*}
  such a form is \emph{invariant} if
  \begin{equation*}
    \gf(g A_1 g^{-1},\ldots,g A_k g^{-1}) = \gf(A_1,\ldots,A_k)
  \end{equation*}
  for all $g \in \GL(r,\FF)$ and $A_j \in M_{r \times r}(\FF)$ \cite[\textsection III.3]{rW08}.} $\gf \in I_k(M_{r \times r}(\FF))$ a differential $2k$-form $\gf(F_\nabla) \in \gO^{2k}(M;\FF)$ with coefficients in $\FF$, defined geometrically in terms of the curvature $F_\nabla$ of the connection $\nabla$.
It turns out that the $2k$-form $\gf(F_\nabla)$ is closed, and hence defines a de Rham cohomology class $[\gf(F_\nabla)] \in H_\text{dR}^{2k}(M;\FF)$; the miracle of Chern-Weil theory is that this class \emph{does not depend on the connection $\nabla$}.\footnote{The \emph{form} $\gf(F_\nabla)$ does depend on $\nabla$, but only up to an exact form.}
Given that there exists a connection on $E$, we can hence view this construction as a way to coax topologically meaningful data (i.e. cohomology classes) out of vector bundles by means of invariant $k$-linear forms.
Since giving an invariant $k$-linear form $\gf$ on $M_{r \times r}(\FF)$ is equivalent to giving an invariant map $\map{\gf}{M_{r \times r}(\FF)}{\FF}$ such that $\gf(A)$ is a degree $k$ homogeneous polynomial in the entries of $A$, we usually consider the Chern-Weil construction as sending these `invariant polynomials' to cohomology classes.

\begin{example}\label{chernclass}\textbf{(Chern classes.)} \cite[Example III.3.1, Definition 3.4]{rW08}
  Let $E \to M$ be a \emph{complex} vector bundle of rank $r$ over a smooth manifold.
  The determinant map
  \begin{equation*}
    \map{\det}{M_{r \times r}(\CC)}{\CC}
  \end{equation*}
  is multiplicative, so that $\det(gAg^{-1}) = \det(A)$ for all $g \in \GL(r,\CC)$ and $A \in M_{r \times r}(\CC)$.
  Consider the normalised characteristic polynomial
  \begin{equation*}
    \det\left(I-\frac{i}{2\gp}tA\right) = \sum_{k=0}^r \gF_k(A) t^k
  \end{equation*}
  of $A$.
  Each $\gF_k(A)$ is a homogeneous polynomial in the entries of $A$, so by the Chern-Weil construction, each $\gF_k$ determines a cohomology class
  \begin{equation*}
    c_k(E) \in H_\text{dR}^{2k}(M;\CC) 
  \end{equation*}
  called the \emph{$k^\text{th}$ Chern class} of the vector bundle $E$.
  The \emph{total Chern class} of $E$ is defined to be the sum
  \begin{equation*}
    c(E) := \sum_{k=1}^r c_k(E) \in H_\text{dR}^\text{ev}(M;\CC),
  \end{equation*}
  living in the even-dimensional de Rham cohomology of $M$ with complex coefficients.
\end{example}

Generalising the Chern-Weil construction to bundles over orbifolds is easy.
Let $E \to \mc{G}$ be a rank $r$ $\FF$-vector bundle over an orbifold groupoid $\mc{G}$ with a connection $\nabla$, and suppose $\gf \in I_k(M_{r \times r(\FF)})$.
Following the construction in the smooth manifold case \cite[\textsection III.3]{rW08}, working with the vector bundle $E \to G_0$, and assuming all objects are $\mc{G}$-equivariant, we see that the resulting $2k$-form $\gf(F_\nabla)$ is $\mc{G}$-equivariant and hence defines a closed differential form
\begin{equation*}
  \gf(F_\nabla) \in \gO^{2k}(\mc{G};\FF).
\end{equation*}
Thus the usual consequences of the Chern-Weil construction follow \emph{mutatis mutandis} in the case of orbifold vector bundles.
In particular, following Example \ref{chernclass} we can define the \emph{$k^\text{th}$ Chern class}
\begin{equation*}
  c_k(E) \in H_\text{dR}^{2k}(\mf{X};\CC)
\end{equation*}
of a rank $r$ complex vector bundle $E$ over an orbifold $\mf{X}$, as well as the \emph{total Chern class}
\begin{equation*}
  c(E) := \sum_{k=1}^r c_k(E) \in H_\text{dR}^\text{ev}(\mf{X};\CC).
\end{equation*}

\begin{example}\label{chernchar}\textbf{(Chern character.)}
  Consider the map $\map{\gf}{M_{r \times r}(\CC)}{\CC}$ given by
  \begin{equation*}
    A \mapsto \tr\left( \exp\left(\frac{A}{2\gp i}\right)\right).
  \end{equation*}
  As with the map defining the Chern class, this can be expanded as a formal power series
  \begin{equation*}
    A \mapsto \sum_{k=0}^\infty \gY_k(A)t^k,
  \end{equation*}
  where each $\gY_k$ is a degree $k$ homogeneous polynomial in the entries of $A$.
  The function $\gf$ is invariant; to see this, let $g$ be in $\GL(r,\CC)$ and write
  \begin{align*}
    \gf(gAg^{-1})
    &= \tr\left(\exp\left(g\frac{A}{2\gp i} g^{-1} \right) \right) \\
    &= \tr\left(g\exp\left(\frac{A}{2\gp i}\right)g^{-1}\right) \\
    &= \tr\left(\exp\left(\frac{A}{2\gp i}\right)g^{-1} g\right) \\
    &= \gf(A),
  \end{align*}
  first using that $\exp(gAg^{-1}) = g\exp(A)g^{-1}$ for invertible $g$, and then using that the trace of a product is invariant under cyclic permutations of the terms of the product.
  If $E \to \mf{X}$ is a complex orbifold vector bundle of rank $r$, the sum $\sum_k \gY_k$ hence determines a sum of cohomology classes.
  This sum is denoted $\Ch(E)$ and called the \emph{Chern character} of $E$.

  Suppose $F \to \mf{X}$ is a second complex orbifold vector bundle.
  As in the case of vector bundles over manifolds, it can be shown that $\Ch(E \oplus F) = \Ch(E) + \Ch(F)$ and $\Ch(E \otimes F) = \Ch(E)\Ch(F)$.
  Furthermore, if $\underline{\CC}$ denotes the trivial complex line bundle on $\mf{X}$, then $\Ch(\underline{\CC}) = 1 \in H_\text{dR}^0(\mf{X};\CC)$.
  Hence the Chern character defines a semiring homomorphism $\Vect_\CC(\mf{X}) \to H_\text{dR}^\text{ev}(\mf{X};\CC)$.
  By the universal property of the Grothendieck group, this induces a ring homomorphism
  \begin{equation*}
    \map{\Ch}{K^0(\mf{X})}{H_\text{dR}^\text{ev}(\mf{X};\CC)}
  \end{equation*}
  when $\mf{X}$ is compact.
\end{example}

\section{$K$-theory of presentable orbifolds}\label{presentablek}

\begin{rmk}
  From now on we will only consider complex vector bundles, so the term `vector bundle' will now exclusively mean `complex vector bundle'.
\end{rmk}

Let $\mf{X}$ be a presentable orbifold, i.e. an orbifold with representing groupoid Morita equivalent to an action groupoid $G \ltimes M$, where $G$ is a compact Lie group acting smoothly and almost freely on a manifold $M$.
The orbifold $\mf{X}$ is compact if and only if the orbit space $G \sm M$ is compact; by \cite[Theorem I.3.1]{gB72}, this holds if and only if $M$ is compact.
As noted in Remark \ref{actionbundles}, vector bundles over $G \ltimes M$ can be identified with $G$-equivariant vector bundles over $M$.
Furthermore, complexes of vector bundles over $G \ltimes M$ can be identified with complexes of $G$-equivariant vector bundles over the $G$-space $M$ (using the terminology of \cite[\textsection 3]{gS68}).
Both of these identifications respect the operations $\oplus$ and $\otimes$, so therefore the $K$-theory $K^0(G \ltimes M)$ of the groupoid $G \ltimes M$ can be identified with the $G$-equivariant $K$-theory $K_G^0(M)$ of the $G$-space $M$ (taking compactly supported equivariant $K$-theory when $M$ is noncompact).\footnote{For the details of equivariant $K$-theory, see \cite{gS68,PR84,mA67}.}

\begin{rmk}
  We have implicitly made use of the fact that $K_G^0(M)$ can be equivalently defined in terms of \emph{smooth} vector bundles rather than continuous vector bundles. Vector bundles over Lie groupoids are, by definition, smooth vector bundles.
\end{rmk}

Since $G \ltimes M$ is Morita equivalent to a representing groupoid for the orbifold $\mf{X}$, we can identify $K^0(\mf{X})$ with $K_G^0(M)$.
In fact, for each $n \in \NN$ we can identify $K^{-n}(\mf{X})$ with $K^{-n}_G(M)$, and so the $K$-theory of presentable orbifolds reduces to the equivariant $K$-theory of manifolds.
This is extremely useful, as all the properties of equivariant $K$-theory (which are well-known) carry over to the $K$-theory of presentable orbifolds.
Compare this to the case of general orbifolds, in which we could not say much at all about $K$-theory (keeping Conjecture \ref{conjecture} in mind).

In particular, since Bott periodicity holds for equivariant $K$-theory, it also holds for presentable orbifold $K$-theory, so that
\begin{equation*}
  K^{-n-2}(\mf{X}) \cong K^{-n}(\mf{X})
\end{equation*}
for all $n \in \NN$.
We hence define $K^n(\mf{X}) := K^{-n}(\mf{X})$ for all $n \in \NN$, and write
\begin{equation*}
  K^\bullet(\mf{X}) := K^0(\mf{X}) \oplus K^1(\mf{X}),
\end{equation*}
so that $K^\bullet$ is a $\ZZ_2$-graded ring.
Consequently, $K$-theory is a $\ZZ_2$-graded cohomology theory when restricted to presentable orbifolds.\footnote{Of course, we could have just said `since equivariant $K$-theory is a $\ZZ_2$-graded cohomology theory, so is presentable orbifold $K$-theory'.}
The covariant Mayer-Vietoris sequence also carries over from compactly supported equivariant $K$-theory: supposing $\mf{X}$ is a presentable orbifold with $|\mf{X}|= X$, and taking an open cover $\{U,V\}$ of $X$, we have an exact sequence

\begin{equation*}
  \xymatrix{
    K^0(\mf{X}|_{U \cap V}) \ar[r] & K^0(\mf{X}|_U) \oplus K^0(\mf{X}|_V) \ar[r] &  K^0(\mf{X}) \ar[d] \\
    K^1(\mf{X}) \ar[u] & K^1(\mf{X}|_U) \oplus K^1(\mf{X}|_V)  \ar[l] & K^1(\mf{X}|_{U \cap V}). \ar[l].
  }
\end{equation*}

\section{The inertia orbifold}\label{inof}

In Section \ref{ineffective} we gave a vague description of the `inertia orbifold' $\mf{X}$ associated to an orbifold $\mf{X}$.
Now we shall give a rigorous definition of $\wtd{\mf{X}}$.
In the next section, we will see the use of this concept in studying the $K$-theory of presentable orbifolds.

Let $\mc{G}$ be a Lie groupoid, and define
\begin{equation*}
  S_\mc{G} := \fp{G_0}{\gD}{(s,t)}{G_1} = \{\gg \in G_1 \mid s(\gg) = t(\gg)\}
\end{equation*}
where $\map{\gD}{G_0}{G_0 \times G_0}$ is the diagonal map.
Since $s$ and $t$ are submersions, $(s,t)$ is also a submersion, and thus $S_\mc{G}$ is a smooth manifold.
Elements of $S_\mc{G}$ are `loops' in $\mc{G}$, i.e. arrows which start and end at the same point.
The manifold $S_\mc{G}$ can be made into a $\mc{G}$-space in the following way.
The anchor $S_\mc{G} \to G_0$ is the map $\map{\pr_1}{\fp{G_0}{\gD}{(s,t)}{G_1}}{G_0}$ which sends a loop to its `basepoint'.
If $\map{\gg}{x}{y}$ is an arrow in $\mc{G}$ and $\gd$ is a loop at $x$, then $\gg\gd\gg^{-1}$ is a loop at $y$.
Defining
\begin{equation*}
  \gQ(\gd,\gg) := \gg\gd\gg^{-1}
\end{equation*}
makes $S_\mc{G}$ into a $\mc{G}$-space.
The resulting action groupoid $\mc{G} \ltimes S_\mc{G}$ is called the \emph{inertia groupoid} of $\mc{G}$ and denoted $\wtd{\mc{G}}$.

Consider the `basepoint' map
\begin{equation*}
  \map{\gb := \pr_1}{\fp{G_0}{\gD}{(s,t)}{G_1}}{G_0}.
\end{equation*}
from before.
This map induces a groupoid homomorphism $\map{\gb}{\wtd{\mc{G}}}{\mc{G}}$ with
\begin{equation*}
  \map{\gb_1 = \pr_1}{\fp{G_1}{s}{\gb}{S_\mc{G}}}{G_1}.
\end{equation*}
If the groupoid $\mc{G}$ is proper, then the map $\map{(s,t)}{G_1}{G_0}$ is proper, and hence so is $\gb$.

Suppose $\ell \in S_\mc{G}$ is a loop in $\mc{G}$ (i.e. an arrow $\ell \in G_1$ with $s(\ell) = t(\ell)$).
Then the isotropy group of $\ell$ in the inertia groupoid $\wtd{\mc{G}}$ is given by the group of arrows $\gd \in \mc{G}$ with $\gd\cdot \ell = \ell$.
By the definition of the $\mc{G}$-action on $S_\mc{G}$, this is the subgroup of arrows $\gd \in \mc{G}_{\gb(\ell)}$ with
\begin{equation*}
  \gd\ell \gd^{-1} = \ell,
\end{equation*}
or equivalently $\gd\ell = \ell \gd$.
Therefore $\wtd{\mc{G}}_\ell$ is precisely the centraliser subgroup $Z(\ell)$ of $\ell$ in the isotropy group $\mc{G}_{\gb(\ell)}$.

\begin{prop}\label{inert}
  Let $G$ be a compact Lie group acting smoothly and almost freely on a manifold $M$.
  Then the inertia groupoid $\wtd{G \ltimes M}$ is weakly equivalent to the disjoint union of action groupoids
  \begin{equation*}
    \bigsqcup_{g \in R} Z(g) \ltimes M^g,
  \end{equation*}
  where $R \subset G$ is any set of representatives for the conjugacy classes of $G$ and $Z(g)$ is the centraliser subgroup of $g$ in $G$.
\end{prop}

\begin{rmk}
  We haven't explicitly defined disjoint unions of Lie groupoids, but it is easy to see how they should be defined.
\end{rmk}

\begin{proof}
  For each $g \in R$, define a Lie groupoid homomorphism
  \begin{equation*}
    \map{\gf_g}{Z(g) \ltimes M^g}{\wtd{G \ltimes M}}
  \end{equation*}
  which sends an arrow $\map{(\gg,x)}{x}{y}$ in $Z(g) \ltimes M^g$ to the arrow
  \begin{equation*}
    \map{(\gg,(x,g))}{(x,g)}{(y,g)}
  \end{equation*}
  in $\wtd{G \ltimes M}$.
  This is well-defined since $\gg$ is in the centraliser $Z(g)$, and is clearly a Lie groupoid homomorphism.
  Define
  \begin{equation*}
    \map{\gf}{\bigsqcup_{g \in R} Z(g) \ltimes M^g}{\wtd{G \ltimes M}}
  \end{equation*}
  to be the disjoint union of these homomorphisms.

  Now we must show that $\gf$ is a weak equivalence.
  To see that the composition
  \begin{equation*}
    \xymatrix{
      \fp{\left(\bigsqcup_{g \in R} M^g\right)}{\gf_0}{s}{\left(\wtd{G \ltimes M}\right)_1} \ar[r]^(.685){\pr_2} & \left(\wtd{G \ltimes M}\right)_1 \ar[r]^(.57){t} & S_{G \ltimes M}
    }
  \end{equation*}
  is surjective, suppose $(k,y)$ is a loop in $S_{G \ltimes M}$ (i.e. $y \in M$ and $k \in G$, with $ky = y$).
  Then $k = \gg g \gg^{-1}$ for some $g \in R$ and $\gg \in G$, and hence the element
  \begin{equation*}
    ((x,M^g),(\gg,(g,x))) \in \fp{\left(\bigsqcup_{g \in R} M^g\right)}{\gf_0}{s}{\left(\wtd{G \ltimes M}\right)_1}
  \end{equation*}
  is mapped by $t \circ \pr_2$ onto $(k,y)$.
  This proves surjectivity.
  To see that the map is a submersion, suppose $g \in R$ and observe that a sufficiently small neighbourhood of $(g,y) \in S_{G \ltimes M}$ is of the form
  \begin{equation*}
    \{(g,z) : z \in N_y\}
  \end{equation*}
  where $N_y$ is a small neighbourhood of $y$ in $M$.
  This holds because the action of $G$ on $M$ is almost free, so using the slice theorem we can locally interpret the action as an action of a finite (hence discrete) group.
  The correspondence
  \begin{equation*}
    (g,z) \mapsto ((z,M^g),(1,z))
  \end{equation*}
  then defines a local section of $t \circ \pr_2$, where $1 \in G$ is the identity element.
  Therefore $t \circ \pr_2$ is a submersion.

  It remains to show the second condition for $t \circ \pr_2$ to be a weak equivalence; this amounts to saying that the set of arrows between two points
  \begin{equation*}
    (x,M^h),(y,M^k) \in \bigsqcup_{g \in R} Z(g) \ltimes M^g
  \end{equation*}
  can be smoothly identified with the set of arrows between the points
  \begin{equation*}
    (h,x), (k,y) \in S_{G \ltimes M}.
  \end{equation*}
  Since $\sqcup_{g \in R} Z(g) \ltimes M^g$ is a disjoint union, there is only an arrow between $(x,M^h)$ and $(y,M^k)$ if $h=k$, in which case the arrows are of the form
  \begin{equation*}
    \map{\gg}{(x,M^h)}{(\gg x,M^h)}
  \end{equation*}
  for $\gg \in Z(h)$.
  Now suppose $\map{\gs}{(h,x)}{(h,y)}$ is an arrow in $\wtd{G \ltimes X}$.
  Then $y = \gs x$ and $h = \gs h \gs^{-1}$, so $\gs$ is in $Z(h)$ and the required identification is apparent.
\end{proof}

Consider a presentable orbifold $\mf{X}$ defined by the action groupoid $G \ltimes X$.
Since $\wtd{G \ltimes X}$ is weakly equivalent to a disjoint union of orbifold groupoids (which must then be an orbifold groupoid), $\wtd{G \ltimes X}$ defines an orbifold.
We denote this orbifold by $\wtd{\mf{X}}$ and call it the \emph{inertia orbifold} associated to $\mf{X}$.
More generally, if $\mc{G}$ is an orbifold groupoid, then $\wtd{\mc{G}}$ is weakly equivalent to an orbifold groupoid; this can be proven similarly to Proposition \ref{inert} by way of Theorem \ref{localaction}.
We will restrict our attention to presentable orbifolds henceforth.

Since weak equivalences induce homeomorphisms of orbit spaces, the following result is an immediate corollary of Proposition \ref{inert}.

\begin{cor}
  Let $\mf{X} = G \sm M$ be a presentable orbifold.
  Then the orbit space $|\wtd{\mf{X}}|$ is homeomorphic to the disjoint union
  \begin{equation}\label{in}
    \bigsqcup_{g \in R} |Z(g) \ltimes M^g|
  \end{equation}
  with notation as in Proposition \ref{inert}.
\end{cor}

The homeomorphism in this corollary yields a nice decomposition of the de Rham cohomology of a presentable orbifold $\wtd{\mf{X}}$.
Using the \v{C}ech-de Rham-Satake Theorem \ref{cdrs}, we can write
\begin{align*}
  H_\text{dR}^\bullet(\wtd{\mf{X}};\CC)
  &\cong H^\bullet(|\wtd{\mf{X}}|;\CC) \\
  &\cong \bigoplus_{(g) \subset G} H^\bullet(|Z(g) \ltimes M^g|;\CC) \\
  &\cong \bigoplus_{(g) \subset G} H_\text{dR}^\bullet(Z(g) \sm M^g;\CC).
\end{align*}
Of course, we obtain the same isomorphisms whether we use real or complex coefficients.
The de Rham cohomology of the inertia orbifold $\wtd{\mf{X}}$ will play an important role in the $K$-theory of the orbifold $\mf{X}$, so this expression will come in handy.

\section{The delocalised Chern character}\label{dcci}

In Section \ref{ccs}, we defined the Chern character homomorphism
\begin{equation*}
  \map{\Ch}{K^0(\mf{X})}{\HdRev(\mf{X};\CC)}
\end{equation*}
when $\mf{X}$ is a compact orbifold.
This can be extended to $\ZZ_2$-graded ring homomorphism
\begin{equation}\label{z2ch}
  \map{\Ch}{K^\bullet(\mf{X})}{\HdR^\bullet(\mf{X};\CC)},
\end{equation}
for all presentable orbifolds $\mf{X}$, where $\HdR^\bullet = \HdRev \oplus \HdRodd$ is viewed as a $\ZZ_2$-graded ring and where we use compactly supported de Rham cohomology.
First we extend the definition to non-compact orbifolds: if $\map{\gs}{E_0}{E_1}$ is a a compactly supported vector bundle over an orbifold $\mf{X}$, then we define
\begin{equation*}
  \Ch\left(E_0 \stackrel{\gs}{\to} E_1\right) := \Ch(E_0) - \Ch(E_1).
\end{equation*}
This defines a ring homomorphism from $K^0(\mf{X}) := K_c^0(\mf{X})$ into the compactly supported de Rham cohomology $\HdRev(\mf{X};\CC)$.
When $\mf{X}$ is compact, then this definition coincides with the usual definition of the Chern character.\footnote{This is because the isomorphism $K_c^0(\mf{X}) \to K^0(\mf{X})$ is defined by sending a compactly supported vector bundle $E_0 \to E_1$ to the $K$-theory class $[E_0] - [E_1]$.}
Now we extend the definition to $K^1$.
Recall that $K^1(\mf{X})$ is defined to be $K^0(\mf{X} \times \RR)$ (whether or not $\mf{X}$ is compact).
The Chern character maps $K^0(\mf{X} \times \RR)$ into $\HdRev(\mf{X} \times \RR; \CC)$.
When $\mf{X} = M$ is a manifold, we have an isomorphism of compactly supported cohomology rings
\begin{equation*}
  \HdR^\bullet(M \times \RR;\CC) \cong \HdR^{\bullet - 1}(M;\CC)
\end{equation*}
given by `integration along the fibre'; a proof of this isomorphism is given in \cite[\textsection I.4]{rBlT82}.
This proof still works when we replace the manifold $M \times \RR$ with the orbifold $\mf{X} \times \RR$, yielding an isomorphism
\begin{equation*}
  \HdR^\bullet(\mf{X} \times \RR;\CC) \cong \HdR^{\bullet - 1}(\mf{X};\CC).
\end{equation*}
Composing the Chern character $K^0(\mf{X} \times \RR) \to \HdRev(\mf{X} \times \RR;\CC)$ with this isomorphism gives the definition of the `odd' Chern character
\begin{equation*}
  \map{\Ch}{K^1(\mf{X})}{\HdRodd(\mf{X};\CC)}.
\end{equation*}
Putting the `even' and `odd' Chern characters together gives the $\ZZ_2$-graded ring homomorphism in \eqref{z2ch}.

For a compact manifold $M$, the Chern character gives rise to an isomorphism
\begin{equation}\label{ci}
  \Ch \colon K^\bullet(M) \otimes_\ZZ \CC \stackrel{\sim}{\to} H_\text{dR}^\bullet(M;\CC)
\end{equation}
of $\ZZ_2$-graded rings (see \cite[Theorem V.3.25]{mK78}.
Consequently the torsion-free behaviour of $K^\bullet(M)$ is completely given by the de Rham cohomology of $M$.
For a compact presentable orbifold $\mf{X}$, the Chern character still gives a map from $K^\bullet(M) \otimes_\ZZ \CC$ to $H_\text{dR}^\bullet(\mf{X};\CC)$, but this map is no longer an isomorphism.
In its place, we have a `delocalised' Chern character $\Ch_\text{deloc}$ which induces an isomorphism
\begin{equation}\label{dci}
  \Ch_\text{deloc} \colon K^\bullet(\mf{X}) \otimes_\ZZ \CC \stackrel{\sim}{\to} H_\text{dR}^\bullet\left(\wtd{\mf{X}};\CC\right);
\end{equation}
note that the de Rham cohomology of $\mf{X}$ has been replaced with the de Rham cohomology of the associated inertia orbifold.
When $\mf{X}$ is a manifold (i.e. when $\mf{X}$ has no singular points), $\Ch_\text{deloc}$ and $\wtd{\mf{X}}$ will coincide with $\Ch$ and $\mf{X}$ respectively, and so \eqref{dci} will reduce to the isomorphism \eqref{ci}.

The remainder of this section will be devoted to constructing $\Ch_\text{deloc}$ and proving the isomorphism \eqref{dci}.
Before doing this we shall take a moment to reflect on the consequences of this isomorphism.
For a presentable orbifold $\mf{X} = G \sm M$, we find by the results of the previous section that
\begin{equation*}
  K^\bullet(\mf{X}) \otimes_\ZZ \CC \cong \bigoplus_{(g) \subset R} H^\bullet(|Z(g) \ltimes M^g|;\CC)
\end{equation*}
where $R$ is a set of representatives for the conjugacy classes in $G$.
Thus the torsion-free behaviour of $K^\bullet(\mf{X})$ can be computed topologically.
However, in contrast with the de Rham cohomology of $\mf{X}$, this behaviour is affected by the \emph{orbifold} structure of $\mf{X}$, not just the underlying topological space $|\mf{X}|$.
Therefore $K$-theory is a truly interesting orbifold invariant, even at the torsion-free level.

Now we shall define the delocalised Chern character.
Let $\mf{X}$ be an orbifold with representing groupoid $\mc{G}$, and suppose $E \to \mf{X}$ is a complex vector bundle represented by a vector bundle $E \to \mc{G}$.
Using the basepoint map $\map{\gb}{\wtd{\mc{G}}}{\mc{G}}$, pull $E$ back to a bundle $\gb^* E$ over $\wtd{\mc{G}}$.
Via the $\mc{G}$-action on $E$, each loop $\ell \in \mc{G}$ induces an isomorphism on the fibre $E_{\gb(\ell)}$.
This fibre is identified with $(\gb^* E)_\ell$.
Hence we can define a canonical automorphism $\gF_E$ of the vector bundle $\gb^* E$ which acts on each fibre $(\gb^* E)_\ell$ as the automorphism on $E_{\gb(\ell)}$ induced by the loop $\ell$.
If we equip the bundle $E$ with a Hermitian metric (by Remark \ref{herm} this can always be done), then each loop $\ell$ acts unitarily on the fibre $(\gb^* E)_\ell$, and so we have an eigenbundle decomposition
\begin{equation}\label{decom}
  \gb^* E = \bigoplus_{\gq \in S^1} E_\gq,
\end{equation}
where each $E_\gq$ is a complex vector bundle over $\wtd{\mc{G}}$ upon which $\gF_E$ acts as multiplication by $\gq \in S^1 \subset \CC$.
Note that this decomposition does not depend on the Hermitian metric on $E$.

\begin{dfn}
  With respect to the eigenbundle decomposition \eqref{decom}, we define the delocalised Chern character of the complex vector bundle $E \to \mf{X}$ by
  \begin{equation*}
    \Chdeloc(E) := \sum_{\gq \in S^1} \gq \Ch(E_\gq) \in H_\text{dR}^\text{ev}\left(\wtd{\mf{X}};\CC\right).
  \end{equation*}
\end{dfn}

When $\mf{X}$ is compact, this defines a map from $K^0(\mf{X})$ to $\HdRev\left(\wtd{X};\CC\right)$.
For a compactly supported vector bundle $\map{\gs}{E_0}{E_1}$ over a noncompact orbifold $\mf{X}$, we define
\begin{equation*}
  \Chdeloc\left(E_0 \stackrel{\gs}{\to} E_1\right) := \Chdeloc(E_0) - \Chdeloc(E_1).
\end{equation*}
This defines a map from $K_c^0(\mf{X})$ into the compactly supported de Rham cohomology $\HdRev\left(\wtd{X};\CC\right)$, and when $\mf{X}$ is compact this coincides with the pre-existing map from $K^0(\mf{X})$.

\begin{prop}
  The delocalised Chern character induces a $\ZZ_2$-graded ring homomorphism $K^\bullet(\mf{X}) \to \HdR^\bullet\left(\wtd{\mf{X}};\CC\right)$.
\end{prop}

\begin{proof}
  We will prove this for compact $\mf{X}$, as the general case is an immediate conseqence of the proof.
  
  First we must show that $\Ch_\text{deloc}$ is a semiring homomorphism from $\Vect_\CC(\mf{X})$.
  Let $\mc{G}$ be a representing groupoid for $\mf{X}$ and suppose $E,F \to \mc{G}$ are two complex vector bundles.
  Since $\mc{G}$ acts on $E \oplus F$ componentwise we can write
  \begin{equation*}
    \gF_{E \oplus F} = \gF_E \oplus \gF_F,
  \end{equation*}
  identifying $\gb^*(E \oplus F)$ with $\gb^* E \oplus \gb^*F$.
  Thus the eigenbundle decomposition of $\gb^*(E \oplus F)$ takes the form
  \begin{equation*}
    \gb^*(E \oplus F) = \bigoplus_{\gq \in S^1} E_\gq \oplus F_\gq,
  \end{equation*}
  and so we have
  \begin{align*}
    \Ch_\text{deloc}(E \oplus F) &= \sum_{\gq \in S^1} \gq \Ch(E_\gq \oplus F_\gq) \\
    &= \sum_{\gq \in S^1} \gq(\Ch(E_\gq) + \Ch(F_\gq)) \\
    &= \Ch_\text{deloc}(E) + \Ch_\text{deloc}(F).
  \end{align*}
  Similarly, the eigenbundle decomposition of $\gb^*(E \otimes F)$ takes the form
  \begin{equation*}
    \gb^*(E \otimes F) = \bigoplus_{\gq \in S^1} \bigoplus_{\substack{\gr,\gs \in S^1\\ \gr\gs = \gq}} E_\gr \otimes F_\gs = \bigoplus_{\gr,\gs \in S^1} E_\gr \otimes F_\gs,
  \end{equation*}
  and so
  \begin{align*}
    \Ch_\text{deloc}(E \otimes F) &= \sum_{\gr,\gs \in S^1} \gr\gs \Ch(E_\gr) \Ch(F_\gs) \\
    &= \Ch_\text{deloc}(E) \Ch_\text{deloc}(F).
  \end{align*}
  Finally, $\mc{G}$ acts trivially on the trivial complex line bundle $\underline{\CC} \to \mc{G}$, and so $\gF_{\underline{\CC}}$ is the identity on $\underline{\CC}$.
  Hence
  \begin{align*}
    \Ch_\text{deloc}(\underline{\CC}) = \Ch(\underline{\CC}) = 1.
  \end{align*}
  Therefore $\Ch_\text{deloc}$ is a semiring homomorphism $\Vect_\CC(\mf{X}) \to H_\text{dR}^\text{ev}\left(\wtd{\mf{X}};\CC\right)$, and thus induces a ring homomorphism
  \begin{equation*}
    \map{\Ch_\text{deloc}}{K^0(\mf{X})}{H_\text{dR}^\text{ev}\left(\wtd{\mf{X}};\CC\right)}.
  \end{equation*}

  We define $\Chdeloc$ on $K^1(\mf{X})$ similarly to the way we defined the Chern character: that is, by using the isomorphism
  \begin{equation*}
    \HdR^\bullet\left(\wtd{\mf{X}} \times \RR;\CC\right) \cong \HdR^{\bullet-1}\left(\wtd{\mf{X}};\CC\right)
  \end{equation*}
  given by integration along the fibre.
  Note that $\Chdeloc$ maps $K^0(\mf{X} \times \RR)$ into $\HdR^\bullet\left(\wtd{\mf{X} \times \RR};\CC\right)$; since the orbifolds $\wtd{\mf{X} \times \RR}$ and $\wtd{\mf{X}} \times \RR$ are equivalent, we have an isomorphism
  \begin{equation*}
    \HdR^\bullet\left(\wtd{\mf{X} \times \RR};\CC\right) \cong \HdR^\bullet\left(\wtd{\mf{X}} \times \RR;\CC\right),
  \end{equation*}
  completing the construction.
\end{proof}

Now we can get to work on proving isomorphism \eqref{dci}.
Our strategy is to prove the isomorphism for possibly non-compact quotients $G \sm M$ where $G$ is finite.
Every orbifold is locally of this form.
Using the Mayer-Vietoris sequence, we can extend this result to compact presentable orbifolds.

\begin{lem}
  Suppose $G$ is a finite group acting smoothly on a possibly non-compact smooth manifold $M$, and let $\mf{X} = G \sm M$ be the corresponding orbifold.
  Then the delocalised Chern character induces an isomorphism
  \begin{equation*}
    K^\bullet(\mf{X}) \otimes_\ZZ \CC \cong \HdR^\bullet\left(\wtd{\mf{X}};\CC\right).
  \end{equation*}
\end{lem}

\begin{proofsketch}
  First suppose $M$ is compact.
  The first part of the proof given by the proof of \cite[Theorem 2]{AS89}, which we briefly recall.
  Suppose $E \to M$ is a $G$-vector bundle.
  Each element $g \in G$ acts fibrewise upon the restriction of $E$ to $M^g$, and so $E|_{M^g}$ has an eigenbundle decomposition
  \begin{equation}\label{atdc}
    E|_{M^g} = \bigoplus_{\gq \in S^1} E_\gq
  \end{equation}
  where $g$ acts on $E_\gq$ as multiplication by $\gq$.
  The correspondence $E \mapsto \sum_\gq \gq E_\gq$ defines a map $\gf_g$ from $\Vect_\CC^G(M)$ into $K_{Z(g)}^0(M^g)$.
  Extending the maps $\gf_g$ to $K^1$, taking the direct sum, and tensoring with $\CC$ gives a homomorphism
  \begin{equation}\label{athom}
    K_G^\bullet(M) \otimes_\ZZ \CC \to \bigoplus_{g \in R} K_{Z(g)}^\bullet(M^g) \otimes_\ZZ \CC.
  \end{equation}
  Where, as usual, $R$ is a set of representatives for the conjugacy classes in $G$.
  The content of \cite[Theorem 2]{AS89} is that this map is an isomorphism.
  Identifying each $K_{Z(g)}^\bullet(M^g)$ with $K^\bullet(Z(g) \ltimes M^g)$ and composing the above map with the sum of the Chern character homomorphisms
  \begin{equation}\label{sumhoms}
    K^\bullet(Z(g) \ltimes M^g) \otimes_\ZZ \CC \to \HdR^\bullet(Z(g) \ltimes M^g;\CC)
  \end{equation}
  yields a homomorphism
  \begin{equation*}
    K^\bullet(\mf{X}) \otimes_\ZZ \CC \cong K_G^\bullet(M) \otimes_\ZZ \CC \to \bigoplus_{g \in R} \HdR^\bullet(Z(g) \ltimes M^g;\CC) \cong \HdR^\bullet\left(\wtd{\mf{X}};\CC\right).
  \end{equation*}
  Upon inspecting the decomposition \ref{atdc} and the map \ref{athom} we see that this homomorphism is precisely that given by $\Chdeloc$.
  By a result of Baum and Connes \cite[Theorem 1.19]{BC88}, the maps in \ref{sumhoms} are isomorphisms, completing the proof when $\mf{X}$ is compact.
  When $\mf{X}$ is not compact, the same argument works (replacing vector bundles with compactly supported vector bundles) since the result in \cite{BC88} does not require compactness of $M^g$.
\end{proofsketch}

\begin{thm}\cite[Proposition 2.7]{HW12}
  Let $\mf{X}$ be a compact presentable orbifold.
  Then the map
  \begin{equation*}
    \map{\Ch_\text{deloc}}{K^\bullet(\mf{X}) \otimes_\ZZ \CC}{H^\bullet_\text{dR}\left(\wtd{\mf{X}};\CC\right)}
  \end{equation*}
  is an isomorphism.
\end{thm}

\begin{proof}
  Since $\mf{X}$ is compact, we can cover the underlying space $X$ of $\mf{X}$ with finitely many open sets $U_i$ such that the restrictions $\mf{X}|_{U_i}$ are orbifolds given by the quotient of a finite group, as in the previous lemma.
  Proceeding by induction (noting that the lemma constitutes the base case), assume that $\{U_1,\ldots,U_{n+1}\}$ is an open cover of $X$ as above, and write
  \begin{equation*}
    A := \bigcup_{i=1}^n U_i \quad \text{and} \quad B := U_{n+1}.
  \end{equation*}
  Then using the covariant Mayer-Vietoris sequences for $K$-theory and compactly supported de Rham cohomology, we have a diagram
  \begin{equation*}
    \scalebox{0.8}{
      \xymatrix{
        K^0(A \cap B) \ar[d] \ar[r] & K^0(A) \oplus K^0(B) \ar[d] \ar[r] & K^0(\mf{X}) \ar[d] \ar[r] & K^1(A \cap B) \ar[d] \ar[r] & K^1(A) \oplus K^1(B) \ar[d] \\
        \HdRev(A \cap B) \ar[r] & \HdRev(A) \oplus \HdRev(B) \ar[r] & \HdRev(\mf{X}) \ar[r] & \HdRodd(A \cap B) \ar[r] & \HdRodd(A) \oplus \HdRodd(B)
      }
    }
  \end{equation*}
  in which the rows are exact and the vertical maps are given by $\Chdeloc$.\footnote{We have abbreviated $K^0(\mf{X}|_{A \cap B})$ by $K^0(A \cap B)$ (and so on) and ommitted references to complex coefficients in our de Rham cohomology.}
  Commutativity of this diagram follows from naturality of the Chern character.
  Upon tensoring the top row with $\CC$ (note that this preserves exactness since $\CC$ is free), all vertical maps except the central map become isomorphisms by the inductive hypothesis.
  This proves the required isomorphism for even $K$-theory.
  The isomorphism for odd $K$-theory is proven in the same way by inspecting the remaining portion of the Mayer-Vietoris sequence.
\end{proof}


\appendix
\chapter{The slice theorem}\label{gpactions}

The slice theorem (first proven in the smooth case in \cite{jK53}) is of fundamental importance in the study of smooth actions of compact Lie groups on manifolds.
We use it (along with some related results) in Sections \ref{effectives} and \ref{examples}, so we have colected these results in one place for easy reference.

To properly state the slice theorem, we need to make a few definitions.
These are adapted from \cite[\textsection 1]{mD11}.

\begin{dfn}
  Let $G$ be a group acting on a topological space $X$ (on the left), and let $H$ be a subgroup of $G$.
  Then $H$ acts on the product $G \times X$ by
  \begin{equation*}
    h(g,x) = (gh^{-1},hx) \qquad (h \in H).
  \end{equation*}
  The orbit space of this action is denoted $G \times_H X$ and called the \emph{twisted product} of $G$ and $X$.
  The orbit of a point $(g,x) \in G \times X$ is denoted $[g,x]$.
  The twisted product is a $G$-space with the action
  \begin{equation*}
    g^\prime[g,x] = [g^\prime g,x] \qquad (g^\prime \in G).
  \end{equation*}
\end{dfn}

\begin{dfn}
  Let $G$ be a Lie group acting smoothly on a manifold $M$, and let $x$ be a point of $M$.
  A \emph{smooth slice} at $x$ is a $G_x$-stable superset $U$ of $x$ in $M$ such that the $G$-map
  \begin{equation*}
    G \times_{G_x} U \to M, \quad [g,u] \mapsto gu
  \end{equation*}
  is a diffeomorphism onto a neighbourhood of the orbit $G(x)$.
  If $U$ is homeomorphic to an open ball in $\RR^n$ for some $n$, then $G \times_{G_x} U$ is called an \emph{equivariant tubular neighbourhood} of $G(x)$.
\end{dfn}

\begin{thm}\label{slice}
  Let $G$ be a compact Lie group acting smoothly on a manifold $M$.
  Then for each point $x \in M$ there exists an equivariant tubular neighbourhood of $G(x)$.
\end{thm}

In fact, we can find an equivariant tubular neighbourhood of the form $G \times_{G_x} S$, where $S$ is a vector space and $G_x$ acts linearly on $S$.
When $G_x$ is finite (as is the case when $G$ acts almost freely on $M$), this essentially comes from the following lemma \cite[Lemma 1]{hC53}.

\begin{lem}\label{orthogonalisation}\textbf{(Orthogonalisation lemma.)}
  Let $M$ be a smooth manifold, let $G$ be a finite group acting smoothly on $M$, and suppose $x \in M$ is fixed by $G$.
  Then there exist coordinates centred at $x$ upon which $G$ acts orthogonally.
\end{lem}

This can be proven by considering the differential action of $G$ on $T_x M$ and using normal coordinates about $x$.

\begin{cor}\label{nowheredense}
  Let $G$ be a compact group acting smoothly, effectively, and almost freely on a manifold $M$.
  Then the fixed point set $M^G$ is nowhere dense in $M$.
\end{cor}

\begin{proof}
  Let $x$ be a point of $M$ and consider the isotropy group $G_x$.
  By the orthogonalisation lemma, there exists a coordinate neighbourhood $N$ centred at $x$ upon which $G_x$ acts linearly.
  If any $g \in G_x$ fixes an open neighbourhood of $x$, then since the only linear map which fixes an open set is the identity and since $G_x$ acts effectively on $N$, $g$ must act as the identity on $N$.
  It follows that any $g \in G$ which fixes an open set $U \subset M$ must be the identity map, and thus the fixed point set $M^G$ contains no open sets.
  Therefore $M^G$ is nowhere dense in $M$.
\end{proof}


\chapter{Fibred products}

Let $\cat{C}$ be a category, and consider a diagram in $\cat{C}$ of the form
\begin{equation}\label{prefp}
  \xymatrix{
    & Y \ar[d]^\gy \\
    X \ar[r]_\gf & Z.
  }
\end{equation}
Recall that the \emph{fibred product} of $X$ and $Y$ over $Z$, denoted $\fp{X}{\gf}{\gy}{Y}$, is an object in $\cat{C}$ along with two morphisms $p_1$, $p_2$ such that the diagram
\begin{equation}\label{fp}
  \xymatrix{
    \fp{X}{\gf}{\gy}{Y} \ar[d]_{p_1} \ar[r]^(.66){p_2} & Y \ar[d]^\gy \\
    X \ar[r]_\gf & Z
  }
\end{equation}
commutes, and such that the triple $(\fp{X}{\gf}{\gy}{Y},p_1,p_2)$ is universal with respect to this property.
It follows that if such a fibred product exists, it is unique up to isomorphism in $\cat{C}$.

Fibred products always exist in the category of sets and functions: in the situation of diagram \ref{prefp}, we can define $\fp{X}{\gf}{\gy}{Y}$ to be the set
\begin{equation*}
  \fp{X}{\gf}{\gy}{Y} := \{(x,y) \in X \times Y \mid \gf(x) = \gy(y)\} \subset X \times Y.
\end{equation*}
In this case $p_1$ and $p_2$ are just the restrictions of the coordinate projections $\pr_1$ and $\pr_2$.
If $A$ is a set and $\map{f}{A}{X}$, $\map{g}{A}{Y}$ are functions such that the diagram
\begin{equation*}\label{fp2}
  \xymatrix{
    A \ar[r]^g \ar[d]_f & Y \ar[d]^\gy \\
    X \ar[r]_\gf & Z
  }
\end{equation*}
commutes, then $\gs := (f,g)$ is the unique morphism making the diagram
\begin{equation}\label{fpx}
  \xymatrix@R=8pt@C=10pt{
    & A \ar@/_/[dddl]_f \ar[dd]^\gs \ar@/^/[dddr]^{g} & \\ \\
    & \fp{X}{f}{g}{Y} \ar[dl]^{p_1} \ar[dr]_{p_2} & \\
    X \ar[dr]_\gf & & Y \ar[dl]^\gy \\
    & Z & \\
  }
\end{equation}
commute.

Assuming now that these sets and functions are topological spaces and continuous maps, we can give $\fp{X}{\gf}{\gy}{Y}$ the subspace topology coming from the product topology on $X \times Y$.
This makes $\fp{X}{\gf}{\gy}{Y}$ is a fibred product in the category of topological spaces and continuous maps, since $(f,g)$ is a continuous map by the definition of the subspace and product topologies.

The situation is more subtle in other categories.
For example, fibred products do not always exist in the category of smooth manifolds and smooth maps, the obstruction being that the topological space $\fp{X}{f}{g}{Y}$ is not always a manifold.
We have the following criterion for a fibred product to exist.

\begin{lem}\label{fibredsub}
  Let $X$, $Y$, and $Z$ be smooth manifolds, and let $\map{\gf}{X}{Z}$ and $\map{\gy}{Y}{Z}$ be smooth maps, either of which is a submersion.
  Then the topological fibred product $\fp{X}{\gf}{\gy}{Y}$ is a manifold with a unique smooth structure.
\end{lem}

This lemma is used repeatedly throughout Chapter \ref{gpoids}; see \cite[\textsection II.2]{sL99} for a proof.


\footnotesize
\bibliographystyle{amsplain}
\bibliography{orbifold}

\end{document}